\documentstyle[12pt,amsfonts,amssymb,amsmath,amscd]{article}
\begin{document}
\title{18 LECTURES ON K-THEORY}
\author{Ioannis P. ZOIS
\\
Centre for Research, Templates and Testing
\\
Public Electricity Corporation
\\
9, Leontariou Street, GR-153 51
\\
Pallini, Athens, Greece
\\
and
\\ 
School of Natural Sciences, Department of Mathematics
\\
The American College of Greece, Deree College
\\
 6, Gravias Street, GR-153 42
\\
Aghia Paraskevi, Athens, Greece
\\
e-mail: i.zois@exeter.oxon.org}

\newtheorem{thm}{Theorem}
\newtheorem{defn}{Definition}
\newtheorem{prop}{Proposition}
\newtheorem{lem}{Lemma}
\newtheorem{cor}{Corollary}
\newtheorem{rem}{Remark}
\newtheorem{ex}{Example}

\newcommand{\Rat}{\mathbb Q}
\newcommand{\Real}{\mathbb R}
\newcommand{\RR}{\Real}
\newcommand{\Rh}{\hat{\Real}}
\newcommand{\Nat}{\mathbb N}
\newcommand{\Complex}{\mathbb C}
\newcommand{\HH}{\mathbb H_3}
\newcommand{\CC}{\Complex}
\newcommand{\Z}{\mathbb Z}

\newcommand{\Ea}{{\mathcal E}}
\newcommand{\Ta}{{\mathcal A}}
\newcommand{\Aa}{{\Ta_\infty}}
\newcommand{\Eb}{{C^*(D_1,D_2,X_2,\Omega)}}
\newcommand{\Tb}{{C^*(D_1,D_2,\hat\Omega)}}
\newcommand{\Ab}{{C^*(D_1,D_2,\Omega)}}
\newcommand{\aco}{\idl}
\newcommand{\aca}{\beta^{\|}}
\newcommand{\acb}{\beta^\perp}
\newcommand{\acc}{{\hat\tau}}

\newcommand{\Uu}{{\cal U}}
\newcommand{\Dd}{{\mathcal D}}
\newcommand{\Oo}{{\mathcal O}}
\renewcommand{\H}{{\mathcal H}}
\newcommand{\NN}{{\bf N}}
\newcommand{\ZZ}{{\bf Z}}
\newcommand{\Pp}{{\mathcal P}}
\newcommand{\Zz}{{\mathcal Z}}
\newcommand{\PP}{{\bf P}}
\newcommand{\LL}{\Lambda}
\newcommand{\LLL}{\Lambda\cup \infty}
\newcommand{\EE}{{\bf E}}
\newcommand{\Bb}{{\mathcal B}}
\newcommand{\Ww}{{\mathcal W}}
\newcommand{\Ss}{{\mathcal S}}

%Spuren
\newcommand{\Tr}{\mbox{\rm Tr}}  %Operatorspur
\newcommand{\TV}{{\mathcal T}}
\newcommand{\TVh}{\hat{{\mathcal T}}}
\newcommand{\tr}{\mbox{tr}}

\newcommand{\Rr}{{\mathcal R}}
\newcommand{\Nn}{{\mathcal N}}
\newcommand{\Cc}{{\mathcal C}}
\newcommand{\Jj}{{\mathcal J}}
\newcommand{\Ff}{{\mathcal F}}
\newcommand{\Ll}{{\mathcal L}}
\newcommand{\id}{{\mbox{\rm id}}}
\newcommand{\idl}{{\mbox{\rm\tiny id}}}
\newcommand{\eva}{{\mbox{\rm ev}}}
\newcommand{\eval}{{\mbox{\rm\tiny ev}}}

\def\essinf{\mathop{\rm ess\,inf}}
\def\esssup{\mathop{\rm ess\,sup}}

\newcommand{\bew}{{\bf Proof:}}
\newcommand{\eb}{\hfill $\Box$}

\newcommand{\x}{{\vec x}}
\newcommand{\y}{{\vec y}}
\newcommand{\n}{{\vec n}}
\renewcommand{\a}{{\vec a}}
\newcommand{\hull}{\Sigma}
\newcommand{\om}{\omega}
\newcommand{\oh}{{\hat{\om}}}

\newcommand{\Af}{{\Aa}_0}
\newcommand{\Tf}{{\Ta}_0}
\newcommand{\Ef}{{\Ea}_0}
\renewcommand{\H}{{\mathcal H}}

\newcommand{\bs}{\bigskip}
\newcommand{\ms}{\medskip}

\newcommand{\erz}[1]{\langle{#1}\rangle}
\newcommand{\pair}[2]{\erz{#1,#2}}
\newcommand{\diag}[2]{\mbox{\rm diag}(#1,#2)}
\newcommand{\tOmega}{{\Omega^s}}
\newcommand{\td}{{d^s}}
\newcommand{\tint}{{\int^s}}
\newcommand{\ttau}{{\tilde\tau}}
\newcommand{\talpha}{{\tilde\alpha}}
\newcommand{\tS}{{\tilde S}}
\newcommand{\ual}{{\underline{\alpha}}}
\newcommand{\hotimes}{{\hat\otimes}}
\newcommand{\K}{{\mathcal K}}
\newcommand{\Ypsilon}{{\Theta}}
\newcommand{\CA}{$C^*$-algebra}
\newcommand{\CF}{$C^*$-field}
\newcommand{\G}{\mathcal G}
\newcommand{\im}{\mbox{\rm im\,}}
\newcommand{\rk}{\mbox{\rm rk\,}}
\newcommand{\cp}{{\rtimes}}
\newcommand{\del}{{\bf \delta}}

\renewcommand{\a}{{\vec a}}
\newcommand{\at}{{\bf \tilde a}}
\newcommand{\cy}{\Phi}
\newcommand{\co}{\cy_\a}
\newcommand{\tco}{\cy_{\at}}

\newcommand{\ttimes}{{\tilde \rtimes}}
\renewcommand{\ss}{{\mathcal R}}
\renewcommand{\Cc}{{\mathcal C}}
\newcommand{\Ri}{(\RR\cup\infty)}
\newcommand{\Rd}{\hat\RR}
\newcommand{\dom}{\mbox{\rm dom}}
\newcommand{\fin}{\mbox{\rm\tiny fin}}

\newcommand{\enn}{\mu}
\newcommand{\ch}{\mbox{\rm ch}}
\newcommand{\supp}{\mbox{\rm supp}}
\newcommand{\chd}{\ch}
%%%%%%%%%%%%%%%%%%%%%%%%%%%%%
\newcommand{\Oh}{\hat\hull}

\newcommand{\hsp}{\RR^{d-1}\times\RR^{\leq 0}}
\newcommand{\IDS}{IDS}
\newcommand{\sv}[2]{\left(\begin{array}{c} #1 \\ #2 \end{array}\right)}

\bibliographystyle{amsalpha}

\maketitle

\newpage

\tableofcontents

\newpage

{\bf Introduction and background Bibliography}\\

Back in 1995 a graduate summer school on \emph{K-Theory} was organised at the University of Lancaster (UK) by the 
\textsl{London Mathematical Society (LMS)}. This summer school covered the three basic branches of K-Theory, namely 
\emph{Topological K-Theory, Analytic K-Theory (K-Homology) and Higher Algebraic K-Theory}. The main lecturers were 
Professors J.D.S. Jones (Warwick), J. Roe (Oxford then, now at Penn. State) and D.G. Quillen (Oxford) respectively.\\ 

The author, a graduate student at the Mathematical Institute of the University of Oxford at that time, had the 
chance to be one of the attendants of that summer school and his personal notes form the backbone of this book. 
Those notes were expanded and polished during the years since the author had the chance to teach a graduate course 
on K-Theory twice in the past at the Universities of Cardiff (UK) and Athens (Greece). This experience  gave the 
motivation to present the book as a group of 18 hourly lectures, 6 for each branch of K-Theory. Subsequently one 
further chapter was added, containing an introduction to \emph{Waldhausen K-Theory} (this last chapter was never taught in class 
though) along with 2 setions, one on twisted K-theory and applications in string/M-Theory and one on the so called gap labeling problem in solid state physics. For convenience we include an appendix which contains some 
results from other areas of mathematics (especially from algebraic topology). A few Propositions have no proofs, this is happens when proofs are either fairly 
straightforward (and in this case they are left as excercises for the reader) or too extensive and complicated and we only exhibit the main ideas and provide references.\\
 
The mathematical bibliography on K-Theory cannot be considered to be extensive. We give an almost complete 
list of relevant books 
available (we have made some use of these books in our notes):\\ 

For topological K-Theory one has the clasic 1967 Harvard notes by M.F. Atiyah (see M.F. Atiyah: \emph{"K-Theory"}, 
Benjamin 1967) 
and a recent book by Efton Park: \emph{"Complex Topological K-Theory"}, Cambridge University Press 2008.
There are also some online notes by A. Hatcher from Cornell University (see http://www.math.cornell.edu/~hatcher/) 
and some online notes by Max Karoubi in Paris 
(see http://people.math.jussieu.fr/~karoubi/KBook.html) along with his book \emph{"K-Theory - An Introduction"}, Springer (1978).\\ 
 
For analytic K-Theory there is the book by N.E. Wegge-Olsen: \emph{"K-Theory and $C^*$ Algebras"}, Oxford University 
Press 1993, the more advanced book by B.Blackadar: \emph{"K-Theory for Operator Algebras"}, Cambridge University Press 
1998, the book by Nigel Higson and John Roe: \emph{"Analytic 
K-Homology"}, Oxford University Press 2000 and the book by F. Larsen,  M. Rrdam and M. Rordam: 
\emph{"An Introduction to K-Theory for $C^*$-algebras"}, Cambridge University Press 2000.
There are also some nice brief online notes by T. Gowers (from Cambridge University): \emph{"K-Theory of Banach Algebras"}  
(see http://www.dpmms.cam.ac.uk/~wtg10/).\\ 

For the Higher Algebraic K-Theory there are the books by Jonathan Rosenberg: \emph{"Algebraic K-Theory and its 
Applications"}, Springer Graduate Texts in Mathematics (1994), the book by V. Srinivas: \emph{"Algebraic K-Theory"}, Birkhauser, Boston (1996) and the book by Hvedri Inassaridze: \emph{"Algebraic K-Theory"}, Kluwer 1995. There are also some on line notes by C. Weibel: \emph{"The K-book: An introduction to Algebraic K-Theory"}.\\ 

\textsl{Thus there is not a concise introduction on K-Theory available in the bibliography covering all basic three branches}.\\

The "fathers" of K-Theory are M.F. Atiyah and A. Grothendieck, arguably the greatest mathematicians of the second 
half of the 20th century. The name was given by Grothendieck, the letter \emph{"K"} stands for the German word 
\emph{"(die) Klasse"} which means \textsl{class} in English. K-Theory is one of the so-called \textsl{generalised (or exotic) 
homology theories} and satisfies 4 out of the 5 Eilenberg-Steenroad homology axioms: it satisfies homotopy, excision, 
additivity and exactness but it does not satisfy the dimension axiom, namely the K-Theory of a point is $\mathbb{Z}$ and 
not zero. During the 50 years since its birth, K-Theory has 
been proven a very useful tool in many areas of mathematics including topology, global analysis, index theory, 
number theory, algebra, (noncommutative) geometry etc. Moreover K-Theory has found many applications in theoretical 
physics as well, for example one can mention the anomaly cancellation in quantum field theories, the so called gap labelling problem in solid state physics, the topological charges of membranes in M-Theory  (where topological charges of membranes are considered to be classes of twisted K-Theory) etc.\\

Mathematical evolutions of K-Theory are the \emph{equivariant K-Theory} of M.F. Atiyah and G.Segal, the \emph{L-Theory} 
used in surgery of manifolds, the \emph{KK-Theory} (Kasparov K-Theory or bivariant K-Theory), the 
\emph{E-Theory} of A. Connes, the \emph{Waldhausen K-Theory} or \emph{"A-Theory"} (which is a topological version of Quillen's Higher 
Algebraic K-Theory) etc. We should also mention the close relation between Higher Algebraic K-Theory and the theory 
of Motives (\emph{motivic cohomology}) by V. Voevodsky.\\

This book is suitable for graduate students, hence we assume that the reader has a good knowledge of algebra 
(in particular groups, modules and homological algebra), geometry (manifold theory) and algebraic topology 
(singular and chain (co)homology for CW-complexes and basic homotopy theory). The end of a proof is denoted by a white box. In the notation used, for example, Proposition 2.3.9 refers to proposition 9 in section (lecture) 3 in chapter 2. A list of suitable 
recommended books with background material follows:\\

{\bf Algebra}\\
$\bullet$ S. MacLane, G. Birckoff: \emph{"Algebra"}, Chelsea, 1988\\ 
$\bullet$ S. Lang: \emph{"Algebra"}, Addison Wesley 1993\\ 
$\bullet$ P. J. Hilton, U. Stammbach: \emph{"A Course in Homological Algebra"}, Springer 1997.\\
$\bullet$ N. Bourbaki, \emph{"Algebra"}, Vol I,II\\
$\bullet$ W.A. Adkins, S.H. Weintraub, \emph{"Algebra"}, Springer, 1992\\
$\bullet$ P.M. Cohn \emph{"Algebra"}, John Wiley, 1989\\
$\bullet$ T.H. Hungerford, \emph{"Algebra"}, Springer, 1980\\

{\bf Geometry}\\
$\bullet$ N. Hitcin online notes on manifolds, see\\
http://people.maths.ox.ac.uk/~hitchin/hitchinnotes/hitchinnotes.html\\
$\bullet$ P. Griffiths and J. Harris: \emph{"Principles of Algebraic Geometry"}, Wiley, 1994\\
$\bullet$ S. Kobayashi, K. Nomizu: \emph{"Foundations of Differential Geometry"}, Vol I,II, Wiley, 1996\\
$\bullet$ J.M. Lee, \emph{"Introduction to Smooth Manifolds"}, Springer, 2002\\

{\bf Topology}\\
$\bullet$ A. Hatcher, on line notes on algebraic topology,\\
see http://www.math.cornell.edu/~hatcher/AT/AT.pdf\\
$\bullet$ J. P. May: \emph{"A concise Course in Algebraic Topology"}\\ see http://www.math.uchicago.edu/~may/CONCISE/ConciseRevised.pdf\\ 
$\bullet$ E.H. Spanier: \emph{"Algebraic Topology"}, Springer 1991\\ 
$\bullet$ A.Dold: \emph{"Lectures on Algebraic 
Topology"}, Springer 1991\\ 
$\bullet$ C.R.F. Maunder: \emph{"Algebraic Topology"}, Dover 1996\\ 
$\bullet$ W.S. Massey: \emph{"A Basic Course in 
Algebraic Topology"}, Springer 1991\\ 
$\bullet$ J. Munkres: \emph{"Elements of Algebraic Topology"}, Westview 1995\\ 
$\bullet$ R. Bott, L. Tu: \emph{"Differential forms in algebraic topology"}, Springer 1982\\
$\bullet$ N. Steenrod: \emph{"The topology of fibre bundles"}, Princeton (1951)\\
$\bullet$ S. Eilenberg and N. Steenrod: \emph{"Foundations of Algebraic Topology"}, Princeton (1952)\\

IPZ, Athens, August 2010\\
http://sites.google.com/site/ipzoisscience/home\\

\newpage

\section{Topological K-Theory}

Topological K-Theory is 
historically the first branch of K-Theory which was developed by Atiyah and Hirzebruch in 1960's. Perhaps the two most famous applications of topological K-Theory is the result of Adams on the maximum number of continuous linearly independent tangent vector fields on spheres and the fact that there are no finite dimensional division algebras over $\mathbb{R}$ in dimensions other than $1,2,4$ and $8$, corresponding to the reals, complex numbers, the quaternions and the  octonions. The basic reference for this chapter are the clasic Harvard notes by M.F. Atiyah from the 1960's. Some other 
books and notes which have appeared subsequently in the literature are based on Atiyah's notes. The key notion is the notion of a vector bundle.\\

\section{Lecture 1 (Topological Preliminaries on Vector Bundles)}

Everyone is familiar with the annulus and the Mobious band. These spaces can be constructed from the circle and the line, the first is their actual Cartesian product whereas the second is their "twisted" product. We start with the fundamental definition of a vector bundle which generalises these two constructions where instead of the circle we take any topological space and instead of the line we use some vector space of finite dimension. We use the symbol ${\mathbb{F}}$ to denote either the field of real ${\mathbb{R}}$ or the complex numbers ${\mathbb{C}}$.\\

{\bf Definition 1.} A \emph{vector bundle} $\xi =(E, \pi , B)$ consists of the following data:\\ 
{\bf 1.} A topological space $E$ (also denoted $E (\xi )$, often assumed to be a manifold) which is called the \emph{total space}.\\
{\bf 2.} A topological space $B$ (often assumed to be a manifold) which is called the \emph{base space} or simply \textsl{base}.\\
{\bf 3.} A continuous surjective map $\pi :E\rightarrow B$ which is called the \emph{projection} satisfying the following properties:\\

{\bf i.} For all $b\in B$, the space $E_b (\xi )=F_b(\xi )=\pi ^{-1}(b)\subset E$ which is called the \emph{fibre} over $b$ has the structure of a \textsl{vector space} of finite dimension say $n$. In this case we say that the vector bundle has \textsl{dimension (or rank)} $n$ (a word of caution here, this dimension may be different from the topological dimension of $E$).\\
{\bf ii.} For every  $b\in B$, the \textsl{local triviality condition} is satisfied, namely if $U_b\subset B$ is a neighbourhood of $b$, then  there is a map
$$h_b:U_{b}\times{\mathbb{F}}^n\rightarrow \pi ^{-1}(U_b)\subset E$$
which is a homeomorphism. These maps like $h_b$ above are called \emph{local trivialisations} and they have the property 
that when restricted to each fibre they give linear isomorphisms of  vector spaces , namely the maps 
$${\mathbb{F}}^n\rightarrow F_b(\xi)$$ 
with
$$v\mapsto h_b(b,v),$$
are linear isomorphisms of vector spaces.\\

Choosing the reals or the complex numbers we get real or complex vector bundles respectively.\\ 

Throughout this chapter, for simplicity vector bundles will be denoted either by small greek letters or by writing only the projection map or even by simply writing the total space (denoted using capital latin letters) when there are no ambiguities.\\

Vector bundle maps can be defined in the obvious way: If $\pi _1:E_1\rightarrow B$ and $\pi _2:E_2\rightarrow B$ are two vector bundles over the same base $B$, then a \emph{vector bundle map} $f:E_1\rightarrow E_2$ is a homeomorphism such that $\pi _2 f=\pi _1$ and when restricted to each fibre it gives an ${\mathbb{F}}$-homomorphism of vector spaces. It is clear that ${\mathbb{F}}$ vector bundles over $B$ and their maps form an additive category. We are more interested in a special kind of vector bundle maps called isomorphisms and the corresponding vector bundles then will be called isomorphic:\\ 

{\bf Definition 2.} Two vector bundles $\xi$, $\eta$ over the same base $B$ will be called \emph{isomorphic} and it will be denoted $\xi\simeq\eta$, if there exists an isomorphism between them, namely a  homeomorphism say $\phi :E(\xi )\rightarrow E(\eta )$ between their total spaces 
which when restricted to each fibre $\phi _b: F_b (\xi )\rightarrow F_b (\eta )$ yields a \textsl{linear isomorphism between 
vector spaces} where $F_b (\xi )$ is an alternative notation for the fibre over the point $b$ of the bundle $\xi$. Obviously isomorphic vector bundles must have the same 
dimension (yet the converse is not true).\\

{\bf Remark 1.} Given any vector bundle $\pi :E\rightarrow B$, there is an alternative way to reconstruct it using the \emph{gluing functions}:  We take an open cover $\{U_a\}$ of $B$ with local trivialisations $h_a :\pi ^{-1}(U_a)\rightarrow U_a \times {\mathbb{F}}^{n}$ and then we reconstruct $E$ as the quotient space of the disjoint union 
$$\bigcup _{a}(U_a\times {\mathbb{F}}^{n})$$
obtained by identifying $(x,v)\in U_a\times {\mathbb{F}}^{n}$ with 
$h_{b}h_{a}^{-1}(x,v)\in U_{b}\times {\mathbb{F}}^{n}$ whenever $x\in U_a\cap U_b$. The functions $h_{b}h_{a}^{-1}$ can be viewed as maps $g_{ab}:U_a\cap U_b \rightarrow GL_n({\mathbb{F}})$ which satisfy the \textsl{cocycle condition} 
$$g_{cb}g_{ba}=g_{ca}$$
on $U_a\cap U_b\cap U_c$. Any collection of gluing functions satisfying the cocycle condition can be used to construct a vector bundle $\pi :E\rightarrow B$.\\

{\bf Remark 2.} One can generalise the above definition of vector bundles by assuming that the model fibres over the points of the base space are in general homeomorphic to some fixed topological space $F$ instead of an $n$-dim vector space thus obtaining the definition of an arbitrary \emph{fibre bundle with fibre} $F$. In this way, fibre bundles provide a generalisation of the Cartesian product.\\

{\bf Remark 3.} Some authors relax the condition that all fibres in a vector bundle should have the same dimension; however by continuity of the local trivialisations, the dimension of the fibres must be locally constant. If the base space is connected (which will be always the case here), then the dimension of the fibres will be constant (and thus we can define the dimension or the rank of a vector bundle).\\

{\bf Examples:}\\

{\bf 1.} Given any topological space $B$, the \emph{trivial vector bundle} over $B$ of dimension $n$ is the vector bundle with total space $B\times{\mathbb{F}}^n$, base space $B$ and the projection $\pi$ is the projection to the first factor
$$\pi :B\times{\mathbb{F}}^n \rightarrow B.$$
For the trivial vector bundle over $B$ of dimension $n$ we may use the alternative notation ${\mathbb{F}}^{n}_{B}$. (Sometimes we may omit the base space alltogether to simplify our notation if no confusion arises).\\  

From this example one can see that an arbitrary vector bundle need not be globally trivial but the triviality condition states that all vector bundles "look like" the product bundle \emph{locally}. In other words one can say that a vector bundle is a continuous family of vector spaces over some (base) space.\\

{\bf 2.} Let $I=[0,1]$ denote the unit interval and let $E$ be the quotient space of $I\times {\mathbb{R}}$ under the identifications $(0,t)\sim (1,-t)$. Then the projection $I\times {\mathbb{R}}\rightarrow I$ induces a map $\pi :E\rightarrow S^1$ which is a $1$-dim vector bundle (these are called in particular \textsl{line bundles}). Since $E$ is homeomorphic to the Mobius band with its boundary circle deleted, this is called the \textsl{Mobius bundle}.\\  

{\bf 3.} Let $M$ be a real differentiable manifold of dimension say $n$. Then its \emph{tangent bundle} 
$TM=\{(x,v)\in M\times{\mathbb{R}}^{n}:v\in T_{x}M\}$, where $T_{x}M$ denotes the tangent space at the point $x$, is a real vector bundle of dimension $n$ with base space $M$ and fibre isomorphic to ${\mathbb{R}}^n$. The construction applies to complex manifolds as well.\\

{\bf 4.} The real projective $n$-space ${\mathbb{R}}P^{n}$ is by definition the space of lines in ${\mathbb{R}}^{n+1}$ passing through the origin. Then ${\mathbb{R}}P^{n}$ can be regarded as the quotient space of $S^n$ (the $n$-sphere) with the antipodal pairs of points identified. The \emph{tautological} (or "canonical") line bundle over ${\mathbb{R}}P^{n}$ has as total space the subset $E$ of ${\mathbb{R}}P^{n}\times {\mathbb{R}}^{n+1}$ defined by
$$E=\{(x,v)\in {\mathbb{R}}P^{n}\times {\mathbb{R}}^{n+1}:v\in x\}.$$    
Then one has an obvious projection map
$$\pi :E\rightarrow {\mathbb{R}}P^{n}$$ 
with $\pi (x,v)=x$ (projection to the first factor), namely each fibre of $\pi$ is the \emph{line} $x$ inside the Euclidean $(n+1)$-space ${\mathbb{R}}^{n+1}$ (thus the dimension of this tautological vector bundle is $1$ and hence it is a line bundle).\\

The tautological line bundle over ${\mathbb{R}}P^{n}$ has also an \emph{orthogonal complement} vector bundle with total space $E^{\perp}$ where
$$E^{\perp}=\{(x,v)\in {\mathbb{R}}P^{n}\times {\mathbb{R}}^{n+1}:v\perp x\}.$$    
One has the obvious projection to the first factor. This vector bundle has dimension $n$.\\
Both these constructions can be applied to the complex case as well.\\
 
{\bf 5.} A natural generalisation of the real projective space is the real Grassmannian $G_{k}({\mathbb{R}}^n)$ which is defined as the space of $k$-dim planes through the origin of ${\mathbb{R}}^n$ (obviously $k<n$). In a similar fashion one can define the \emph{canonical} $k$-dim vector bundle over $G_{k}({\mathbb{R}}^n)$ consisting of pairs $(x,v)$ where $x$ is a "point" in the Grassmannian (a $k$-dim subspace) and $v$ is a vector in $x$.
This has an \emph{orthogonal complement} too which is an $(n-k)$-dim vector bundle over the Grassmannian $G_{k}({\mathbb{R}}^n)$. There is also the complex version of them.\\

{\bf Definition 3.} Given some vector bundle $\pi :E\rightarrow B$, we take a subspace $A\subset B$. Then $\pi :\pi ^{-1}(A)\rightarrow A$ is clearly a vector bundle called the \emph{restriction} of $E$ over $A$ and it will be denoted $E|_{A}$.\\

{\bf Definition 4.} A \emph{vector subbundle} of some vector bundle $\pi :E\rightarrow B$ is a subspace $E_0\subset E$ which intersects every fibre in a vector subspace so that $\pi :E_0\rightarrow B$ is again a vector bundle.\\
   
{\bf Definition 5.} Given two vector bundles $\pi _{1} :E_1\rightarrow B_1$ and  $\pi _{2} :E_2\rightarrow B_2$, the \emph{product vector bundle} is $\pi _1 \times \pi _2:E_1\times E_2\rightarrow B_1\times B_2$ where the fibres are the Cartesian products of the form $\pi _1^{-1}(b_1)\times\pi _2^{-1}(b_2)$ and  $h_a\times h_b$ are the local trivialisations where $h_a:U_{a}\times{\mathbb{F}}^n\rightarrow \pi _{1}^{-1}(U_a)$ and $h_b:U_{b}\times{\mathbb{F}}^m\rightarrow \pi _{2}^{-1}(U_b)$ are local trivialisations for $E_1$ and $E_2$ respectively.\\

{\bf Definition 6.} Given any vector bundle $(E,\pi ,B)$, a \emph{section} on $E$ is a continuous map $s:B\rightarrow E$ such that for all $b\in B$ one has $\pi (s(b))=b$.\\

Every vector bundle has a canonical section, the \emph{zero section} whose value on each fibre is zero. One often identifies the zero section with its image, a subspace of $E$ which projects homeomorphically onto $B$ by $\pi$. One can sometimes distinguish non-isomorphic vector bundles by looking at the complement of the zero section since any vector bundle isomorphism $h:E_1\rightarrow E_2$ must take the zero section of $E_1$ to the zero section of $E_2$, hence their complements must be homeomorphic. At the other extreme
 from the zero section, one has a \emph{nowhere vanishing section}. Isomorphisms also take nowhere vanishing sections to nowhere vanishing sections. Clearly the trivial bundle has such a section yet not all vector bundles have one. For instance we know that the tangent bundle of the sphere $S^n$ has a nowhere vanishing section if and only if $n$ is odd. From this it follows that the tangent bundle of $S^n$ for $n$ even is not trivial (for the proof see for example \cite{hatcher}). \\

In fact an $n$-dim vector bundle $\pi :E\rightarrow B$ is isomorphic to the trivial bundle if and only if it has $n$ sections $s_1,...,s_n$ such that the vectors $s_1(b),...,s_n(b)$ are linearly independent in each fibre $\pi ^{-1}(b)$. One direction is evident since the trivial bundle certainly has sections and isomorphisms take linearly independent sections to linearly independent ones. Conversely, if one has $n$ linearly independent sections $s_i$, the map $B\times {\mathbb{F}}^{n}\rightarrow E$ given by $h(b,t_1,t_2,...,t_n)=\sum _{i}t_is_i(b)$ is a linear isomorphism in each fibre; moreover it is continuous since its composition with a local trivialisation $\pi ^{-1}(U)\rightarrow U\times {\mathbb{F}}^n$ is continuous. Hence $h$ is an isomorphism since it maps fibres isomorphically to fibres.\\
  
It follows immediately from the definition of a vector bundle that for any $b\in B$ there exists a neighborhood $U$ of $b$ and sections $s_1,...,s_n$ of $E$ over $U$ such that $s_1(x),...,s_n(x)$ form a basis for the fibre $F_x=\pi ^{-1}(x)$ over $x\in U$ for all $x\in U$. We say that $s_1,...,s_n$ form a \emph{local basis} at $b$ and any section of $E$ can be written as $s(x)=\sum_{i}a_i(x)s_i(x)$ where $a_i(x)\in{\mathbb{F}}$. Clearly $s$ is continuous of all the $a_i$ functions are.\\

Let us consider some fixed topological space $X$ as base space; we would like to study the set $VB(X)$ of all 
(real or complex) vector bundles over $X$ of finite dimension. In this set we can define the following operations:\\

\emph{1. Direct sum}: If $\xi ,\eta$ are two vector bundles over $X$, we can define their direct sum $\xi\oplus\eta$ 
which is a new vector bundle over $X$ with total space
$$E(\xi\oplus\eta )\subset E(\xi )\times E(\eta )$$
where
$$E(\xi\oplus\eta )=\{(v,w):\pi _{\xi}(v)=\pi _{\eta}(w)\}.$$
The dimension of this new vector bundle is equal to the sum 
of the dimensions of
$\xi ,\eta$ and every fibre is the direct sum of the corresponding fibres (recall that the fibres are vector spaces and hence 
one can form their direct sum) 
$$F_{x}(\xi\oplus\eta )=F_{x}(\xi )\oplus F_x(\eta ).$$

\emph{2. Pull-back}: If $f:X\rightarrow Y$ is a continuous map and $\xi\rightarrow Y$ a vector bundle over $Y$, then we can 
define the vector bundle $f^*(\xi )\rightarrow X$ which is called the pull-back of $\xi$ over $X$ by $f$ 
which is a new vector bundle over $X$ with total space 
$$E(f^{*}(\xi ))\subset X\times E(\xi )$$
with
$$E(f^{*}(\xi ))=\{(x,e):\pi _{\xi}(e)=f(x)\}.$$
Its fibre is defined as follows:
$$F_x(f^*(\xi ))=F_{f(x)}(\xi ).$$

\emph{3. Tensor product}: Let $\pi _1:E_1\rightarrow X$ and  $\pi _2:E_2\rightarrow X$ be two vector bundles over the same base $X$. Then we can form their \emph{tensor product} $E_1\otimes E_2$ which is a vector bundle over $X$. As a set, the total space $E_1\otimes E_2$ is the disjoint union of the tensor product of the vector spaces (fibres) $\pi _{1}^{-1}(x)\otimes\pi _{2}^{-1}(x)$. The structure of a vector bundle over $X$ can be given using the gluing functions: Let $\{U_a\}$ be an open cover of $X$ so that both $E_i$, where $i=1,2$ are trivial over each $U_a$ and thus one can obtain gluing functions $g_{ba}^{i}:U_a\cap U_b \rightarrow GL_{n_i}({\mathbb{F}})$ for each $E_i$. Then the gluing functions for the tensor bundle $E_1\otimes E_2$ are obtained via the tensor product functions $g_{ba}^{1}\otimes g_{ba}^{2}$ which assign to each $x\in U_a\cap U_b$ the tensor product of the two matrices $g_{ba}^{1}(x)$ and $g_{ba}^{2}(x)$.\\

\emph{4. Quotient bundles}: Given a vector bundle $\pi :E\rightarrow X$ and a vector subbundle $E_0\subset E$ of $E$ (which is another vector bundle over $X$), we can form the \textsl{quotient bundle} $E/E_0\rightarrow X$ where the total space has the quotient topology and each fibre over an arbitrary point $x\in X$ is the quotient of the corresponding fibres $F_x(E)/F_x(E_0)$; clearly the dimension of the quotient bundle $E/E_0$ equals the difference of the dimensions of $E$ and $E_0$, namely
$$dim(E/E_0)=dimE-dimE_0.$$ 

{\bf Exercise 1.} \textsl{Define the exterior products of vector bundles}.\\

{\bf Exercise 2.} \textsl{Prove that the pull-back respects both the direct sum and the tensor product, namely} $f^{*}(E_1\oplus E_2)=f^{*}(E_1)\oplus f^{*}(E_2)$ and $f^{*}(E_1\otimes E_2)=f^{*}(E_1)\otimes f^{*}(E_2)$.\\ 

It is obvious that the direct sum of trivial vector bundles will be another trivial vector bundle; a nontrivial vector bundle which becomes trivial after taking the direct sum with a trivial vector bundle will be called \emph{stably trivial}; what is not obvious is that the direct sum of nontrivial vector bundles can be trivial; under certain assumptions, the later is always the case:\\

{\bf Proposition 1.} For any vector bundle $\pi :E\rightarrow B$ where $B$ is compact and Hausdorff, there exists another vector bundle $E'\rightarrow B$ such that $E\oplus E'$ is trivial.\\

For the proof of this proposition we need a definition and a Lemma:\\

{\bf Definition 7.} An \emph{inner product} on a vector bundle $\pi :E\rightarrow B$ is a map $<,>:E\oplus E\rightarrow {\mathbb{R}}$ which restricts on each fibre to an inner product (i.e. a positive definite symmetric bilinear form).\\ 

[\emph{Note:} We treat the real case for brevity here; it is straightforward to get the compelx case as well].\\

Inner products always exist when the base space is compact and Hausdorff (or more generally paracompact, namely a topological space which is Hausdorff and every open cover has a subordinate partition of unity which is constructed using Urysohn's Lemma). Inner products can be obtained as pull-backs of the standard Euclidean inner product in ${\mathbb{R}}^n$ by the local trivialisations and then extend it to the whole of $B$ using a partition of unity.\\

{\bf Lemma 1.} If $\pi :E\rightarrow B$ is a vector bundle over a paracompact space $B$ and $E_0\subset E$ is a vector subbundle, then there is another vector subbundle $E_{0}^{\perp}\subset E$ such that $E_{0}\oplus E_{0}^{\perp}\simeq E$.\\

{\bf Proof of Lemma 1:} Suppose $E$ has dim $n$ and $E_0$ has dim $m<n$. We choose an inner product on $E$ and we define  
$E_{0}^{\perp}$ to be the subspace of $E$ whose fibres consist of vectors orthogonal to the vectors of $E_0$. We claim that the natural projection $E_{0}^{\perp}\rightarrow B$ defines a vector bundle over $B$ (thus making $E_{0}^{\perp}$ a subbundle and not just a subspace of $E$). If this is so, then we have our result, namely $E\simeq E_{0}\oplus E_{0}^{\perp}$ the isomorphism being $(v,w)\mapsto v+w$ (since isomorhisms of vector bundles restrict to linear isomorphisms of vector spaces on each fibre).

We have to check the local triviality condition on $E_{0}^{\perp}$. Since $E_0$ has dimension $m$, it has $m$ linearly independent local sections $b\mapsto (b,s_i(b))$ near each point $b_0\in B$. We extend this set of $m$ linearly independent local sections of $E_0$ to a set of $n$ linearly independent local sections of $E$ by choosing $s_{m+1},...,s_n$ first in the fibre $\pi ^{-1}(b_0)$, then taking the same vectors for all nearby fibres since if $s_1,...,s_m,...,s_n$ are linearly independent at $b_0$, they will remain independent for nearby $b$'s by the continuity of the determinant. Apply the Gram-Schmidt orthogonalisation process to $s_1,...,s_m,...,s_n$ using the inner product in each fibre to obtain new sections $s'_{i}$. The explicit formulae for the orthogonalisation process show that the new local sections are continuous and the first $m$ of them give a basis of $E_0$ on each fibre. The new sections $s'_{i}$ allow us to define a local trivialisation $h:\pi ^{-1}(U)\rightarrow U\times{\mathbb{R}}^n$ with $h(b,s'_{i}(b))$ equal to the ith standard basis vector of ${\mathbb{R}}^n$. This $h$ carries $E_0$ to $U\times{\mathbb{R}}^m$ and $E_{0}^{\perp}$ to $U\times{\mathbb{R}}^{n-m}$ and hence $h|_{E_{0}^{\perp}}$ is a local trivialisation for $E_{0}^{\perp}$. $\square$\\

{\bf Proof of Proposition 1:} In order to give the motivation of the construction, let us assume that the result holds and hence $E$ is a subbundle of the trivial bundle $B\times{\mathbb{R}}^N$ for some large enough $N$, thus one has an inclusion  $E\hookrightarrow B\times{\mathbb{R}}^N$; composing this inclusion with the projection onto the second factor one obtains a map $E\rightarrow{\mathbb{R}}^N$ which is a linear injection on each fibre; the idea of the proof is to reverse the logic, namely first we construct a map $E\rightarrow{\mathbb{R}}^N$ which is a linear injection on each fibre and next we shall show that this gives an embedding of $E$ in $B\times{\mathbb{R}}^N$ as a direct summand.\\

Each point $x\in B$ has a neighborhood $U_x$ over which $E$ is trivial. By Urysohn's Lemma there is a map $\phi _x:B\rightarrow [0,1]$ that is $0$ outside $U_x$ and nonzero at $x$. By letting $x$ vary, the sets $\phi ^{-1}_x (0,1]$ form an open cover of $B$ and since $B$ is compact, this cover has a finite subcover. We relabel the corresponding $U_{x}'s$ and $\phi _{x}'s$ $U_j$ and $\phi _j$ respectively. We define $g_j:E\rightarrow{\mathbb{R}}^n$ via $g_j(v)=\phi _{j}(\pi (v))[p_{j}h_{j}(v)]$ where $p_{j}h_{j}$ is the composition of a local trivialisation $h_j:\pi ^{-1}(U_j)\rightarrow U_j\times{\mathbb{R}}^n$ with the projection $p_j$ to ${\mathbb{R}}^n$. Then $g_j$ is a linear injection on each fibre over $\phi _{j}^{-1}(0,1]$; hence if we take the $g_{j}$'s as the coordinates  of a map $g:E\rightarrow{\mathbb{R}}^N$ with ${\mathbb{R}}^N$ a product of copies of ${\mathbb{R}}^n$, then $g$ is a linear injection on each fibre.\\

The map $g$ is the second coordinate of a map $f:E\rightarrow B\times{\mathbb{R}}^N$  with first coordinate $\pi$. The image of $f$ is a subbundle of the product $B\times{\mathbb{R}}^N$ since projection of       
${\mathbb{R}}^N$ onto the ith ${\mathbb{R}}^n$ factor gives the second coordinate of a local trivialisation over $\phi _{j}^{-1}(0,1]$. Thus $E$ is isomorphic to a subbundle of $B\times{\mathbb{R}}^N$ so by the Lemma 1 there exists a complementary subbundle $E'$ with $E\oplus E'$ isomorphic to  $B\times{\mathbb{R}}^N$. $\square$\\

\textsl{In the sequel, we shall assume that the topological space $X$ is compact and Hausdorff and our vector bundles are 
complex}.\\

Let $Vect_{k}(X)$ denote the set of \textsl{isomorphism classes of vector bundles over} 
$X$ of dimension $k$, let $Vect(X)$ denote the set of \textsl{isomorphism classes of vector bundles over} 
$X$ of finite dimension and let $VB(X)$ denote the set of \textsl{all vector bundles over} $X$ of finite dimension. Then the direct sum $\oplus$ defined above defines an addition on $Vect(X)$:
$$[\xi ]+[\eta ]=[\xi\oplus\eta ].$$
This operation is associative, commutative, there exists a neutral element (a zero) but one cannot define symmetric 
elements (the "opposite" isomorphism class of bundles), thus 
$(Vect(X),+)$ is a(n additive abelian) \emph{semi-group}.\\

There is a trick which can concoct a group out of a semi-group which goes back to Grothendieck: This is a certain symmetrisation process in the same way that we define ${\mathbb{Z}}$ from the (additive) semi-group ${\mathbb{N}}$. More concretely, to an arbitrary semi-group $M$ we associate a new group denoted $Gr(M)$ defined as the quotient of $M\times M$ by the equivalence relation
$$(a,b)\sim (c,d)\Leftrightarrow\exists e\in M : a+d+e=b+c+e.$$ 

There is an alternative technique which does the same job: If $M$ is an arbitrary semi-group, then 
we consider the diagonal homomorphism of semi-groups 
$$\Delta :M\rightarrow M\times M.$$
Then $Gr(M)$ (the Grothendieck or universal group of the semi-group $M$) is the set of cosets of $\Delta (M)$ in $M\times M$ which is the quotient semi-group; yet 
the interchange of factors in the Cartesian product $M\times M$ 
induces symmetric elements in $Gr(M)$ and thus it is promoted to a group.\\

Back to our K-Theory, we have thus the following key definition due to Alexander Grothendieck:\\

{\bf Definition 8.} We define the Abelian group $KX$ as the universal group (or the "Grothendieck group" denoted $Gr$) of the semi-group 
$Vect(X)$, namely
$$KX=Gr[Vect(X)].$$
(This is in fact the 0th K-group of the topological K-Theory as we shall see later). [Its customarily to denote $KX$ and $KOX$ and the complex and the real case respectively].\\

{\bf Definition 9.} Two vector bundles $E_1$ and $E_2$ over the same base space $X$ are called \emph{stably isomorphic} and it will be denoted $E_1\simeq _{S}E_2$ if there exists some $n\in{\mathbb{N}}$ so that $E_1\oplus {\mathbb{C}}^{n}_{X}\simeq E_2\oplus {\mathbb{C}}^{n}_{X}$.\\

{\bf Definition 10.} Two vector bundles $E_1$ and $E_2$ over the same base space $X$ are called \emph{similar} and it will be denoted $E_1\sim E_2$ if there exist some $n,m \in{\mathbb{N}}$ so that $E_1\oplus {\mathbb{C}}^{n}_{X}\simeq E_2\oplus {\mathbb{C}}^{m}_{X}$.\\

It is really straightforward to prove that both $\simeq _S$ and $\sim$ are equivalence relations in $VB(X)$. On equivalnce classes of both kinds the direct sum is well defined and moreover it is associative and commutative. A neutral (zero) lement exists which is the class of ${\mathbb{C}}^{0}_{X}$.\\ 

{\bf Definition 8'.} An equivalent definition for the Abelian group $KX$ then is that $KX$ is the quotient space of $VB(X)$ by the equivalence relation of \emph{stable isomorphism}, namely
$$KX=VB(X)/\simeq _S.$$

Thus "stability" encodes the Grothendieck trick to get an Abelian group from a semi-group.\\ 

{\bf Proposition 2.} If the base space $X$ is connected, compact and  Hausdorff, the  set of equivalence classes under similarity $VB(X)/\sim$ forms another Abelian group under direct sum $\oplus$ denoted $\tilde{K(X)}$.\\ 

{\bf Definition 11.} The group $\tilde{K(X)}$ is called the \emph{reduced} K-Theory of $X$ (the reduced $0$th K-group of $X$).\\

{\bf Proof of Proposition 2:} It will saffice to prove the existence of inverses; this is done by showing that for each vector bundle $\pi :E\rightarrow X$ there exists a vector bundle $E'\rightarrow X$ such that $E\oplus E'\simeq {\mathbb{C}}_{X}^{n}$ for some $n$. Yet this is precisely Proposition 1 above. (We assume that $X$ is connected and thus all fibres have the same dimension, it is not hard to prove that the Proposition 2 still holds even if we relax the connectedness assumption). $\square$\\

If we consider as an extra operation the tensor product between vector bundles, then $Vect(X)$ becomes a semi-ring and then $KX$ is a ring.\\

For the direct sum operation on stably isomorphic K-classes, only the zero element, $[{\mathbb{C}}_{X}^{0}]$ can have an inverse since $E\oplus E'\simeq _S {\mathbb{C}}_{X}^{0}$ implies  $E\oplus E'\oplus  {\mathbb{C}}_{X}^{n} \simeq _S {\mathbb{C}}_{X}^{n}$ for some $n$ which can only happen if both $E$ and $E'$ are zero dimensional. However, even though inverses do not exist, we do have the cancellation property that $E_1\oplus E_2\simeq _S E_1\oplus E_3\Rightarrow E_2\simeq _S E_3$ over a compact space $X$, since we can add to both sides of $E_1\oplus E_2\simeq _S E_1\oplus E_3$ a bundle $E_{1}'$ so that $E_1\oplus E_{1}'\simeq {\mathbb{C}}_{X}^{n}$ for some $n$.\\ 

Thus for compact $X$, K-classes can be represented as "formal differences" of vector bundles of the form $E-E'$ where $E$ and $E'$ are honest real vector bundles over $X$ with the equivalence relation 
$$E_1 - E_{1}'=E_2 -E_{2}'$$
if and only if
$$E_1\oplus E_{2}'\simeq _S E_2\oplus E_{1}'.$$
For this condition to be well-defined we need the cancellation property and thus $X$ has to be compact. With the obvious addition rule
$$(E_1-E_{1}')+(E_2 -E_{2}')=(E_1\oplus E_2)-(E_{1}'\oplus E_{2}')$$
then $KX$ is a group. The zero element is the equivalence class of $E-E'$ for any $E$ and the inverse of $E-E'$ is $E'-E$. Note moreover that every element of $KX$ can be represented also as a difference
$$E-{\mathbb{C}}_{X}^{n}$$
since if we start with $E-E'$ we can add to both $E$ and $E'$ a bundle $E''$ such that $E'\oplus E''\simeq {\mathbb{C}}_{X}^{n}$ for some $n$.\\ 

There is a natural homomorphism $KX\rightarrow \tilde{KX}$ which sends $E-{\mathbb{C}}_{X}^{n}$ to the corresponding $\sim$-class of $E$ which is well defined since if $E-{\mathbb{C}}_{X}^{n}=E' -{\mathbb{C}}_{X}^{m}$ in $KX$, then $E\oplus {\mathbb{C}}_{X}^{m}\simeq _S E'\oplus {\mathbb{C}}_{X}^{n}$, hence $E\sim E'$. This map is obviously surjective and its kernel consists of elements $E-{\mathbb{C}}_{X}^{n}$ with $E\sim {\mathbb{C}}_{X}^{0}$, hence $E\simeq _S {\mathbb{C}}_{X}^{m}$ for some $m$, so the kernel consists of elements of the form 
${\mathbb{C}}_{X}^{m}-{\mathbb{C}}_{X}^{n}$. This subgroup of $KX$ is isomorphic to ${\mathbb{Z}}$. In fact, restriction of vector bundles to a base point $x_0\in X$ defines a homomorphism $KX\rightarrow K(x_0)\simeq{\mathbb{Z}}$ which restricts to an isomorphism on the subgroup $\{{\mathbb{C}}_{X}^{m}-{\mathbb{C}}_{X}^{n}\}$. One thus has a splitting
$$KX\simeq \tilde{KX}\oplus{\mathbb{Z}}$$
depending on the choice of $x_0$.\\

Let $f: X\rightarrow Y$ be a continuous map and let $E\rightarrow Y$ be a vector bundle over $Y$. Then as we saw earlier, one can construct the pull-back bundle $f^{*}(E)$ of $E$ by $f$ over $X$. Thus any such map $f: X\rightarrow Y$ induces a map $f^*:Vect(Y)\rightarrow Vect(X)$ where $\xi\mapsto f^*\xi$, 
which preserves direct sums, i.e. it is a homomorphism of semi-groups and thus it gives a group homomorphism
$f^*:K(Y)\rightarrow K(X)$ between abelian groups. This is often referred to as 
\emph{the wrong way functoriality in K-Theory}.\\

{\bf Theorem 1.} Suppose $f_0,f_1 :X\rightarrow Y$ are homotopic maps and $\xi =(E,\pi ,Y)$ is a vector bundle over $Y$. Then
$f_{0}^*\xi $ is isomorphic to $f_{1}^*\xi $.\\

Recall that $f_0,f_1 :X\rightarrow Y$ homotopic means that there exists a continuous $F:X\times I\rightarrow Y$, 
where $I=[0,1]$, so that $f_0 (x)=F(x,0)$ and $f_1 (x)=F(x,1)$.\\  

{\bf Proof of Theorem 1:} Let $F:X\times I\rightarrow Y$ be a homotopy from $f_0$ to $f_1$. The restrictions of  $F^{*}(E)$ over $X\times \{0\}$ and $X\times \{1\}$ are $f_{0}^*\xi $ and $f_{1}^*\xi $ respectively. Hence, it will saffice to prove the following:\\

{\bf Lemma 2.} The restrictions of a vector bundle $E\rightarrow X\times I$ over $X\times \{0\}$ and $X\times \{1\}$  are isomorphic if $X$ is compact and Hausdorff (more generally if $X$ is paracompact).\\

{\bf Proof of Lemma 2:} We shall use two preliminary facts:\\
{\bf i.} A vector bundle $\pi :E\rightarrow X\times [a,b]$ is trivial if its restrictions over $X\times [a,c]$ and $X\times [c,b]$ are both trivial for some $c\in (a,b)$. To see this, let these restrictions be $E_1 =\pi ^{-1}(X\times [a,c])$ and $E_2 =\pi ^{-1}(X\times [c,b])$ and let $h_1:E_1\rightarrow X\times [a,c]\times{\mathbb{C}}^{n}$ and      
$h_2:E_2\rightarrow X\times [c,b]\times{\mathbb{C}}^{n}$ be isomorphisms. These isomorphisms may not agree on $\pi ^{-1}(X\times\{c\})$ but they can be made to agree by replacing $h_2$ by its composition with the isomorphism $X\times [c,b]\times{\mathbb{C}}^{n}\rightarrow X\times [c,b]\times{\mathbb{C}}^{n}$ which on each slice $X\times \{x\}\times{\mathbb{C}}^{n}$ is given by $h_{1}h_{2}^{-1}:X\times\{c\}\times{\mathbb{C}}^{n}\rightarrow X\times\{c\}\times{\mathbb{C}}^{n}$. Once $h_1$ and $h_2$ agree on $E_1\cap E_2$ they define a trivialisation of $E$.\\

{\bf ii.} For a vector bundle $\pi :E\rightarrow X\times I$ there exists an open cover $\{U_a\}$ of $X$ so that each restriction $\pi ^{-1}(U_a\times I)\rightarrow U_a\times I$ is trivial. This is so because for each $x\in X$ one can find open neighborhoods $U_{x,1},...,U_{x,k}$ in $X$ and a partition $0=t_0<t_1<...<t_k=1$ of $[0,1]$ such that the bundle is trivial over $U_{x,i}\times [t_{i-1},t_i]$ using compactness of $[0,1]$. Then by the fact (i) the bundle is trivial over $U_a\times I$ where $U_a=U_{x,1}\cap ...\cap U_{x,k}$.\\

We now come to the proof of the Lemma 2: By fact (ii) we can choose an open cover $\{U_a\}$ of $X$ so that $E$ is trivial over each $U_a\times I$. We deal first with the simpler case where $X$ is compact Hausdorff: In this case a finite number of $U_a$'s cover $X$. We relabel these as $U_i$, for $i=1,2,...,m$. As shown in Proposition 1 there is a corresponding partition of unity by functions $\phi _i$ whose support is contained in $U_i$. For $i\geq 1$, let $\psi _i =\phi _1 +...+\phi _i$, so in particular $\psi _0=0$ and $\psi _m=1$. Let $X_i$ be the graph of $\psi _i$, i.e. the subspace of $X\times I$ consisting of points of the form $(x,\psi _i(x))$ and let $\pi _i:E_i\rightarrow X_i$ be the restriction of the bundle $E$ over $X_i$. Since $E$ is trivial over $U_i\times I$, the natural projection homeomorphism $X_i\rightarrow X_{i-1}$ lifts to a homeomorphism $h_i:E_i\rightarrow E_{i-1}$ which is the identity outside $\pi _{i} ^{-1}(U_i)$ and which takes each fibre of $E_i$ isomorphically onto the corresponding fibre of $E_{i-1}$. Thus the composition $h=h_1h_2...h_m$ is an isomorphism from the restriction of $E$ over $X\times \{1\}$ to the restriction over $X\times \{0\}$.\\

In the general case where $X$ is only paracompact, there is a countable cover $\{V_i\}_{i\geq 1}$ of $X$ and a partition of unity $\{\phi _i\}$ with $\phi _i$ supported in $V_i$ such that each $V_i$ is a disjoint union of open sets each contained in some $U_a$. This means that $E$ is trivial over each $V_i\times I$. As before we let $\psi _i =\phi _1 +...+\phi _i$ and let $\pi _i:E_i\rightarrow X_i$ be the restriction of $E$ over the grapf of $\psi _i$ and we construct $h_i:E_i\rightarrow E_{i-1}$ using the fact that $E$ is trivial over $V_i\times I$. The infinite composition $h=h_1h_2...$ is then a well-defined isomorphism from the restriction of $E$ over $X\times \{1\}$ to the restriction over $X\times \{0\}$ since near each point $x\in X$ only finitely many $\phi _i$'s are nonzero, hence there is a neighborhood of $x$ in which all but finitely many $h_i$'s are the identity. $\square$\\

{\bf Corollary 1.} If $X$ is contractible (namely homotopic to a point) then $Vect(X)=\mathbb{N}$ and $KX=\mathbb{Z}$.\\

{\bf Proof:} $X$ contractible means that it is homotopic to a point, hence a vector bundle over a point means just a vector 
space; isomorphism classes of vector spaces are characterised by their dimension, thus $Vect(*)=\mathbb{N}$ and then the 
universal (or Grothendieck) group $Gr({\mathbb{N}})={\mathbb{Z}}$. $\square$\\

{\bf Corollary 2.} Vector bundles over contractible spaces are trivial.\\

{\bf Proof:} Straightforward. $\square$\\

[Theorem 1 holds for real vector bundles as well. It also holds for arbitrary fibre bundles].\\

We defined vector bundles over topological spaces and we defined various equivalence relations among them (isomorphism, stable isomorphism, similarity). An important question is the \emph{classification of vector bundles} over some fixed topological space say $X$: Classification up to similarity can be achieved via the reduced K-Theory, classification up to stable isomorphism can be done via K-Theory but the problem of classification up to isomorphism is still largely an open question. One very useful tool in tackling this problem are some cohomological objects called \textsl{"characteristic classes"} which we shall meet later. They give a partial answer to the question of classification of vector bundles up to isomorphism, in fact they work very well in low dimensions. For the moment we shall rephrase the problem of classification in terms of a standard concept of algebraic topology, the idea of \textsl{homotopy classes of maps} along with the notions of the \textsl{classifying space} and the \textsl{universal bundle}.  Thus we shall construct a $k$-dim vector bundle (called the universal bundle) $E_k\rightarrow G_k$ with the property that all $k$-dim vector bundles over a compact Hausdorff space can be obtained as pull-backs of this single bundle (in particular this can be generalised for paracompact spaces and the case $k=1$ gives the line bundle over the infinite projective space in the real case).\\

Recall that $G_k({\mathbb{C}}^n)$ denotes the \emph{Grassmannian} of $k$-dim vector subspaces of ${\mathbb{C}}^n$ for nonnegative integers $k\leq n$ (the set of $k$-dim planes passing through the origin). One can gine a topology to the Grassmannian using the \emph{Stiefel manifold} $V_k({\mathbb{C}}^n)$, the space of orthonormal $k$-frames in ${\mathbb{C}}^n$, in other words the $k$-tuples of of orthonormal vectors in ${\mathbb{C}}^n$. This is a subspace of the product of $k$ copies of the unit sphere $S^{n-1}$, namely the subspace of orthonormal $k$-tuples. It is a closed subspace since orthogonality of two vectors can be expressed by an algebraic equation. Hence the Stiefel manifold is compact (since the products of spheres are compact). There is a natural surjection $V_k({\mathbb{C}}^n)\rightarrow  G_k({\mathbb{C}}^n)$ sending a $k$-frame to the subspace it spans and thus $G_k({\mathbb{C}}^n)$ can be topologised by giving it the quotient topology with respect to this surjection and hence the Grassmannian is also compact.\\

The inclusions ${\mathbb{C}}^{n}\subset {\mathbb{C}}^{n+1}\subset ... $  give inclusions $G_k({\mathbb{C}}^n)\subset G_k({\mathbb{C}}^{n+1})\subset ...$ and we let 
$$G_k({\mathbb{C}}^{\infty})=\bigcup _{n}G_k({\mathbb{C}}^n).$$
We give $G_k({\mathbb{C}}^{\infty})$ the weak (or direct limit topology), so a set is open if and ony if it intersects each 
$G_k({\mathbb{C}}^{n})$ in an open set.\\

There are canonical $n$-dim vector bundles over $G_k({\mathbb{C}}^n)$: We define
$$E_k({\mathbb{C}}^n)=\{(l,v)\in G_k({\mathbb{C}}^n)\times{\mathbb{C}}^{n}:v\in l\}.$$
The inclusions ${\mathbb{C}}^{n}\subset {\mathbb{C}}^{n+1}\subset ... $  give inclusions $E_k({\mathbb{C}}^n)\subset E_k({\mathbb{C}}^{n+1})\subset ...$ and we set
$$E_k({\mathbb{C}}^{\infty})=\bigcup _{n}E_k({\mathbb{C}}^n)$$
 again with the direct limit topology.\\

{\bf Lemma 3.} The projection $p:E_k({\mathbb{C}}^n)\rightarrow G_k({\mathbb{C}}^n)$ with $p(l,v)=l$ defines a vector bundle both for finite and infinite $n$.\\

{\bf Proof:} Suppose first that $n$ is finite. For $l\in G_k({\mathbb{C}}^n)$, let $\pi _l:{\mathbb{C}}^n\rightarrow l$ be an orthogonal projection and let $U_l=\{l'\in G_k({\mathbb{C}}^n) : \pi _l(l')$ has dim $k\}$. In particular $l\in U_l$. We shall show that $U_l$ is open in $G_k({\mathbb{C}}^n)$ and that the map $h:p^{-1} (U_l)\rightarrow U_{l}\times l\simeq U_l\times {\mathbb{C}}^k$ defined by $h(l',v)=(l',\pi _{l}(v))$ is a local trivialisation of $E_k({\mathbb{C}}^n)$.\\
For $U_l$ to be open is equivalent to its preimage in $V_k({\mathbb{C}}^n)$ being open. This preimage consists of orthonormal frames $v_1,...,v_k$ such that $\pi _l(v_1),...,\pi _l(v_k)$ are independent. Let $A$ be the matrix of $\pi _l$ with respect to the standard basis in the domain ${\mathbb{C}}^n$ and any fixed basis in the range $l$. The condtition on $v_1,...,v_k$ is then that the $k\times k$ matrix with columns $Av_1,...,Av_k$ have nonzero determinant. Since the value of this determinant is obviously a continuous function of $v_1,...,v_k$, it follows that the frames $v_1,...,v_k$ yielding a nonzero determinant form an open set in $V_k({\mathbb{C}}^n)$.\\
It is clear that $h$ is a bijection which is a linear isomorphism on each fibre. We need to check that $h$ and $h^{-1}$ are continuous. For $l'\in U_l$ there is a unique invertible linear map $L_{l'}:{\mathbb{C}}^n\rightarrow {\mathbb{C}}^n$ restricting to $\pi _l$ on $l'$ and the identity on $l^{\perp}=ker(\pi _l)$. We claim that $L_{l'}$ regarded as a $n\times n$ matrix depends continuously on $l'$ namely we can write 
$$L_{l'}=AB^{-1}$$
where $B$ sends the standard basis to $v_1,...,v_k,v_{k+1},...,v_n$ with $v_1,...,v_k$  an orthonormal basis for $l'$ and $v_{k+1},...,v_n$ a fixed basis for $l^{\perp}$ and $A$ sends the standard basis to $\pi _{l}(v_1),...,\pi _l(v_k),v_{k+1},...,v_n$.\\
Both $A$ and $B$ depend continuously on $v_1,...,v_k$. Since matrix multiplication and matrix inversion are continuous operations, it follows that the product $L_{l'}=AB^{-1}$ depends continuously on $v_1,...,v_k$. Yet since $L_{l'}$ depends only on $l'$ and not on the basis $v_1,...,v_k$ for $l'$, it follows that $L_{l'}$ depends continuously on $l'$ since $G_k({\mathbb{C}}^n)$ has the quotient topology from $V_k({\mathbb{C}}^n)$. Since we have $h(l',v)=(l',\pi _l(v))=(l', L_{l'}(v))$, we see that $h$ is continuous. Similarly $h^{-1}(l',w)=(l',L_{l'}^{-1}(w))$ and $L_{l'}^{-1}$ depends continuously on $l'$, matrix inversion being continuous, so $h^{-1}$ is continuous. This completes the proof for finite $n$.\\

For the infinite case one takes $U_l$ to be the union of the $U_l$'s for increasing $n$. The local trivialisations $h$ constructed above for the finite case then fit together to give a local trivialisation over this $U_l$ whereas continuity is manifest since we use the weak topology. $\square$\\

We shall be interested in the case $n=\infty$ now and to simplify our notation we shall write $G_k$ for $G_k({\mathbb{C}}^{\infty})$ and similarly we shall write $E_k$ for $E_k({\mathbb{C}}^{\infty})$. As we have already done previously, we denote by $[X,Y]$ the set of homotopy classes of maps $f:X\rightarrow Y$.\\

{\bf Theorem 2.} For compact Hausdorff $X$ (more generally if $X$ is paracompact), the map $[X,G_k]\rightarrow Vect _k(X)$ with $[f]\mapsto f^{*}(E_k)$ is a bijection.\\

{\bf Proof:} The key observation is the following: For a $k$-dim vector bundle $\pi :E\rightarrow X$, an isomorphism $E\simeq f^{*}(E_k)$ is equivalent to a map $g:E\rightarrow {\mathbb{C}}^{\infty}$ that is a linear injection on each fibre. To see this, suppose first that we have a map $f:X\rightarrow G_k$ along with an isomorphism $E\simeq f^{*}(E_k)$. Then we have a commutative diagram  
$$\begin{CD}
E @>\simeq >> f^{*}(E_k) @>\tilde{f}>> E_k  @>\pi>> {\mathbb{C}}^{\infty}\\
  @VpVV   @VVV     @VVV  @.\\
X @>>id> X @>>f> G_k @. 
\end{CD}$$\\
where $\pi (l,v)=v$. The composition across the top row is a map $g:E\rightarrow {\mathbb{C}}^{\infty}$ which is a linear injection on each fibre since both $\tilde{f}$ and $\pi$ have this property. Conversely, given a map $g:E\rightarrow {\mathbb{C}}^{\infty}$ which is a linear injection on each fibre, we define $f:X\rightarrow G_k$ by letting $f(x)$ be the $k$-plane $g(p^{-1}(x))$ and this clearly yields a commutative diagram as above.\\
To show the surjectivity of the map $[X,G_k]\rightarrow Vect_{k}(X)$, suppose $p:E\rightarrow X$ is a $k$-dim vector bundle and let $\{U_i\}$ be an open cover of $X$ such that $E$ is trivial over each $U_i$. Then by paracompactness of $X$ there is a countable open cover $\{U_a\}$ of $X$ such that $E$ is again trivial over each $U_a$ and there is a partition of unity $\{\phi _a\}$ with $\phi _a$ supported in $U_a$. Let $g_a:p^{-1}(U_a)\rightarrow {\mathbb{C}}^k$ be the composition of a trivialisation $p^{-1}(U_a)\rightarrow U_a\times{\mathbb{C}}^k$ with the projection onto ${\mathbb{C}}^k$. The map $(\phi _{a}p)g_a$ with $v\mapsto\phi _a(p(v))g_{a}(v)$ extends to a map $E\rightarrow {\mathbb{C}}^k$ which is zero outside $p^{-1}(U_a)$. Near each point of $X$ only finitely many $\phi _a$'s are nonzero and at least one $\phi _a$ is nonzero, so these extended $(\phi _{a}p)g_a$'s constitute the coordinates of a map $g:E\rightarrow ({\mathbb{C}}^k)^{\infty}={\mathbb{C}}^{\infty}$ which is a linear injection on each fibre.\\
For injectivity, if one has isomorphisms $E\simeq f_{0}^{*}(E_k)$ and          
$E\simeq f_{1}^{*}(E_k)$ for two maps $f_0,f_1:X\rightarrow G_k$, then these give maps $g_0,g_1:E\rightarrow {\mathbb{C}}^{\infty}$ which are linear injections on fibres, as in the first part of the proof. Then the claim is that $g_0$ and $g_1$ are homotopic through maps $g_t$ which are again linear injections on fibres. If this is so, then $f_0$ and $f_1$ will be homotopic via $f_t(x)=g_t(p^{-1}(x))$.\\
To prove the claim, the first step is to construct a homotopy $g_t$ by taking the composition of $g_0$ with the homotopy $L_t:{\mathbb{C}}^{\infty}\rightarrow {\mathbb{C}}^{\infty}$ defined by $L_t(x_1,x_2,...)=(1-t)(x_1,x_2,...)+t(x_1,0,x_2,...)$. For each $t$ this is a linear map whose kernel is easily computed to be $0$, so $L_t$ is injective. Composing $L_t$ with $g_0$ moves the image of $g_0$ into the odd-numbered coordinates. Similarly we can homotope $g_1$ into the even-numbered coordinates. We keep denoting the new $g$'s by $g_0$ and $g_1$ and we let $g_t=(1-t)g_0+tg_1$. This is linear and injective on fibres $\forall t$ since $g_0$ and $g_1$ are. $\square$\\

This theorem indicates that $k$-dim vector bundles over some fixed topological space are classified by homotopy classes of maps into $G_k$ and because of this $G_k$ is called the \emph{classifying space} for $k$-dim vector bundles and $E_k\rightarrow G_k$ is called the \emph{universal vector bundle}.  The truth however is that this theorem is of limited usefulness in enumerating all the different vector bundles (of fixed dim) over a given space since explicit calculations of $[X, G_k]$ are usually beyond technical reach. Its importance is due more to its theoretical implications since among other things it can reduce the proof of a general statement to the special case of the universal bundle.\\

The above construction can be applied to the real case as well; there is also a version for oriented real vector bundles (see \cite{hatcher}).\\

We would like to close this section by describing a construction of complex vector bundles over spheres using the so called \emph{clutching functions}. This is a technical point which we shall use in the next chapter.\\

We want to construct a (real or complex, we shall take the complex case here) vector bundle $E\rightarrow S^k$ over the $k$-sphere $S^k$.  We write the sphere as the union of its upper $D^{k}_{+}$   
 and lower hemispheres $D^{k}_{-}$ with $D^{k}_{+}\cap D^{k}_{-}=S^{k-1}$. Given a map $f:S^{k-1}\rightarrow GL_n({\mathbb{C}})$, let $E_f$ be the quotient of the disjoint union $D^{k}_{+}\times {\mathbb{C}}^{n}\cup D^{k}_{-}\times {\mathbb{C}}^{n}$ obtained by identifying $(x,v)\in\partial D^{k}_{-}\times {\mathbb{C}}^{n}$ with $(x,f(x)(v))\in \partial D^{k}_{+}\times {\mathbb{C}}^{n}$. There is then a natural projection $E_f\rightarrow S^k$ and this is an $n$-dim vector bundle. The map $f$ used above is called the \textsl{clutching function} (since it does essentially what the clatch does in vehicles). In fact one can prove that the map $\Phi: f\mapsto E_f$ gives a bijection 
$$\Phi :[S^{k-1},GL_{n}({\mathbb{C}})]\rightarrow Vect_{n}(S^k)$$
(One can find the proof of this statement in \cite{hatcher}; in fact this is an interesting special case of Theorem 1 above). This bijection does not quite work for the real case since $GL_{n}({\mathbb{R}})$ is not path connected. However it works for oriented real vector bundles (we refer to \cite{hatcher} for more details).\\

\newpage

\section{Lecture 2 (Homotopy, Bott Periodicity and Cohomological Properties)}

We assume complex vector bundles with compact base. If the base is the disjoint union
$$X =\coprod X_a ,$$
then
$$KX=\prod KX_a .$$
Moreover we saw that if $G_k =\bigcup _n G_k({\mathbb{C}}^n)$, then
$$Vect _{k}(X)=[X,G_k]$$
(to be precise we mean equality of the cardinalities of the above sets).\\

We shall study the homotopic interpretation of the group $KX$.\\

Suppose that the base $X$ has a base point $x_0$. Then we denote by $\tilde{K}$ the reduced K-Theory (see Definition 2 below)
$$\tilde{K}X=Ker(KX\rightarrow K(x_0)).$$
Yet $K(*)={\mathbb{Z}}$ whereas the dimension of the vector bundle defines a map 
$ÊX \rightarrow\mathbb{Z}$ where 
$[\xi ]-[\eta ]\mapsto dim(\xi )-dim(\eta )$.\\

{\bf Definition 1.} Two vector bundles $\xi$ and $\eta$ over $X$ are called \emph{stably isomorphic} 
if there exist integers $N,M\in\mathbb{Z}$ such that
$$\xi\oplus{\mathbb{C}}_{X}^{M}=\eta\oplus{\mathbb{C}}_{X}^{N}.$$ 

{\bf Definition 2.} The \emph{reduced K-group} $\tilde{K}X$ of $X$ is defined as the set of 
stably isomorphism classes of vector bundles over $X$.\\

We know that for the Grassmannians one has the following inclusions:
$$G_{k}({\mathbb{C}}^n)\hookrightarrow G_{k}({\mathbb{C}}^{n+1})$$
and taking the inductive limit we get
$$G_{k}=\lim _{n\rightarrow\infty}G_{k}({\mathbb{C}}^{n})= \bigcup _{n}G_{k}({\mathbb{C}}^n).$$
Moreover
$$G_{k}({\mathbb{C}}^n)\hookrightarrow G_{k+1}({\mathbb{C}}^{n+1})$$
thus
$$G_{k+1}=\lim _{n\rightarrow\infty}G_{k+1}({\mathbb{C}}^{n+1})= \bigcup _{n}G_{k+1}({\mathbb{C}}^{n+1}).$$
Hence
$$G_k\hookrightarrow G_{k+1},$$
which enables us to define the inductive limit 
$$G_{\infty}=\lim _{k\rightarrow\infty}G_k=\bigcup _{k}G_k .$$
Intuitively, $G_{\infty}$ is the double inductive limit of $G_{k}({\mathbb{C}}^n)$ when both 
variables $k$ and $n$ tend to infinity.\\

Then one has the following:\\

{\bf Theorem 1.}
$$\tilde{K}X=[X, G_{\infty}].$$

Next we shall study vector bundles over the suspension $SX$ of a space $X$.\\

Let $X$ be a topological space; we denote by $CX$ the cone of $X$ which is defined as
$$CX=\frac{X\times [0,1]}{X\times 0}.$$
Then the suspension $SX$ of $X$ is defined by
$$SX=C_{+}X\cup _{X}C_{-}X.$$
Suppose $f:X\rightarrow GL_k(\mathbb{C})$ is a map. We define the set
$$V(f)=\frac{C_{+}X\times{\mathbb{C}}^{k}\cup C_{-}X\times{\mathbb{C}}^k}{(x,v)\equiv (x,f(x)v)}.$$ 
Then $V(f)\rightarrow SX$ is a vector bundle.\\

{\bf Lemma 1.} Let $f_0,f_1:X\rightarrow GL_k(\mathbb{C})$ be homotopic maps. Then 
$$V(f_0)\simeq V(f_1).$$
Moreover one has the following important result:\\

{\bf Theorem 2.}
$$Vect_{k}(SX)=[X, GL_k(\mathbb{C})].$$
If we consider the inclusion
$$GL_k({\mathbb{C}})\hookrightarrow GL_{k+1}({\mathbb{C}})$$
given by
$$A\mapsto
\left( \begin{array}{cc}A & 0\\
0 & 1\end{array} \right)$$
and set
$$GL{\mathbb{C}}=\lim_{k\rightarrow\infty}GL_{k}({\mathbb{C}}) =\bigcup_{k}GL_{k}({\mathbb{C}}),$$
then one has an extension of Theorem 2:\\

{\bf Theorem 3.}
$$\tilde{K}(SX)=[X,GL\mathbb{C}].$$

For the proofs of Theorems 2 and 3 one can see Atiyah's notes.\\

{\bf Remark:}\\
We mentioned that
$$\tilde{K}X=[X,G_{\infty}]$$ 
whereas
$$\tilde{K}SX=[X,GL\mathbb{C}].$$ 
Yet
$$\tilde{K}SX=[SX,G_{\infty}]=[X,\Omega G_{\infty}],$$ 
where $\Omega G_{\infty}$ the loop space of $G_{\infty}$, namely
$$\Omega G_{\infty}=Maps_{*}(S^1 ,G_{\infty}).$$   
Thus
$$\Omega G_{\infty}\simeq GL\mathbb{C},$$
which means that the \emph{infinite Grassmannian $G_{\infty}$ is a classifying space of the infinite group $GL\mathbb{C}$}:
$$G_{\infty}=BGL\mathbb{C}.$$ 

The basic result in K-Theory is \textit{Bott periodicity} which we shall study next. In its simplest form it states that there is an isomorphism between $KX\otimes K(S^2)$ and $K(X\times S^2)$.\\

We denote by $S^{2}X$ the second suspension of $X$, in other words $S^{2}X=S(SX)$. Inductively one can define the nth suspension $S^{n}X$ of a space $X$. Recall from topology that for the suspensions of spheres one has the following results:
$$S^{2n}(S^{1})=S^{2n+1}$$
and
$$S^{2n-2}(S^{2})=S^{2n}.$$

As our first version of Bott periodicity we mention the following:
$$\tilde{K}(S^{2}X)\simeq\tilde{K}X$$
whereas
$$[S^1,G_{k}({\mathbb{C}}^n)]=\pi _{1}(G_{k}({\mathbb{C}}^n))=0$$
because Grassmannians are contractible.\\

Furthermore
$$[S^2,G_{k}({\mathbb{C}}^n)]=\pi _{2}(G_{k}({\mathbb{C}}^n))=\mathbb{Z},$$
$$\tilde{K}S^1=0,$$
$$\tilde{K}S^2=\mathbb{Z},$$
$$\tilde{K}S^{2n+1}=0$$
and
$$\tilde{K}S^{2n}=\mathbb{Z}.$$
Now we come to the cohomological properties of K-Theory.\\

Let $X$ be compact and $A\subset X$ is a compact subspace. We choose some $n\in\mathbb{N}$ and we have:\\

{\bf Definition 3.} The negative powers of K-Theory are defined as follows:
$$\tilde{K}^{-n}X=\tilde{K}(S^{n}X).$$
For the relative K-Theory one has the following definition:\\

{\bf Definition 4.} The relative K-groups are defined as follows: 
$$\tilde{K}^{-n}(X,A)=\tilde{K}^{-n}(X/A)=\tilde{K}(S^{n}(X/A)).$$
Moreover one has that
$$K^{-n}X=\tilde{K}^{-n}(X,\emptyset ),$$
where we define
$$(X,\emptyset )=X^{+}$$
and $X^{+}$ denotes $X$ with a disjoint base point adjoined, i.e. 
$$\tilde{K}(X^{+})=KX.$$
Using Bott periodicity we can extend the above definition and get 
$$K^{n}(X,A)=K^{n-2}(X,A),\forall n\in\mathbb{Z}.$$
To summarise:\\

{\bf 1.} The groups $K^{n}(X,A)$ can be defined for any $n$ whereas $K^{n}(X,A)\simeq K^{n+2}(X,A)$.\\ 

{\bf 2.} Essentially, due to Bott periodicity, there are only two K-groups, $K^{0}(X,A)$ and $K^{1}(X,A)$ where
$K^{0}(X,A)=[X/A,G_{\infty}]$ and $K^{1}(X,A)=[X/A,GL\mathbb{C}]$.\\

Given a map $f:(X,A)\rightarrow (Y,B)$, there exists an induced map which is a homomorphism 
$$K^{n}(X,A)\xleftarrow{f^*} K^{n}(Y,B)$$
(this is the wrong way funcoriality in K-Theory) where at the same time given
$$X/A\leftarrow X\leftarrow A$$
we can get an exact K-Theory sequence of a pair
$$\begin{CD}
K^{0}(X,A) @>>> K^0X @>>> K^0A\\
  @AAA   @.     @VV\delta V\\
K^1A @<<< K^1X @<<< K^1(X,A)
\end{CD}$$\\
\\

\newpage

\section{Lecture 3 (Products and Bott Periodicity Revised)}

Let $X$ be compact and connected. From the set of isomorphism classes of vector bundles $Vect(X)$ over $X$ we defined $K^{0}X=[X,{\mathbb{Z}}\times G_{\infty}]$ whereas the reduced 0th K-group is defined by
$$\tilde{K}^{0}X=Ker(K^{0}X\rightarrow K^{0}(*))=[X,G_{\infty}].$$
Furthermore we saw that
$$\tilde{K}^{1}X=\tilde{K}^{0}(SX)=[X,GL\mathbb{C}]$$
whereas for a compact subspace $A\subset X$ and some $n\in\mathbb{N}$ one has
$$\tilde{K}^{-n}X=\tilde{K}^{0}(S^{n}X),$$
$$\tilde{K}^{-n}(X,A)=\tilde{K}^{-n}(X/A),$$
$$K^{-n}X=K^{-n}(X,\emptyset),$$
where $X^{+}=X/\emptyset$ and 
$$K^{n}(X,A)=K^{n-2}(X,A)$$
using Bott periodicity.\\

We list the properties of the groups $K^{n}(X,A)$ (which follow from the Eilenberg-Steenroad axioms):\\

{\bf 1.} The correspondence $(X,A)\mapsto K^{n}(X,A)$ is a functor from the category of topological spaces to the category of abelian groups.\\

{\bf 2.} Homotopic maps induce the same homomorphisms between abelian groups.\\

{\bf 3.} There is an exact sequence for the K-Theory of pairs (see the end of the previous lecture).\\

{\bf 4.} The exact sequence is natural.\\

{\bf 5.} Excision is satisfied:
$$K^{n}(X,A)\simeq K^{n}(X/A,*).$$

{\bf 6.} The dimension axiom fails
$$K^{n}(*)=0, n\in{\mathbb{N}}, odd$$
whereas
$$K^{n}(*)={\mathbb{Z}}, n\in{\mathbb{N}}, even.$$
Thus we have a \emph{generalised (or exotic)} cohomology theory.\\

Given 
$$A\rightarrow X\rightarrow X/A,$$
the exact sequence for pairs was
$$\begin{CD}
K^{0}(X,A) @>>> K^0X @>>> K^0A\\
  @AAA   @.     @VV\delta V\\
K^1A @<<< K^1X @<<< K^1(X,A)
\end{CD}.$$\\
\\
If $A\subset X$, then
$$SA\simeq\frac{X\cup CA}{X}$$
while
$$SA\simeq\frac{X\cup CA}{X}\leftarrow X\cup CA\rightarrow X/A.$$
Furthermore
$$K^{-n-1}A=K^{-n}(SA)\simeq K^{-n}(X\cup CA,X)\rightarrow K^{-n}(X\cup CA)\simeq K^{-n}(X/A)$$
along with the maps
$$K^{-1}A\rightarrow K^{0}(X/A)$$
and
$$K^{0}A\rightarrow K^{1}(X/A).$$
The proof of exactness is tedious.\\

Next we shall try to give a geometric interpretation of the conecting maps $\delta$ which appear in the K-Theory exact sequence.\\

For $K^{0}(X,A)$, we start with a pair $(\xi ,\phi )$ where $\xi$ is a vector bundle over $X$ and $\phi$ is a trivialisation of $\xi /A$. There are two equivalence relations between these pairs, isomorphism and homotopy, thus one can define the class $[\xi ,\phi ]\in K^{0}(X,A)$.\\

Consider the map   
$$\delta :K^{1}A\rightarrow K^{0}(X,A)$$
with some
$$f:A\rightarrow GL_{k}(\mathbb{C}).$$
Define
$$\frac{X\times{\mathbb{C}}^{k}\coprod CA\times{\mathbb{C}}^{k}}{glue via f}$$
and obtain that
$$X\cup CA\simeq X/A.$$

Now let us study the products.\\

Let
$$Vect(X)\times Vect(Y)\rightarrow Vect(X\times Y)$$
where
$$(\xi ,\eta )\mapsto \xi\otimes\eta$$
is the external tensor product. For fibres one has 
$$F_{x,y}(\xi\otimes\eta )=F_{x}(\xi )\otimes F_{y}(\eta ),$$
since ç $\xi\otimes\eta\rightarrow X\times Y$ is a new vector bundle.\\

Then one has an induced map
$$K^{n}(X,A)\otimes K^{m}(Y,B)\rightarrow K^{n+m}(X\times Y,A\times Y\cup X\times B).$$
{\bf 1.} Extend the above to a map
$$K^{0}X\otimes K^{0}Y\rightarrow K^{0}(X\times Y)$$
with
$$([\xi _1]-[\eta _1])\otimes ([\xi _2]-[\eta _2]):=[\xi _1\otimes\xi _2 ]+[\eta _1\otimes\eta _2 ]-[\eta _1\otimes\xi _2 ]-
[\xi _1\otimes\eta _2 ].$$
Recall that the \emph{smash product} in algebraic topology $X\wedge Y$ where $X$ has a base point $x_0$ 
and $Y$ has a base point $y_0$ is defined by
$$X\wedge Y=\frac{X\times Y}{x_0\times Y\cup X\times y_0}.$$
For example, the suspension of a space $X$ is the smash product with the circle:
$$S^1\wedge X\simeq SX$$
and
$$S^n\wedge X\simeq S^{n}X.$$
{\bf 2.} Let us focus on the reduced K-Theory: We extend our map now to
$$\tilde{K}^{0}X\otimes\tilde{K}^{0}Y\rightarrow\tilde{K}^{0}(X\wedge Y).$$
We define the smash product as
$$\tilde{K}^{0}(X\wedge Y)=Ker(K^{0}(X\times Y)\rightarrow K^{0}X\oplus K^{0}Y).$$
We use the obvious maps 
$$K^{0}(X\times Y)\rightarrow K^{0}X$$
from
$$X=X\times y_0\rightarrow X\times Y$$
along with
$$K^{0}(X\times Y)\rightarrow K^{0}Y$$
from
$$Y=Y\times x_0\rightarrow X\times Y.$$
Hence we obtain:
$$\tilde{K}^{0}(S^{n}(X/A))\otimes\tilde{K}^{0}(S^{m}(Y/B))\rightarrow\tilde{K}^{0}(S^{n}(X/A)\wedge S^{m}(Y/B))=
\tilde{K}^{0}(S^{n+m}\wedge (X/A)\wedge (Y/B)),$$
where
$$\frac{X\times Y}{A\times Y\cup X\times B}=\frac{X}{A}\wedge\frac{Y}{B}.$$
As an example consider the special case
$X=Y$ êáé $A=B=\emptyset$. We get:
$$K^{n}X\otimes K^{m}X\rightarrow K^{n+m}(X\times X)\xrightarrow{\Delta ^*} K^{n+m}X$$
where the last map 
$$\Delta ^{*}:K^{n+m}(X\times X)\rightarrow K^{n+m}X$$
is induced by the map
$$\Delta :X\rightarrow X\times X$$
with
$$\Delta (x)=(x,x).$$
Hence we have the following maps
$$K^{0}X\otimes K^{0}X\rightarrow K^{0}X,$$
$$K^{0}X\otimes K^{1}X\rightarrow K^{1}X,$$
along with
$$K^{1}X\otimes K^{1}X\rightarrow K^{0}X,$$
and this is esentially Bott periodicity.\\

We can also define
$$K^{*}X=K^{0}X\oplus K^{1}X$$
which is a graded commutative ring.\\

We are now in a position to see a second revised version of Bott periodicity:\\

The map
$$\tilde{K}^{0}(S^2)\otimes\tilde{K}^{0}X\rightarrow\tilde{K}^{0}(S^{2}X)$$
is an isomorphism (recall that $\tilde{K}^{0}S^2={\mathbb{Z}}$).\\

Bott periodicity is perhaps the most useful tool in computations. Let us mention some examples:\\
 
$$\tilde{K}^{1}(S^{2n+1})={\mathbb{Z}}=\tilde{K}^{0}(S^{2n})$$
while
$$\tilde{K}^{1}(S^{2n})=0=\tilde{K}^{0}(S^{2n-1}).$$
Moreover
$$[S^{2n+1},GL{\mathbb{C}}]=\pi _{2n+1}(GL{\mathbb{C}})=\pi _{2n+1}(GL_{N}({\mathbb{C}})), N>2n+2),$$
whereas
$$[S^{2n},GL{\mathbb{C}}]=0.$$

If we denote by $U_N$ the group of unitary $N\times N$ complex matrices, this is onviously a subgroup of $GL_{N}({\mathbb{C}})$), then one has (for $N>2n+1$):
$$\pi _{2n+1}(U_N)={\mathbb{Z}},$$
whereas
$$\pi _{2n}(U_N)=0.$$

For the special orhtogonal group $SO_N$ we have (for $N>i+1$):
$$\pi _{i}(SO_N)={\mathbb{Z}}/2, i\equiv 1 (8),$$
$$\pi _{i}(SO_N)=0, i\equiv 2 (8),$$
$$\pi _{i}(SO_N)={\mathbb{Z}}, i\equiv 3 (8),$$
$$\pi _{i}(SO_N)=0, i\equiv 4 (8),$$
$$\pi _{i}(SO_N)=0, i\equiv 5 (8),$$
$$\pi _{i}(SO_N)=0, i\equiv 6 (8),$$
$$\pi _{i}(SO_N)={\mathbb{Z}}, i\equiv 7 (8) $$
and
$$\pi _{i}(SO_N)={\mathbb{Z}}/2, i\equiv 1 (8).$$

\newpage

\section{Lecture 4 (Chern Character and Chern Classes)}

We saw that $K^{n}(X,A)$ is a multiplicative generalised cohomology theory. Moreover we identified the maps 
$$K^{n}(X,A)\times K^{m}(Y,B)\rightarrow K^{n+m}(X\times Y, X\times B\cup A\times Y).$$
It is true that
$$K^{*}(X,A)\otimes\mathbb{Q}$$
is an ordinary cohomology theory.\\

We assume that $X$ is a "nice" topological space (for example a finite simplicial complex or finite CW-complex).\\

Then we define:
$$H^{ev}(X;{\mathbb{Q}})=\bigoplus_{n\geq 0}H^{2n}(X;{\mathbb{Q}})$$
and
$$H^{odd}(X;{\mathbb{Q}})=\bigoplus_{n\geq 0}H^{2n+1}(X;{\mathbb{Q}}).$$
One then has the following fundamental result:\\

{\bf Theorem 1.} There exist natural isomorphisms
$$ch:K^{0}X\otimes{\mathbb{Q}}\rightarrow H^{ev}(X;{\mathbb{Q}})$$
and
$$ch:K^{1}X\otimes{\mathbb{Q}}\rightarrow H^{ïdd}(X;{\mathbb{Q}})$$
which preserve products. These isomorphisms are called \emph{Chern characters}.\\ 

{\bf Proof:} There is a scetch of the proof at the end of this lecture $\square$.\\

To continue we shall need the characteristic classes of complex line bundles. Let $L\rightarrow X$ 
be a complex line bundle (we shall denote it simply by $L_X$ below) and let 
$$Vect_{1}X=[X,G_1],$$ 
where 
$$f_{L}:X\rightarrow G_{1}=\bigcup _{n}G_{1}({\mathbb{C}}^n)=\bigcup _{n}{\mathbb{C}}P^{n-1}:={\mathbb{C}}P^{\infty}.$$
Then
$$H^{*}({\mathbb{C}}P^{\infty};{\mathbb{Z}})={\mathbb{Z}}[u],$$
where
$$u\in H^{2}({\mathbb{C}}P^{\infty};{\mathbb{Z}}).$$
{\bf Definition 1.} The first $Chern$ class of $L$ is 
$$c_1{L}=f_{L}^{*}(u)\in H^2(X;{\mathbb{Z}}).$$ 
{\bf Lemma 1.} The map
$$c_{1}:Vect_{1}(X)\rightarrow H^2(X; {\mathbb{Z}})$$ 
is a group isomorphism.\\

Given two line bundles $L_1$ and $L_2$, we can form their tensor product $L_1\otimes L_2$ where ${\mathbb{C}}_{X}$ is the unit for the tensor product $\otimes$. We denote by $L_{X}^{*}:=Hom(L_{X},{\mathbb{C}})$ the duall complex line bundle. Then
$$L^*\otimes L\simeq{\mathbb{C}}_X .$$
Furthermore
$$c_{1}(L_1\otimes L_2)=c_{1}L_{1}+c_{2}L_{2}.$$
Given a vector bundle $\xi\rightarrow X$ of dimension $n$, we consider the bundle $\pi :{\mathbb{P}}(\xi )\rightarrow X$ 
with fibre $\pi ^{-1}(x)={\mathbb{P}}(F_x(\xi ))\simeq{\mathbb{C}}P^{n-1}$. The induced map 
$$\pi ^{*}:H^{*}(X)\rightarrow H^{*}({\mathbb{P}}(\xi ))$$
is a ring homomorphism which makes $H^{*}({\mathbb{P}}(\xi ))$ a module over the ring $H^{*}(X)$ (with integer coefficients).\\

Next consider the map
$$H^{*}(X)\otimes H^{*}({\mathbb{P}}(\xi ))\rightarrow H^{*}({\mathbb{P}}(\xi )),$$
with
$$x\otimes y\mapsto\pi ^{*}(x)y.$$
There exists a line bundle $L_{\xi}$ over ${\mathbb{P}}(\xi )$ with $1\in H^{0}({\mathbb{P}}(\xi ))$, 
$c_{1}(L_{\xi})=u\in H^{2}({\mathbb{P}}(\xi ))$ and we form $u^2$, $u^3$ etc. Then one has the following:\\

{\bf Theorem 2.}  $H^{*}({\mathbb{P}}(\xi ))$ is a free module over $H^{*}(X)$ with basis 
$1, u, u^2,...u^{n-1}$.\\

{\bf Proof:} We describe the basic idea.\\

Case 1: $\xi =X\times{\mathbb{C}}^{1}$, then
${\mathbb{P}}(\xi )=X\times{\mathbb{C}}P^{n-1}$ whereas $H^{*}(X\times{\mathbb{C}}P^{n-1})=H^{*}(X)\otimes H^{*}({\mathbb{C}}P^{n-1})$ 
can be computed using the Kunneth formula in algebraic topology. To conclude the proof we use  a Mayer-Vietories inductive argument. $\square$.\\

{\bf Definition 2.} (\emph{Chern classes}). Let $\xi\rightarrow X$ be a vector bundle, we form    
${\mathbb{P}}(\xi )$ and take the class $u\in H^{2}({\mathbb{P}}(\xi );{\mathbb{Z}})$. Then
$$-u^{n}=c_{1}(\xi )u^{n-1}+c_{2}(\xi )u^{n-2}+...+c_{k}(\xi )u^{n-k}+...+c_{n}(\xi ).$$ 
The above equation defines Chern classes $c_{i}(\xi )\in H^{2i}(X;{\mathbb{Z}})$.\\

{\bf Problem:} Check the following properties of Chern classes:\\
{\bf 1.} (pull-backs) $c_{i}(f^{*}(\xi ))=f^{*}(c_{i}(\xi ))$.\\
{\bf 2.} $(Whitney$ $sum$ $formulae)$\\
$c(\xi )=1+c_{1}(\xi )+c_{2}(\xi )+...+c_{n}(\xi )$,\\
$c(\xi\oplus\eta )=c(\xi )c(\eta )$ and\\
$c_{k}(\xi\oplus\eta )=\sum_{i+j=k}c_{i}(\xi )c_{j}(\eta )$.\\

{\bf Lemma 2.} (Splitting Principle) Let $\xi\rightarrow X$ be a complex vector bundle. Then there exists a space
$F(\xi )$ along with a map $f:F(\xi )\rightarrow X$ such that:\\
{\bf i.} $f^{*}:H^{*}(X)\rightarrow H^{*}(F(\xi ))$ is injective and\\
{\bf ii.} $f^{*}(\xi )\simeq L_{1}\oplus L_{2}\oplus ...\oplus L_{n}$, where $n$ is the dimension of $\xi$.\\
(Note that the same holds for K-Theory, namely $f^{*}:K^{*}(X)\rightarrow K^{*}(F(\xi ))$ is also injective).\\

{\bf Proof:} $f_{1}:{\mathbb{P}}(\xi )\rightarrow X$ is injective by Lemma 1 above. Moreover 
$L_{\xi }\subset f^{*}(\xi )\Rightarrow f^{*}(\xi )=L_{\xi}\oplus\tilde{\xi}$ and to complete the proof we use induction. $\square$.\\

Next we would like to construct the Chen character $ch:K^{0}(X)\rightarrow H^{ev}(X; {\mathbb{Q}})$ which satisfies the desired properties (see the Problem above; we ommit the square brackets $[ \xi ]$ to simplify our notation):\\
$\bullet$ $ch (\xi\oplus\eta )=ch(\xi )+ch (\eta )$\\
$\bullet$ $ch (\xi\otimes\eta )=ch(\xi )ch (\eta )$.\\

What should $ch(L)$ be for a line bundle $L$?\\

Suppose $c_{1}(L)=x$ and $c_{1}(\tilde{L})=y$. Then $c_{1}(L\otimes\tilde{L})=x+y$. This is reminiscent of powerseries $F(x)$ (think of $ch(L)=F(x)$) 
where $F(x+y)=F(x)F(y)$, hence 
$$F(x)=e^x.$$ 
We thus end up with the following:\\

{\bf Definition 3.} For a complex line bundle $L$ we define 
$$ch(L)=e^{c_{1}(L)}=1+c_{1}(L)+\frac{c_{1}^{2}(L)}{2}+...$$
For an arbitrary complex vector bundle $E$ we have:
$$E=L_1\oplus L_2\oplus ... \oplus L_n$$ 
with $c_{1}(L_1)=x_1, c_{1}(L_2)=x_2$,...,$c_{1}(L_n)=x_n$. Then
$$ch(E)=e^{x_1}+e^{x_2}+...+e^{x_n}$$
and hence
$$ch(E)=c(L_1)c(L_2)...c(L_n)=(1+x_1)(1+x_2)...(1+x_n)=1+c_1(E)+c_2(E)+...+c_n(E)$$
since $c(L_i)=1+x_i$.\\

Thus
$$c_i(E)=\sigma _i(x_1,x_2,...x_n)$$  
where $\sigma _i$ is the $i-th$ elementary symmetric function.\\

For example $c_1(E)=x_1+x_2+...+x_n$ while $c_2(E)=\prod _{i<j}x_{i}x_{j}$, hence
$$e^{x_1}+e^{x_2}+...+e^{x_n}=n+S_1(\sigma _1,\sigma _2,...,\sigma _n)+S_2(\sigma _1,\sigma _2,...,\sigma _n)+...$$

{\bf Definition 4.}
$$ch(\xi )=\dim\xi +S_{1}(c_{1}(\xi ),c_{2}(\xi ),...,c_n(\xi ))+S_{2}(c_{1}(\xi ),c_{2}(\xi ),...,c_n(\xi ))+ ...$$

As an example let us work out the case of $S_2$:
$$x_1^2+x_2^2+...+x_n^2=(x_1+x_2+...+x_n)^2-2\prod_{i<j}(x_ix_j)=\sigma _1^2-2\sigma _2.$$
Then
$$ch(\xi )=\dim\xi +c_1(\xi )+\frac{1}{2}[c_1^2(\xi )-2c_2(\xi )]+...$$
Recall that the two desired properties $ch (\xi\oplus\eta )=ch(\xi )+ch (\eta )$ and
$ch (\xi\otimes\eta )=ch(\xi )ch (\eta )$ follow from the splitting principle.\\

We defined
$$ch:K^0(X)\rightarrow H^{ev}(X;{\mathbb{Q}}).$$

{\bf Definition 5.}
$$ch([\xi ]-[\eta ])=ch(\xi )-ch(\eta ).$$
Then

$$\begin{CD}
K^{1}X @>>> H^{odd}(X;{\mathbb{Q}})\\
  @|        @|\\
\tilde{K}^0(SX) @>>ch> \tilde{H}^{ev}(SX;{\mathbb{Q}})
\end{CD}$$\\
\\

We finally come to the promised proof of Theorem 1 at the beginning of this lecture:\\

{\bf Proof of Theorem 1:}\\

 Step 1. If $X=S^2$, then $\tilde{K}^{0}(S^2)={\mathbb{Z}}$ and 
$K^{0}(S^2)={\mathbb{Z}}\oplus{\mathbb{Z}}$. Furthermore $ch(1)=1$ and $ch(\xi )=1+u$ because  
$H^{ev}(S^2;{\mathbb{Q}})={\mathbb{Q}}\oplus{\mathbb{Q}}$ (1 corresponds to the first factor $\mathbb{Q}$ and $u$ to the second).\\

Step 2. We have an isomorphism when $X=S^{2n}$ because Chern character comutes with products.\\

Step 3. We have an isomorphism when $X=S^{2n+1}$ (by direct computation).\\

Step 4. Suppose $X$ is a simplicial complex with $N$ simplices and $A\subset X$ is obtained by removing one top dimensional simplex. But then $X/A=S^n$ and we use exact sequence and the 5-Lemma. $\square$.\\

\newpage

\section{Lecture 5 (Operations in K-Theory, Symmetric Products and Adams Operations)}

In the previous section we defined the Chern characters:
$$ch: K^0(X)\rightarrow H^{ev}(X;{\mathbb{Q}})$$
and
$$ch: K^1(X)\rightarrow H^{odd}(X;{\mathbb{Q}})$$
and we stated the basic result that the map 
$$ch:K^{*}(X)\otimes{\mathbb{Q}}\rightarrow H^{*}(X;{\mathbb{Q}})$$
is an isomorphism.\\

Now we shall study \emph{operations in K-Theory}.\\

Let $V$ be a complex vector space of finite dimension say $n$. From linear algebra we know the exterior powers 
$$\wedge ^k V=(V\otimes V\otimes ... \otimes V)/(x_1\otimes x_2 \otimes ... \otimes x_k -sign (\sigma )
x\otimes _{\sigma (1)}x\otimes _{\sigma (2)}...\otimes _{\sigma (k)}x).$$
where we have taken $k$ factors in the tensor product. Clearly $\wedge ^0 V=\mathbb{C}$, $\wedge ^1 V=V,...,\wedge ^k V,...,\wedge ^n V$, where $\dim\wedge ^n V=1$ and more generally
$\dim\wedge ^k V=(n!)/[k!(n-k)!]$.\\

We can extend the above to vector bundles $\xi$ over some space $X$ in a straightforward way: We form the exterior powers
$\wedge ^k \xi$ which are also vector bundles over $X$ with fibre
$$F_x (\wedge ^k \xi )=\wedge ^k (F_x (\xi )).$$
We want to construct a map
$$\wedge ^k :Vect(\xi )\rightarrow K^0X.$$
If $V,W$ are two vector spaces, then
$$\wedge ^k (V\oplus W)=\sum _{i+j=k}\wedge ^i V\otimes\wedge ^j W$$ 
and the corresponding relation for vector bundles will be 
$$\wedge ^k (\xi\oplus\eta )=\sum _{i+j=k}\wedge ^i (\xi )\otimes\wedge ^j (\eta ).$$ 
Observe the analogy with Chern classes:
$$c_k :Vect(X)\rightarrow H^{ev}(X;{\mathbb{Q}})$$
where
$$c_k (\xi\oplus\eta )=\sum_{i+j=k}c_i (\xi )c_j (\eta ).$$
Given a vector bundle $\xi$, we can construct $\wedge _t (\xi )\in K^0(X)[[t]]$ which is a formal power series in $t$ and coefficients from the group $K^0(X)$, in other words
$$\wedge _t \xi =\sum_{k=0}^{\infty}t^k \wedge ^k \xi .$$
Hence
$$\wedge ^k (\xi\oplus\eta )=\sum _{i+j=k}\wedge ^i (\xi )\otimes\wedge ^j (\eta )$$ 
becomes
$$\wedge _t (\xi\oplus\eta )=\wedge _t (\xi )\wedge _t (\eta )$$
which means that we have a map 
$$\wedge _t :Vect(X)\rightarrow K^0(X)[[t]].$$
Our convention is that $\wedge ^0 (\xi )=1\in K^0(X)$
thus $\wedge _t (\xi )$ is a formal power series beginning with $1$. Whence $\wedge _t (\xi )$ is a unit in $K^0(X)[[t]]$.\\

We denote by $G(K^0(X)[[t]])$ the group of units in the ring 
$K^0(X)[[t]]$, namely the group of invertible elements. Thus one gets the following diagram

$$\begin{CD}
Vect(X) @>\wedge _t>> G(K^0(X)[[t]])\\
  @VVV        @AAA\\
K^0(X) @= K^0(X)
\end{CD}$$\\
\\

Given $x\in K^0(X)$, we define $\wedge ^i (x)$ to be the coefficient 
of  $t^i$ in $\wedge ^i_t (x)$. 
If $x=[\xi ]-[\eta ]$, then $\wedge _t (x)=\wedge _t (\xi )\wedge _t (\eta )^{-1}$.\\

For example, if $L$ is a line bundle, then $\wedge _t (L)=1+tL$ and 
$$\wedge _{-t}(-L)=\frac{1}{1-tL}=1+tL+t^2L^2+...$$
Let us mention here that it is difficult to make sense of the value for $t=$something (that's why we are talking about formal power series), the coefficients however make perfect sense (hopefully).\\

Next we recall the symmetric powers from linear algebra
$$S^k(V)=(V\otimes V\otimes ...\otimes V)^{\Sigma _k}$$
where we assume $k$-factors in the tensor product.
Similarly one can form the symmetric powers of vector bundles
$S^k (\xi )$ and then we define
$$S_t (\xi )=\sum_{k=0}^{\infty}t^k S^k (\xi )$$
where our convention is that $S^0 (\xi )=1$. Moreover from linear algebra one has that 
$$S^k(V\oplus W)=\sum_{i+j=k}S^i V\otimes S^j W.$$
By the same method we can get the following commutative diagram

$$\begin{CD}
Vect(X) @>S_t>> G(K^0(X)[[t]])\\
  @VVV        @AAA\\
K^0(X) @= K^0(X)
\end{CD}$$\\
\\

For example, if $L$ is a line bundle, then
$$S_t (L)=1+tL+t^2 L^2+...$$
because
$$S^k L=L^{\otimes k}.$$
Hence
$$S_t (L)=\wedge _{-t}(-L),$$
and thus
$$\wedge _{-t}(L)S_t (L)=1\Leftrightarrow\wedge _t (L)S_{-t}(L)=1\Rightarrow\wedge _t (\xi )S_{-t}(\xi )=1$$
where the arbitrary vector bundle $\xi$ can be written as a sum of line bundles.\\

From the K-Theory splitting principle we deduce that for any vector bundle $\xi$ we have that
$$\wedge _t (\xi )S_{-t}(\xi )=1$$
hence
$$\wedge _{t}([\xi ]-[\eta ])=(\wedge _t \xi )(S_{-t}\eta )$$
and
$$\wedge ^k ([\xi ]-[\eta ])=\sum_{i+j=k}(-1)^j\wedge ^i(\xi )S^j (\eta ).$$
We make an important remark: If $f:X\rightarrow Y$ is a map, then for some $x\in K^0 (Y)$ we have 
$$f^*(\wedge ^i (x))=\wedge ^i (f^* x)$$   
where $\wedge ^i :K^0(X)\rightarrow K^0(X)$.\\

Now we shall study the \emph{Adams operations} in K-Theory.\\

Suppose $x\in K_{0}(X)$ and $\Psi _t (x)\in K^0(X)[[t]]$.\\

We define
$$\Psi _t (x):=\dim (x)-\frac{t}{\wedge _{-t}(x)}\frac{d}{dt}\wedge _{-t}(x)=\dim (x)-t\frac{d}{dt}\log (\wedge _t(x)).$$
All coefficients are integral multiplets of elements in $K^0(X)$.\\

Moreover we define the quantities $\Psi ^k(x)$ via the equation  
$$\sum t^k \Psi ^k (x)=\Psi _t (x).$$
{\bf Definition 1.} The maps
$$\Psi ^k :K^0(X)\rightarrow K^0(X)$$
are called \emph{Adams operations} in K-Theory.\\

{\bf Example:}\\
 
$$\Psi ^1 -\wedge ^1 =0$$
$$\Psi ^2 -\Psi ^1\wedge ^1 +2\wedge ^2 =0$$
$$\Psi ^3 -\Psi ^2\wedge ^1 +\Psi ^1\wedge ^2 -3\wedge ^2=0$$
In general
$$\Psi ^k -\Psi ^{k-1}\wedge ^1 \pm ... \pm k\wedge ^2 =0.$$
Adams operations relate power sums $\Phi ^k (x_1, x_2, ..., x_n)=x_1^k + x_2^k +...+x_n^k$ and elementary symmetric functions $\sigma _i (x_1, x_2, ..., x_n)$, namely
$\Psi ^k\leftrightarrow\Phi ^k$ and $\wedge ^i\leftrightarrow\sigma _i$.\\

\textsl{We summarise the basic properties of Adams operations}:\\

{\bf Proposition 1.} The Adams operations have the following properties:\\

{\bf 1.} 
$$\Psi ^k (x+y)=\Psi ^k (x)+\Psi ^k (y).$$
{\bf 2.} For a line bundle $L$,
$$\Psi ^k (L)=L^k$$ 
{\bf 3.}
$$\Psi ^k (xy)=\Psi ^k (x)\Psi ^k (y)$$ 
{\bf 4.}
$$\Psi ^k (\Psi ^l (x))=\Psi ^{kl} (x).$$ 
{\bf 5.} If $p$ is prime, then 
$$\Psi ^p (x) = x^p +py=x^p mod p$$ 
{\bf 6.} If $u\in \tilde{K}^0 (S^{2n})$, then
$$\Psi ^k (u)=k^n u.$$

{\bf Proof:} Let us prove some of them. For the first we have 

$$\Psi _t (x+y)=\dim (x+y)-\frac{t}{\wedge _{-t}(x+y)}\frac{d}{dt}\wedge _{-t}(x+y)=\dim (x+y)-
\frac{t}{\wedge _{-t}(x)\wedge _{-t}(y)}\times$$
$$\times\frac{d}{dt}(\wedge _t (x)\wedge _{-t}(y))$$
and then we use the Leibniz rule to get the desired result $\Psi _t(x)+\Psi _t (y)$.\\

For the second we have:
$$\wedge _{-t}(L)=1-tL$$
thus
$$\frac{d}{dt}\wedge _{-t}(L)=-L,$$
and consequently
$$\Psi _t (L)=1+\frac{tL}{1-tL}=1+tL+t^2 L^2\Rightarrow \Psi ^k (L)=L^k.$$

For the fourth we briefly have: By the first property (additivity), if the propety holds for line bundles, then by the splitting principle it will also hold for any vector bundle.\\

For the fifth, suppose $x=L_1+L_2+...+L_n$ where the $L_i's$ are line bundles, then
$$\Psi ^p (x)=L_1^p+L_2^p+...+L_n^p=(L_1+L_2+...+L_n)^p mod p.$$

For the last one, let $h$ be a generator of $K^0(S^2)$ and let $J$ denote the Hopf line bundle over $S^2={\mathbb{C}}P^1$. Then
$$h=[J]-1\in \tilde{K}^{0}(S^2).$$
Calculate
$$\Psi ^k(h)=\Psi ^k [J]-\Psi ^k (1)=[J^k]-1$$
hence
$$[J^k]=k[J]-(k-1)\in K^0(S^2).$$
Using the Chern character
$$ch:K^0(S^2)={\mathbb{Z}}\oplus{\mathbb{Z}}\rightarrow H^{ev}(X; {\mathbb{Q}})$$
along with the above relation
$$[J^k]=k[J]-(k-1)\in K^0(S^2)$$
we obtain
$$[J^k]-1=k[J]-k.$$
Then the map
$$\tilde{K}^{0}(S^2)\otimes\tilde{K}^{0}(S^2)\rightarrow\tilde{K}^0 (S^4)$$
where
$$S^4=S^2\wedge S^2.$$
The generator of $\tilde{K}^0(S^{2n})$ is $h\otimes ...\otimes h$ and then
$$\Psi ^k(h\otimes ...\otimes h)=k^n (h\otimes ...\otimes h).$$
$\square$\\

\newpage

\section{Lecture 6 (Applications of Adams operations: Non-Existence of Hopf Invariant $1$ Maps)}

\textsl{As an application of the Adams operations we shall prove the non-existence of Hopf invariant $1$ maps}.\\ 

Let
$$S^{4n-1}\xrightarrow{f} S^{2n}\rightarrow S^{2n}\cup _{f} C^{4n}=X.$$
We calculate that 
$$\tilde{K}^0(X)={\mathbb{Z}}\oplus{\mathbb{Z}}$$ 
using the exact sequence for the pair $(X,S^{2n})$ and $X/S^{2n} =S^{4n}$.\\

Pick some generators $p:X\rightarrow X/S^{2n}=S^{4n}$ and $p^* :\tilde{K}^0(S^{4n})\rightarrow \tilde{K}^0(X)$ where 
$u_{4n}$ is a generator of $\tilde{K}^0(S^{4n})={\mathbb{Z}}$ and $y=p^{*}(u_{4n})\in\tilde{K}^0(X)$.\\

Next we choose an element $x\in\tilde{K}^0(X)$ such that
$$i^*(x)=u_{2n}$$
where
$$i^* :\tilde{K}^0(X)\rightarrow\tilde{K}^0(S^{2n})$$
and $u_{2n}\in\tilde{K}^{0}(S^{2n})$ and we consider 
$$x^2 : i^*(x^2)=0.$$

{\bf Definition.} The \emph{Hopf invariant} of the map $f$ is the integer $\lambda _f$ defined by the equation
$$x^2=\lambda _f y.$$
One can prove that $\lambda _f$ is independent of the choice of $x$.\\

One then has the following result:\\

{\bf Theorem 1.} If $\lambda _f$ is odd, then $n=1,2,4$.\\

{\bf Proof:} Recall that
$$\tilde{K}^0(X)={\mathbb{Z}}\oplus{\mathbb{Z}}$$
with generators, say $x$ and $y$ (one for each copy of $\mathbb{Z}$). Then:
$$\Psi ^2(x)=2^nx+ay$$
and
$$\Psi ^3(x)=3^nx+by$$
while
$$\Psi ^k(y)=\Psi ^k (p^* (u_{4n}))=p^* ( \Psi ^k (u_{4n}))=p^* (k^{2n}u_{4n})=k^{2n}y,$$
where
$$\Psi ^2 (x)\mapsto\Psi ^2 (u_{2n})=2^n u_{2n}$$
and
$$\Psi ^3 (x)\mapsto\Psi ^3 (u_{2n})=3^n u_{2n}.$$
Moreover we compute:
$$\Psi ^6(x)=\Psi ^3 (\Psi ^2 (x))=6^nx+(2^n b + 3^{2n}a)y$$
$$\Psi ^6(x)=\Psi ^2 (\Psi ^3 (x))=6^nx+(2^{2n} b + 3^{n}a)y,$$
thus
$$2^n b+3^{2n}a=2^{2n}b+3^n a\Leftrightarrow 2^n (2^n -1)b=3^n (3^n -1)a.$$
But\\

$\Psi ^2 (x)=x^2$ $mod$ $2=\lambda _f y$ $mod$ $2 = y$ $mod$ $2\Rightarrow a$ $odd.$\\

Whence from
$$2^n (2^n -1)b=3^n (3^n -1)a$$
we deduce that $2^n$ divides $3^n -1$. By direct calculation we obtain taht $n=1,2,4$. $\square$.\\

There are three classical Hopf maps:\\

$S^3\rightarrow S^2$, $S({\mathbb{C}}^2)=S^3\rightarrow S^2 ={\mathbb{C}}P^1$, $(z_1, z_2)\mapsto [z_1, z_2]$\\

$S^7\rightarrow S^4$, $S({\mathbb{H}}^2)=S^7\rightarrow S^4 ={\mathbb{H}}P^1$, $(q_1, q_2)\mapsto [q_1, q_2]$\\

$S^{15}\rightarrow S^8$, $S({\mathbb{O}}^2)=S^{15}\rightarrow S^8 ={\mathbb{O}}P^1$, $(k_1, k_2)\mapsto [k_1, k_2]$\\

where $\mathbb{H}$ denotes the quaternions while $\mathbb{O}$ denotes the Caley numbers (the octonions).
\\

\textsl{It follows from Theorem 1 above that there no other such maps}.\\

An analogous question is this: \textsl{In the sequence of numbers}\\

\emph{reals} $\rightarrow$ \emph{complex} $\rightarrow$ \emph{quaternions} $\rightarrow$ 
\emph{Caley numbers}\\

\textsl{are there any other numbers?}\\

\textsl{The answer is also negative}.\\

We end this lecture with the following remark: The circle $S^1$ is (topologically) the group $(U(1)$.\\

The sphere $S^3$ is (topologically) the group of unit quaternions.\\

The sphere $S^7$ is topologically the group of unit Caley numbers.\\

The space $1\times S^7 \subset S^7 \times S^7\rightarrow S^7$ is also a "group", i.e. it is an h-space whereas
$S^7 \times 1\subset S^7 \times S^7\rightarrow S^7$, in other words we want to see if the spheres $S^i$ are \emph{parallelisable spaces}, namely if their tangent bundles are trivial.\\

\newpage

\section{Analytic K-Theory (K-Homology)}

One motivation for analytic K-Theory could be the question \emph{"what is infinity minus infinity"?} Moreover the 
\textsl{index} of a 
\textsl{Fredholm operator} is one of the most useful definitions in mathematics since it gave rise to various 
\emph{index theorems}, 
arguably the most central result in mathematics during the second half of the 20th century. At the same time, the K-Theory of $C^*$-algebras is an important ingredient in A. Connes' \emph{noncommutative geometry}.\\

\section{Lecture 1 (Some preliminaries from Functional Analysis, the Index of Fredholm Operators)}

From topological K-Theory we know that $K_0(*)=\mathbb{Z}$, namely K-Theory does not satisfy the dimension axiom, hence 
it is a generalised (or exotic) homology theory.\\

Let $V_0,V_1$ be two finite dimensional vector spaces (over some field 
 ${\mathbb{F}}$, usually the complex numbers). It is straightforward to compute the difference of their dimensions 
which we shall denote $[V_0]-[V_1]$. (For simplicity we shall often write $[V]$ instead of $\dim V$).\\

Let $T:V_0\rightarrow V_1$ a linear operator. We expand:
$$0\rightarrow KerT \rightarrow V_0 \xrightarrow{T} V_1 \rightarrow cokerT\rightarrow 0,$$
where
$$cokerT=V_{1}/ImT=V_{1}/T(V_{0})=(ImT)^{\perp}.$$

{\bf Theorem 1} (Linear Algebra).
$$[KerT]-[V_{0}]+[V_{1}]-[cokerT]=0.$$
(The theorem holds for any short exact sequence of finite dimensional vector spaces).\\
 
We can rewrite the relation of the theorem in a more convenient form:
$$[V_{0}]-[V_{1}]=[KerT]-[cokerT].$$

If we try to generalise the above relation to \emph{infinite} dimensional vector spaces (for example Hilbert spaces), 
we see that the LHS has no meaning since it gives $\infty -\infty$, but the RHS may, in some favourable cases, 
give a meaningful result (provided that the operator $Ô$ is a "nice" operator). Hence the difference  
$\infty -\infty$ may give, in some cases, a finite result.\\

In this chapter, unless otherwise stated explicitly, all operators are assumed to be linear.\\

Let us briefly recall the definition of a Hilbert space:\\

{\bf Definition 1.} A Hilbert space $H$ is an infinite dimensional complete (complex) vector space with an inner product. 
Complete means that every sequence $(x_n)\in H$ with $\sum ||x_n||<\infty$, converges in $H$.\\

{\bf Definition 2.} A linear map (linear opearator) $T:H_0\rightarrow H_1$ is called \emph{bounded} if there exists 
some positive real number $c$ such that
$$||Tx||\leq c||x||, \forall x\in H.$$
The "best" such positive real number is called the \emph{norm} of $T$ and it is denoted $||T||$. In particular, one can define
$$||T||=sup <Tx,y>$$
where $||x||=||y||=1$.\\

The set of all bounded linear operators $T:H\rightarrow H$ will be denoted $B(H)$.\\

Given a Hilbert space, we can define the notion of a \emph{Hilbert basis} and with respect to some Hilbert bases in the 
Hilbert spaces $H_0, H_1$ above, every 
(bounded) linear operator $Ô$ can be represented by an $\infty\times\infty$ matrix (basically using the same recipe 
as for representing linear maps between finite dimensional vector spaces by matrices). Conversely, given an 
$\infty\times\infty$ matrix, it is a lot more difficult to tell whether it defines a bounded linear operator.\\

We know however that an $\infty\times\infty$ matrix is indeed the representation of some bounded linear operator if (but not 
only if) it has 
a finite number of non-zero elements. In this case it is called a \emph{finite rank operator}. In other words, finite rank 
operators are bounded (but the converse does not necessarily hold).\\

{\bf Definition 3.} A linear operator $T$ is called \emph{compact} if $\forall\epsilon >0$, there exists 
a finite rank operator $F$ such that
$$||T-F||<\epsilon.$$

The set of all compact operators $T:H\rightarrow H$ will be denoted $K(H)$.\\

$B(H)$ has a natural $C^{*}$-algebra structure and $K(H)$ is an ideal in $B(H)$ (closed, 2-sided $*$-ideal), thus one can form the quotient
$$Q(H)=B(H)/K(H)$$
which is another $C^{*}$-algebra, called the \emph{Calkin algebra} (see also Lecture 5 later).\\

For example, the identity operator $I:H\rightarrow H$ is {\bf not} compact.\\

{\bf Proposition 1.} Let $B(H)$ denote the set of all bounded linear operators from $H$ to itself. 
The subset of all invertible operators in $B(H)$ is open.\\

{\bf Proof:} It saffices to prove that the identity operator $I$ is an internal point of the set of invertibles, 
because if $||A||<1$, then $(I-A)^{-1}=I+A+A^2 +...$. $\square$.\\

{\bf Definition 4.} A bounded linear operator $T:H_0\rightarrow H_1$ between two Hilbert spaces 
 is called \emph{Fredholm} if its kernel $KerT$ is of finite dimension, its image $ImT$ is a closed subspace 
of $H_1$ and the quotient $H_1/ImT$ has also finite dimension.\\

One has the following fundamental result:\\

{\bf Theorem 1.} (Atkinson) The following are equivalent:\\

{\bf a.} $T$ is Fredholm.\\
  
{\bf b.} $T$ is invertible modulo compacts, namely there exists an operator 
 $S:H_1\rightarrow H_0$ such that $ST-I$ and $TS-I$ are compact on $H_0$ and $H_1$ respectively. Alternatively, the term invertible modulo compacts means that the image $\pi (T)$ of $T$ under the canonical projection $\pi :B(H)\rightarrow Q(H)$ is invertible.\\

{\bf Proof:}\\
{\bf a. $\Rightarrow$ b.}\\
$H_0$ is the disjoint union of the subspaces $(KerT)^{\perp}$ and $KerT$ whereas $H_1$ is the disjoint union of 
$ImT$ and $(ImT)^{\perp}\simeq H_1/ImT$, hence there is a linear map $(KerT)^{\perp}\rightarrow ImT$. Appeal to 
the closed graph theorem (see Appendix) of functional analysis and deduce that the above continuous 
bijection has a continuous inverse $S:(ImT)\rightarrow (ImT)^{\perp}$. Then $ST-I$ and $TS-I$ are finite rank 
operators.\\

{\bf b. $\Rightarrow$ a.}\\
Suppose that we have an operator $S$ such that both $ST-I$ and $TS-I$ are compact. 
We restrict the first to $KerT$ and then
$$(ST-I)|_{KerT}=-I|_{KerT},$$
so $I|_{KerT}$ is compact, thus $KerT$ is finite dimensional.\\

Similarly the other operator $TS-I$ will give that $cokerT$ is also finite dimensional provided that we know that 
$ImT$ is closed. Then
$$||ST-I-F||<1/2,$$
where $F$ is of finite rank. Then if $x\in KerF$,
$$||STx-x||\leq (1/2)||x||,$$ 
so
$$||STx||\geq (1/2)||x||,$$
hence
$$||Tx||\geq \frac{1}{2||S||}x,$$
namely $T$ is bounded below on $KerF$, thus $T(KerF)$ is closed. But
$$ImT=T(KerF)+T((KerF)^{\perp})$$
is closed. $\square$.\\

We denote $Fred(H)$ the set of all Fredholm operators from $H$ to itself.\\

{\bf Definition 5.} If $Ô\in Fred(H)$, then we define the \emph{Index} of $T$, denoted $IndexT$ 
(or  $IndT$ for short), as the difference
$$IndT:=dim(KerT)-dim(cokerT)=[KerT]-[cokerT]\in\mathbb{Z}.$$
Clearly the index of a Fredholm operator is an integer.\\

\textsl{Intuitively one can say that the index measures how far an operator is from being invertible since for invertible 
operators the index vanishes}.\\

A basic property of the index is that if $T_0,T_1\in Fred(H)$, then $T_0T_1 \in Fred(H)$ and
$$Ind(T_0T_1)=IndT_0 + IndT_1.$$ 
The proof of the above property can be deduced by the snake (or serpent) Lemma in Homological Algebra. The following sequence is exact:\\

$$0\rightarrow KerT_{1}\rightarrow KerT_{0}T_{1}\rightarrow KerT_{0}\rightarrow cokerT_{1}\rightarrow cokerT_{0}T_{1}\rightarrow cokerT_{0}\rightarrow 0.$$

{\bf Example:} Let $H$ have a basis $e_0,e_1,e_2,...$ and suppose $U:H\rightarrow H$ is defined by:
$$Ue_0=e_1$$
$$Ue_1=e_2$$
 $$ etc $$
Then $KerU=\{0\}$ and $dim(cokerU)=1$. Hence $IndU=-1$.\\

From the above example we deduce that the index depends much on the spaces $H_0, H_1$ and little on the operator 
$Ô$ between them.\\

{\bf Proposition 2.} If $T_0, T_1\in Fred(H)$ can be linked by a continuous path $T_t$, $t\in [0,1]$ 
of Fredholm operators, then 
$$IndT_0 = IndT_1.$$
{\bf Proof:} We need only prove that the map $t\mapsto IndT_t$ is locally constant. We look near $t=t_0$.\\
{\bf a.} Perhaps $T_{t_0}$ is bijective (hence invertible). Then nearby $T_t$'s are invertible as well, so the index 
is constant (zero in fact) near $t=t_0$.\\
{\bf b.} Perhaps $T_{t_0}$ is surjective, suppose $dim(KerT_{t_0})=n$. Consider
$$S_t :H\rightarrow H\oplus{\mathbb{C}}^n,$$ 
$$x\mapsto (T_{t}x,P_{KerT_{t_0}}x).$$ 
These operators are Fredholm and  
$$IndS_{t}=IndT_{t}-n.$$

But $S_{t_0}$ is bijective, so the index of $S_t$ is constant near $t_0$.\\
{\bf c.} Perhaps neither of the above holds. Then we choose an orthonormal basis $e_i$ and let $Q_n$ be the 
projection onto $H_n$, the subspace of $H$ generated by the $e_{n+1},e_{n+2},...$ We claim that
$$Q_n T_{t_0}:H\rightarrow H_n$$
is surjective for large enough $n$. Hence
$$Ind(Q_n T_{t_0})=IndT_{t_0}+m$$
and apply case (b.) above. $\square$.\\

Proposition 2 in other words states that the set of Fredholm operators is an open subset in $B(H)$ and the index is locally constant in $Fred(H)$.\\ 

{\bf Corollary 1.} If $T$ is Fredholm and $K$ compact, then $T+K$ is Fredholm and moreover
$$Ind(T+K)=IndT.$$
{\bf Proof:} Consider $T_t=T+tK$ and apply Proposition 2. $\square$.\\

\newpage

\section{Lecture 2 (Index of Toeplitz Operators, Winding Number and K-Homology)}

We know that
$$Fred(H)\simeq{\mathbb{Z}}\times G_{\infty}.$$
If  $M$ is a compact, oriented, Riemannian manifold, then we can form
$$L^2(M)=\{f:M\rightarrow{\mathbb{C}}| \int_{M}|f|^2 d\mu <\infty \}.$$
This is a Hilbert space with inner product
$$<f,g>=\int_{M}f\bar{g}d\mu .$$

For example, let $L^2 (S^1)$ be the set of square integrable complex functions on the circle $S^1$. We think of the circle 
as $[0,2\pi ]/(0=2\pi )$, namely the closed interval $[0,2\pi ]$ with the end points identified.
Then think of the Fourier series: A complete orthonormal basis of $L^2 (S^1)$ is given by the functions 
$$e_n (x)=\frac{1}{\sqrt{2\pi}}e^{inx},n\in\mathbb{Z},$$
i.e. if $f\in L^2 (S^1)$, then 
$$f(x)=\sum c_n e^{inx},$$ 
where
$$c_n =\frac{1}{2\pi}\int _{0}^{2\pi}f(x)e^{-inx}dx.$$
We can also regard the circle as a subset of the complex numbers, $S^1\subset\mathbb{C}$. Then the basis become 
$$e_n=\frac{1}{\sqrt{2\pi}}z^n.$$
Complex analysis suggests considering the subspace $H^2 (S^1)=\{f\in L^2 (S^1)|f=\sum_{n=0}^{\infty}c_n z^n\}$, 
\emph{the Hardy space}, namely the functions with only positive Fourier coefficients. The orthogonal projection 
$$P:L^2 (S^1)\rightarrow H^2 (S^1)$$
is called the \textsl{Hardy projection}.\\

Notice that $L^2 (S^1)$ is a module over $C(S^1)$, (the space of continuous complex valued functions 
on the circle), namely given an $f\in C(S^1)$, we can define a \textsl{multiplication operator} 
$M_f:L^2 (S^1)\rightarrow L^2 (S^1)$ via the relation:
$$(M_f g)(x)=f(x)g(x).$$  
This is a \emph{bounded} operator with
$$||M_f||=sup|f|.$$
\textsl{The non-locality of operators can be measured by the commutator of it with multiplications}.\\

{\bf Lemma 1.} If $f\in C(S^1)$, then the commutator
$$[P,M_f]=PM_f - M_f P$$
is a compact operator (hence the Hardy projection is non-local but not very far from being local).\\

{\bf Proof:} Suppose first that $f$ is a polynomial in $z$, $z^{-1}$, say
$$f=a_{-m} z^{-m}+...+a_m z^m.$$
Then, relative to the basis $\{e_n\}$ mentioned above, $f$ can be represented by a matrix. By commuting with $M_f$, 
we get the Lemma (we get a finite rank operator, and hence compact). So if $f$ is polynomial, then the commutator 
$[P,M_f]$ is a finite rank operator and hence by the \textsl{Stone-Weierstrass theorem} (see Appendix) we 
deduce that polynomials are dense  in $C(S^1)$. By approximation we get the result in full generality. $\square$.\\

{\bf Definition 1.} Let $f\in C(S^1)$. Then the ï \emph{Toeplitz operator} 
$T_f :H^2 (S^1)\rightarrow H^2 (S^1)$ is the operator $T_f =PM_f$. Then $f$ is called the \emph{symbol} of $T_f$.\\

{\bf Proposition 1.} $T_f T_g -T_{fg}$ is compact.\\

{\bf Proof:} $PM_f PM_g -P^2 M_f M_g =P[M_f ,P]M_g $ which is compact. $\square$.\\

{\bf Corollary 1.} If $f\in C(S^1)$ is invertible, then $T_f$ is Fredholm.\\

{\bf Proof:} $T_{f^{-1}}$ is an inverse modulo compacts. $\square$.\\

{\bf Key Question:} \textsl{What is the index of} $T_f$?\\

Let $f:S^1\rightarrow {\mathbb{C}}^{x}={\mathbb{C}}-\{0\}$ be a map such that the curve (graph of $f$) never hits the origin. 
Recall that the \emph{winding number}, denoted $wn(f)$, of $f$ is the only topological invariant of such maps.\\

{\bf Theorem 1.} $$IndT_f = -wn(f).$$ 

{\bf Proof:}  By deformation invariance of the index, the correspondence 
$$f\mapsto IndT_f$$
defines a map
$$[S^1 ,{\mathbb{C}}^{x}]={\mathbb{Z}}\rightarrow{\mathbb{Z}}.$$
Each homotopy class contains a representative $z^m$ and we just check what the answer is on this representative.\\

For example, consider $f(x)=z$. Then $H^2(S^1)$ is spanned by $e_0,e_1,e_2,...$ whereas $M_z (e_i)=e_{i+1}$. Then 
$T_z =T_f$ is just $T_z(e_i)=e_{i+1},i=0,1,2,...$. Thus $T_z$ is the unilateral shift with index $-1$. 
$\square$.\\

Next, we consider instead operators $T_f$ where $f:S^1 \rightarrow M_n({\mathbb{C}})$, then such operators are 
Fredholm provided that $f$ is invertible, namely if $f:S^1 \rightarrow GL_n({\mathbb{C}})$. Just as before, we get maps 
$[S^1, GL_n ({\mathbb{C}}]\rightarrow\mathbb{Z}$ with $f\mapsto IndT_f$. But 
$$GL_n({\mathbb{C}})\xrightarrow{i} GL_{n+1}({\mathbb{C}}),$$
where the inclusion $i$ given explicitly by
$$f\mapsto 
\left( \begin{array}{cc}f & 0\\
0 & 1\end{array} \right)=\tilde{f}.$$
Then
$$T_{\tilde{f}}=\left( \begin{array}{cc}T_f & 0\\
0 & 1\end{array} \right).$$
Hence the following diagram commutes:

$$\begin{CD}
[S^1, GL_{n} ({\mathbb{C}}] @>ind>> \mathbb{Z}\\
 @ViVV     @|\\
[S^1, GL_{n+1} ({\mathbb{C}}] @>ind>> \mathbb{Z} 
\end{CD}$$\\
\\

In the limit $n\rightarrow\infty$ we get an \textsl{index} map (in fact a group homomorphism) 
$$K^{-1}(S^1)\rightarrow\mathbb{Z}.$$

The above can be generalised for any compact space $X$ (instead of the circle $S^1$) equipped with the following data:\\ 

{\bf a.} A Hilbert space $H$ which is a module over $C(X)$ (namely equipped with an algebra homomorphism 
$C(X)\rightarrow B(H)$ with $f\mapsto Mf$).\\

{\bf b.} A self-adjoint projection $P$ on $H$.\\

{\bf c.} $PM_f -M_f P$ compact $\forall f\in C(X)$.\\

The above data define an \emph{odd K-cycle} for the space $X$, namely an element of $K^{-1}(X)$.\\

Associated to such K-cycles there is a group homomorphism
$$K^{-1}(X)\rightarrow\mathbb{Z}.$$
This suggests  that K-cycles should generate a "dual" to K-Theory, namely \emph{K-Homology}.\\

{\bf Note:} By considering an arbitrary (possibly noncommutative) algebra $A$ instead of the commutative algebra 
$C(X)$ of some space $X$, one gets the definition of K-cycle a la Connes in Noncommutative Geometry.\\

\newpage

\section{Lecture 3 (The slant product and the pairing between K-Theory and K-Homology)}

{\bf Definition 1.} Let $X$ be a compact topological space. An \emph{odd K-cycle} of $X$ (odd because we deal with 
$K_1(X)$ as we shall see shortly) is defined by the following data:\\

{\bf a.} A Hilbert space $H$ on which $C(X)$ acts.\\

{\bf b.} A projection $P$ on $H$ (projection means that $P=P^2=P^*$).\\

{\bf c.} $\forall f\in C(X), [P,f]$ is compact.\\

Odd K-cycles thus give a homomorphism 
$$K^{-1}(X)\rightarrow\mathbb{Z}$$
(which is given by the index of the corresponding Toeplitz operator).\\

{\bf Definition 2.} A Ê-cycle is \emph{degenerate} if $[P,f]=0, \forall f$ which induces the zero map 
$K^{-1}(X)\rightarrow\mathbb{Z}$.\\

K-cycles have the following properties:\\

$\bullet$ K-cycles can be added by direct sum.\\

$\bullet$ There is a notion of \emph{equivalence} for K-cycles: \textsl{Unitary equivalence}, namely 
a unitary operator $U:H\rightarrow H'$ and \textsl{homotopy equivalence}, namely the existence of a homotopy family 
$(H,P_t), t\in [0,1]$ of K-cycles with $t\mapsto P_t$.\\

\textsl{For a given $X$, the set of equivalence classes of K-cycles forms a semi-group under direct sum
and we can take the Grothendieck group of this semi-group denoted $K_1 (X)$ (which is an Abelian group)}.\\

Thus we get a \emph{pairing}
$$K_1 (X)\otimes K^{-1}(X)\rightarrow\mathbb{Z}.$$
 
{\bf Lemma 1.} Degenerat K-cycles are zero in $K_1 (X)$.\\

{\bf Proof:}  Let $(H,P)$ be degenerate. Form $(\oplus ^{\infty}H,\oplus ^{\infty}P)$ which is a K-cycle.\\

Let $X=[(H,P)]$ in $K_1 (X)$ and $Y=[(H^{\infty},P^{\infty}]$. Then $X+Y=Y\Rightarrow X=0$. $\square$.\\

We should have the analogue of the \emph{slant product} 
$$K_1 (X)\otimes K^{-1}(X\times Y)\rightarrow K^{0}(Y).$$
Elements of $K_1(X)$ are odd K-cycles of the form $(P,H)$ where $C(X)$ acts on $Ç$. Elements of 
$K^{-1}(X\times Y)$ are maps 
$X\times Y\rightarrow GL_n ({\mathbb{C}})=Y\rightarrow Maps (X\rightarrow GL_n ({\mathbb{C}}))$
namely $y\mapsto g_y\mapsto T_{g_y}\in Fred(P,H)$, hence we get a map $Y\rightarrow Fred(H)$, where $Fred(H)$ is the set 
of Fredholm operators on $H$.\\

The key point in the construction of slant product is to show how a Fredholm family over $Y$, namely a continuous map 
$Y\rightarrow Fred(H)$ gives rise to a K-theory class in $K^0(Y)$. The requirement is that when $Y=*$ a point, 
this should be the \emph{index}.\\

The idea is to consider the maps $\{KerT_y :y\in Y\}\rightarrow Y$ and 
$\{cokerT_y :y\in Y\}\rightarrow Y$. We want the index $IndT_y$ to be the formal difference of vector bundles  
$[KerT_y]-[cokerT_y]$. Unfortunately these might not be vector bundles since the dimensions may jump.\\ 

Recall however from previous work on deformation invariance of the index that we can find a finite codimesnion projection 
$Q_m$ such that $Q_m T_y$ is onto the range of $Q_m$ for any $y\in Y$.\\

Now the kernels $Ker(Q_m T_y)$ do form a vector bundle. Near $y_0$, the operator $Q_m T_{y_0}$ is onto the range 
of $Q_m$. Consider the map
$$S_y :H\rightarrow (Range Q_m)\oplus (Ker(Q_m T_{y_0})).$$
The map $S_{y_0}$ is invertible, hence $S_y$ is invertible $\forall y\in U$, where $U$ is a neighborhood of  
$y_0$. Consider next
$$S^{-1}_{y}(O\oplus Ker(Q_m T_{y_0}))=KerS_y.$$
Hence $S_y$ restricted to $U$, gives a bijection between the family of vector spaces 
$KerS_{y}\rightarrow U$ and the constant family $O\oplus Ker(Q_m T_{y_0})\rightarrow U$. Hence $KerS_y$ forms 
a vector bundle.\\

To summarise then, we have "defined" a map
$$[Y,Fred(H)]\rightarrow K^0 (Y)=[Y,{\mathbb{Z}}\times G_{\infty}].$$

{\bf Theorem 1.} The map
$$Fred(H)\rightarrow{\mathbb{Z}}\times G_{\infty}$$
is a weak homotopy equivalence.\\

{\bf Proof:}  The result follows from the above construction $\square$.\\

Consider the space
$$Fred^{S}(H)=\{ \left( \begin{array}{cc}T & 0\\
0 & 1\end{array} \right)| T\in Fred(H)\}.$$

{\bf Theorem 2.} The map
$$[X,Fred^{S}(H)]\rightarrow K^{0}(Y)$$
is an isomorphism. This map is nothing other than the index.\\

{\bf Proof:}  We have to prove that the map is surjective and injective.\\

For surjectivity: Given $V,W$ two vector bundles over $Y$, we embed them in a big enough trivial bundle (because the Hilbert 
space is infinite dimensional):
$$V\hookrightarrow{\mathbb{C}}^N\hookrightarrow H$$
and
$$W\hookrightarrow{\mathbb{C}}^N\hookrightarrow H.$$
The Hilbert space $H$ can be written as the direct sum of subspaces $V_y\oplus V_y^{\perp}$ and $W_y\oplus W_y^{\perp}$ 
where $T_y$ 
is the linear map between them.\\

For injectivity: Let $(T_y)$ be a family of index zero. Then using a "cutoff" as the $Q_m$ before, we can deform it to a 
family with constant cokernel of dimension $m$. Since the index is zero, $(KerT_y)$ is a stably trivial vector bundle of 
dimension $m$. By taking $m$ large enough, we may assume that $(KerT_y)$ is a trivial vector bundle. The Hilbert space can 
be written as the direct sum $(KerT_y)\oplus (KerT_y)^{\perp}$ and ${\mathbb{C}}^m\oplus S$, where $S$ 
some fixed subspace. The linear map between them is $U_y+T_y$ which is invertible. Hence $T_y$ can be deformed 
to a family of invertibles.\\

Next recall the following fact: The group $GLH$ is contractible (this is difficult to prove). 
For our stable version we need only know that the space
$$\{\left( \begin{array}{cc}U & 0\\
0 & 1\end{array} \right)| U\in GLH\}$$
is contratible. To complete the proof we recall that from the Whitehead Lemma 
$$\left( \begin{array}{cc}U & 0\\
0 & U^{-1}\end{array} \right)\sim 
\left( \begin{array}{cc}1 & 0\\
0 & 1\end{array} \right),$$
(where $\sim$ means homotopic), hence the following classes can be homotoped 
$$U_y\oplus 1\oplus 1\oplus ...$$
êáé
$$U_y\oplus U_y^{-1}\oplus U_y\oplus ...$$
$\square$.\\

\newpage

\section{Lecture 4 (Bott Periodicity)}

Recall that in the previous section we defined the slant product:

$$K_{1}X\otimes K^{-1}(X\times Y)\rightarrow K^{0}Y,$$
where
$$x\otimes y\mapsto x\backslash y.$$

\emph{Claim}:
$$K_{1}X\otimes K^{-1}(X\times Y)\otimes K^{0}Z\rightarrow K^{0}(Y\times Z),$$
with
$$x\otimes y\otimes z\mapsto (x\backslash y)\otimes z=x\backslash (y\otimes z).$$
This is obvious for:\\ 
What is the pairing $K^{-1}W\otimes K^{0}Z\rightarrow K^{-1}(W\times Z)$?\\
The answer comes from the fact that an automorphism of any bundle $V$ over $W$ 
defines an element $V\oplus V'$ of the group $K^1 (W)$.\\ 

We introduce the idea of K-Theory for locally compact spaces (locally compact means that each point has a compact 
neighbourhood).\\

Let $X$ be a locally compact topological space and let $X^{+}$ denote its 1-point compactification, namely
$X^{+}=X\cup \{\infty \}$.\\

{\bf Definition 1.}
$$K^{i}X=\tilde{K}^{i}X^{+}:=Ker(K^{i}X^{+}\rightarrow K(*)).$$
This definition makes sense because by functoriality, a map $X\rightarrow Y$ 
induces a map $X^{+}\rightarrow Y^{+}$.\\

Thus Bott periodicity becomes a theorem relating $K^{*}X$ with $K^{*}({\mathbb{R}}\times X)$.\\

Recall the elementary fact:
$$K^{0}({\mathbb{R}}\times X)=K^{-1}X.$$
Moreover
$$({\mathbb{R}}\times X)^{+}=SX/A$$
where $SX$ is the \textsl{suspension} of $X$ and $A$ denotes the two vertices of the suspension $SX$.\\

{\bf Theorem 1.} 
$$K^{-1}({\mathbb{R}}\times X)\simeq K^{0}X.$$

{\bf Proof:}  Let's denote the above isomorphism by $\beta$ and let's try to describe it: Let 
$b\in K^{-1}\mathbb{R}$ be the generator corresponding to the function
$z\mapsto 1/z$ on $S^1$. Then
$$\beta :K^{0}X\rightarrow K^{-1}({\mathbb{R}}\times X),$$
with
$$x\mapsto b\otimes x.$$
Note that
$$K^{0}(*)\xrightarrow{\beta} K^{-1}{\mathbb{R}}\rightarrow K^{0}({\mathbb{R}}^{2})=\tilde{K}^{0}(S^2)$$
where $\beta (1)$ is the tautological line bundle.\\

In the case where $S^2={\mathbb{C}}P^1$, we choose homogeneous coordinates $(z_0,z_1)$, the North Hemisphere consists of 
points with homogeneous coordinates $(1,z_1)$ whereas the South Hemisphere consists of points with homogeneous coordinates
$(z_0,1)$ and the boundary circle (the equator) consists of points with coordinates $(z_0/z_1)$ and thus we obtain the 
Hopf line bundle.\\
Next
$$\beta :K^{0}X\rightarrow K^{-1}({\mathbb{R}}\times X)=K^{0}({\mathbb{R}}^{2}\times X),$$
with
$$x\mapsto b\otimes x.$$
We define the map 
$$\alpha :K^{-1}({\mathbb{R}}\times X)\rightarrow K^{0}X.$$
We also have maps
$$K^{-1}({\mathbb{R}}\times X)\rightarrow K^{-1}(S^1\times X)$$ 
and
$$a\backslash :K^{-1}(S^1\times X)\rightarrow K^{0}X,$$
where $a\in K_1 (S^1)$, hence $\alpha (y)=a\backslash y$. To continue, we shall need the following\\ 

{\bf Lemma 1.} 
$$\alpha (\beta (x))=x, \forall x\in K^{0}X.$$
{\bf Proof of Lemma:} $a\backslash (b\otimes x)=(a\backslash b)\otimes x=1\otimes x =x$. $\square$.\\

To complete the proof of the theorem that $\beta$ is an isomorphism, it will be enough to prove that it is surjective, 
namely we would like to know that every $y\in K^{0}({\mathbb{R}}^{2}\times X)$ is of the form $b\otimes (-)$.\\
  
Consider the map
$$K^{0}({\mathbb{R}}^{2}\times X)\rightarrow K^{0}({\mathbb{R}}^{2}\times X\times{\mathbb{R}}^2),$$
with
$$y\mapsto y\otimes b.$$
Notice that $y=(y\otimes b)\backslash a$ by Lemma 1. Yet $y\otimes b=b\otimes\tilde{y}$, where $\tilde{y}$ is the image 
of $y$ under the map
$${\mathbb{R}}^2\times X\rightarrow X\times{\mathbb{R}}^2,$$
with
$$(v,p)\mapsto (p,-v).$$
The matrix 
$$\left( \begin{array}{cc}0 & -I\\
I & 0\end{array} \right)$$
which acts on ${\mathbb{R}}^2\oplus{\mathbb{R}}^2$ is homotopic to the identity, thus we get 
$$y=(y\otimes b)\backslash a=(b\otimes\tilde{y})\backslash a=b\otimes (\tilde{y}\backslash a).$$
$\square$.\\

The above is known as \emph{Atiyah's trick}.\\

We can generalise: $X\times{\mathbb{R}}^{2k}$ is the total space of a trivial complex vector bundle over $X$. Let 
$V$ be any complex vector bundle over $X$. \textsl{Can we generalise Bott periodicity to compute $KV$ in terms of $KX$?}\\ 

The answer is affirmative.\\

There is a map
$$\beta :KX\rightarrow KV$$ 
which is an isomorphism, $\beta$ is a product with a \emph{Thom class} for $V$, namely a class of $KV$ 
which restricts to each fibre to the usual Bott generator.\\ 

Suppose $V=L$ is a complex line bundle. The Thom class is the tautological line bundle over 
${\mathbb{P}}(L\oplus 1)$. The proof of the Thom isomorphism is very similar to Atiyah's trick.
See the following reference: M.F. Atiyah: \emph{"Bott Periodicity and the Index of Elliptic Operators"}, 
Q.J. of Maths, Oxford, May 1968.\\

\newpage

\section{Lecture 5 (K-Theory of Banach and $C^*$-algebras)}

Recall that a \emph{Banach Algebra} is a complete normed algebra.\\

A $C^*$-algebra is an involutive Banach algebra such that
$$||xx^*||=||x||^2.$$

{\bf Examples:}\\

{\bf a.} The set $B(H)$ of \textsl{bounded operators} on a Hilbert space $H$.\\

{\bf b.} The set $K(H)$ of \textsl{compact operators} on a Hilbert space $H$ (this is an ideal of $B(H)$).\\

{\bf c.} The \textsl{Calkin algebra}
$$Q(H)=B(H)/K(H).$$

{\bf d.} If $X$ is a compact Hausdorff space, then $C(X)$ (the set of continuous complex functions on $X$) 
equipped with the supreme norm
$$||f||=sup|f|$$
and involution
$$g^{*}=\bar{g}$$
where $\bar{g}$ denotes complex conjugate, is a $C^*$-algebra.\\

{\bf Aside Note 1:} If $X$ is only locally compact and Hausdorff, then consider $C_{0}(X)$ (the set of complex continuous 
functions vanishing at infinity). This is also a $C^*$-algebra \emph{without unit}. Moreover
$$C_0(X)=Ker(C(X^{+})\rightarrow{\mathbb{C}}),$$ 
with
$$f\mapsto f(\infty ),$$
where $X^{+}$ is the 1-point compactification of $X$.\\

Every non-unital $C^*$-algebra can be "1-compactified" in this way: 
$$J\rightarrow J^{+}\rightarrow{\mathbb{C}},$$
where
$$J^{+}=\{j+\lambda 1:\lambda\in{\mathbb{C}}\}.$$
{\bf Aside Note 2:} Because from Gelfand's theorem we know that the topology of a compact Hausdorff space $X$
is captured by $C(X)$ which is a commutative $C^*$-algebra, every commutative $C^*$-algebra corresponds to a 
compact Hausdorff space (its spectrum).  Hence \textsl{the topological K-Theory is the K-Theory of commutative 
$C^*$-algebras} and thus the study of arbitrary $C^*$-algebras corresponds to \emph{noncommutative topology}.
This is the starting point of A. Connes' famous \emph{noncommutative geometry}.\\

Now let $A$ be a unital Banach algebra and let $P_{n}(A)$ denote the set of projections ($P^2 =P$) in $M_{n}(A)$, the set 
of $n\times n$ matrices with entries from $A$. We denote by $P(A)$ the limit 
$$P(A)=\lim _{n\rightarrow\infty }P_{n}(A).$$

{\bf Definition 1.} 
 $$K_{0}A=Gr[\pi _{0}(P(A))],$$
where $Gr$ is the Grothendieck group of the semi-group $\pi _{0}(P(A))$.\\

{\bf Definition 2.} 
 $$K_{1}A=\pi _{0}(GLA).$$

Two alternative definitions are the following:\\

{\bf Definition 1'.} We define $K_{0}(A)$ to be the abelian group defined by generators and relations as follows:\\
$\bullet$ A generator $[p]$ for each projection $p\in M_n(A)$ ($\forall n\in {\mathbb{N}}$).\\
$\bullet$ Relations as follows:\\
{\bf i.} $[0]=0$\\
{\bf ii.} If $p,q$ are homotopic through projections, then $[p]=[q]$.\\
{\bf iii.} $[p\oplus q]=[p]+[q]$.\\

{\bf Definition 2'.} We define $K_{1}(A)$ to be the abelian group defined by generators and relations as follows:\\
$\bullet$ A generator $[u]$ for each unitary $u\in M_n(A)$ ($\forall n\in {\mathbb{N}}$ where unitary means $u^{*}u=uu^{*}=1$).\\
$\bullet$ Relations as follows:\\
{\bf i.} $[1]=0$\\
{\bf ii.} If $u,v$ are homotopic through unitaries, then $[u]=[v]$.\\
{\bf iii.} $[u\oplus v]=[u]+[v]$.\\

{\bf Lemma 1.} 
 $$K_{0}A=K_{0}^{alg}A,$$
namely the definition we gave above for the 0th K-group is the same as the algebraic 0th K-group definition. Moreover in
$P(A)$, homotopy and conjugacy generate the same equivalence relation.\\

BUT 

$$K_{1}A\neq K_{1}^{alg}A,$$
since clearly
$$\pi _{0}(GLA)\neq \pi _{1}(BGLA^{+}).$$
If $J$ is a non-unital $C^*$-algebra, we denote by $J^{+}$ the algebra $J$ with a unit attached to it; then we define the K-groups of $J$ as follows:
$$K_{i}J=Ker(K_{i}(J^{+})\rightarrow K_{i}{\mathbb{C}}),i=0,1.$$
We hope the reader can distinguish between the $"+"$ of the 1-point compactification and the $"+"$ of Quillen's plus 
construction.\\

\textsl{Basic relations:}\\

{\bf 1.} $$K_{i}(C(X))=K^{-i}X.$$

{\bf 2.} If $K(H)$ denotes the set of compact operators on a Hilbert space $H$, then 
$$K_{0}(K(H))=\mathbb{Z}$$
whereas
$$K_{1}(K(H))=0.$$
Recall that
$$K(H)=\lim _{n\rightarrow\infty}M_{n}{\mathbb{C}},$$
hence
$$K_{*}(K(H))=\lim _{n\rightarrow\infty}K_{*}(M_{n}{\mathbb{C}})=
\lim _{n\rightarrow\infty}K_{*}{\mathbb{C}}=K_{*}{\mathbb{C}}.$$

{\bf 3.} If $B(H)$ denotes the set of bounded operators on a Hilbert space $H$, then
$$K_{i}(B(H))=0, i=0,1.$$

Moreover one has that
$$K_{0}({\mathbb{C}})={\mathbb{Z}}$$
whereas
$$K_{1}({\mathbb{C}})=0.$$

Next we ask the question: \textsl{What should be the $C^*$-algebra K-Theory version of the exact sequence of a pair} 
$(X,A)$?\\

Let $U=X/A$. Then there is an exact sequence 
$$0\rightarrow C_{0}(U)\rightarrow C(X)\rightarrow C(A)\rightarrow 0.$$
One has the following\\

{\bf Theorem 1.} To any $C^*$-algebra short exact sequence
$$0\rightarrow J\rightarrow A\rightarrow A/J\rightarrow 0,$$
is associated an exact sequence of K-Theory groups

$$\begin{CD}
K_{1}J @>>> K_1A @>>> K_1(A/J)\\
  @AAA   @.     @VV\partial V\\
K_0(A/J) @<<< K_0A @<<< K_0J
\end{CD}$$\\
\\

where the vertical arrows
$$K_{0}(A/J)\rightarrow K_{1}J$$
and
$$\partial :K_{1}(A/J)\rightarrow K_{0}J$$
make the above diagram commute.\\

{\bf Proof:} Omitted, see the bibliography. $\square$.\\

Let $u\in M_{n}(A/J)$ be invertible. Then $\partial [u]$ should obstruct lifting to an invertible in 
$M_{n}A$.\\

We know we can always lift 
$$\left( \begin{array}{cc}u & 0\\
0 & u^{-1}\end{array} \right),$$
(more precisely its class), to an invertible matrix $R\in M_{2n}A$. Let $P\in M_{2n}A$,
$$P=\left[ \begin{array}{cc}1 & 0\\
0 & 0\end{array} \right]\in M_{2n}(J^{+})$$
and let $Q=R^{-1}PR\in M_{2n}(J^{+})$. Then
$$[P]-[Q]\in K_{0}J.$$
An interesting special case is this: Take $J=K(H), A=B(H)$ hence $A/J=Q(H)$. Then the boundary map $\partial$ is a 
generalised \emph{index}.\\

We shall return to K-Homology (of $C^*$-algebras) next.\\

{\bf Definition 3.} 
Let $A=C(X)$, where $X$ is a compact and metrisable topological space. A Hilbert space 
$Ç$ is an \emph{$X$-module} if there is a representation $\rho :A\rightarrow B(H)$.\\

{\bf Definition 4.} 
An \emph{essential equivalence} between two $X$-modules $H,H'$ is a unitary operator
$u:H\rightarrow H'$ such that
$$u\rho (a)u^*-\rho '(a)\in K$$
is a compact operator.\\

{\bf Definition 5.} 
An $X$-module is \emph{big} if the map
$$\rho :A\rightarrow B(H)\rightarrow Q(H)$$
is injective.\\

Then one has the following important\\

{\bf Theorem 2.} (Voiculescu) \textsl{Big modules are absorbing}, namely 
if $H,H'$ are two $X$-modules and 
$H$ is big, then $H\oplus H'$ is essentially equivalent to $H$.\\ 

{\bf Proof:} See the bibliography. $\square$\\

{\bf Corollary 1.} All big modules are essentially equivalent.\\

Let $X$ be some topological space (compact and metrisable). Choose a big $X$-module $H_X$.\\ 

{\bf Definition 6.} The \emph{Paschke dual} $D(X)$ of $X$ is the $C^*$-algebra
$$D(X)=\{T\in B(H_X):[T,\rho (a)]compact\}, \forall a\in A=C(X).$$
Next we define:\\

{\bf Definition 7.} 
  $$\tilde{K}_{i}X=(\tilde{K}^{i}A=) K_{1-i}(D(X)).$$
There is an exact sequence of a pair $(X,Y)$ of $C^*$-algebras in K-Homology:\\  

{\bf Definition 8.} 
$$D(X,Y)=\{T\in D(X):T\rho (a)\in K(H), \forall a\in C_{0}(X/Y)\}.$$
Because $D(X,Y)\lhd D(X)$, one gets a short exact sequence of $C^*$-algebras:
$$0\rightarrow D(X,Y)\rightarrow D(X)\rightarrow D(X)/D(X,Y)\rightarrow 0.$$

\newpage

\section{Lecture 6 ($C^*$-algebra extensions and K-Homology)}

Let $H$ be a Hilbert space and $N\in B(H)$ a bounded operator.\\

{\bf Definition 1.} $N$ is called  \emph{normal} if
$$NN^*=N^*N.$$

If $N$ is normal, then $N$ generates a commutative $C^*$-subalgebra of $B(H)$, in fact one can show that 
it is isomorphic to $C(\sigma (N))$, where
$$\sigma (N)=\{\lambda\in{\mathbb{C}}:\nexists (N-\lambda I)^{-1}\}$$ 
is the \emph{spectrum} of $N$.\\

{\bf Task:} \textsl{Try to classify normal operators up to unitary equivalence modulo compacts}.\\

Berg proved that these are classified entirely by their \textsl{essential spectrum} 
$\sigma _{e}(N)$, namely
$$\sigma (\pi (N)),\pi :B(H)\rightarrow Q(H)=B(H)/K(H),$$
as a consequence of Voiculescu's theorem.\\

{\bf Definition 2.} An operator $N\in B(H)$ is called \emph{essentially normal} 
if $NN^*-N^*N\in K(H)$ is compact.\\

Next, Brown, Douglas, Fillmore asked themselves the following question: \emph{Is the essential spectrum the only invariant 
of essential equivalence for essentially normal operators?}\\

\emph{The answer is negative}.\\

As a counterexample consider the unilateral shift $U$: We know that $\sigma _{e}(U)=S^1$ and $IndU=-1$ but 
$IndN=0, \forall N$ normal.\\

For any $\lambda$ which is not in the essential spectrum of $N$, the index $Ind(N-\lambda I)\in{\mathbb{Z}}$ is an invariant 
of essential equivalence. So we may ask: If\\
{\bf 1.} $\sigma _{e}(N_1)=\sigma _{e}(N_2)=X$,\\
{\bf 2.} $Ind(N_1 -\lambda _{i}I)=Ind(N_2 -\lambda _{i}I)$, for $\lambda _{i}\in{\mathbb{C}}-X$,\\
are the operators $N_1$ and $N_2$ essentially equivalent?\\

The answer to this question is affirmative.\\

{\bf Definition 3.} We define the set $BDF(X)$ as the set which consists of 
the essential equivalence classes of essentially normal operators with essential spectrum $X$.\\

Then, 
\textsl{$BDF(X)$ equals the set of $C^*$-algebra extensions}\\

$$\begin{CD}
0 @>>> K(H) @>>> A @>>> C(X) @>>> 0\\
@.   @|     @VViV     @VViV    @.\\
0 @>>> K(H) @>>> B(H)  @>>>  Q(H)  @>>> 0
\end{CD}$$\\
\\

where $A$ is the algebra generated by $N$ and $K(H)$ and $i$ denotes the inclusions 
$A\xrightarrow{i} B(H)$, $C(X)\xrightarrow{i} Q(H)$.\\

\textsl{The set $BDF(X)$ is also equal to the set of monomorphisms $C(X)\rightarrow Q(H)$}.\\

Moreover we know the following facts:\\

$\bullet$ $BDF(X)$ is a semi-group under direct sum.\\

$\bullet$ From Voiculescu's theorem, any normal operator defines an identity element.\\ 

{\bf Examples:}\\

{\bf 1.} $BDF(I)=0$, üðïõ $I=[0,1]\subset{\mathbb{C}}$.\\

{\bf 2.} $BDF(S^1)={\mathbb{Z}}$.\\

{\bf Theorem 1.} (Averson) Given any homomorphism $\alpha :C(X)\rightarrow Q(H)$, one can find another homomorphism 
$\phi :C(X)\rightarrow B(H\oplus H)$ such that if 
$$\phi =\left( \begin{array}{cc}\phi _{11}& \phi _{12}\\
\phi _{21} & \phi _{22}\end{array} \right),$$
then
$$\pi\circ\phi _{11}=\alpha .$$

Observe that $\phi _{12}, \phi _{21}$ must be compact (namely compact operator valued) 
and $\phi _{22}$ must be another $*$-homomorphism, say $\beta$, thus $\alpha\oplus\beta$ is liftable to 
$B(H)$.\\

{\bf Proof:} The basic idea is this:
$$\phi (xx^*)=\phi (x)\phi (x^*),$$
ïðüôå
$$\left( \begin{array}{cc}\phi _{11}(xx^*) & \phi _{12}(xx^*)\\
\phi _{21}(xx^*) & \phi _{22}(xx^*)\end{array} \right)=
\left( \begin{array}{cc}\phi _{11}(x)& \phi _{12}(x)\\
\phi _{21}(x) & \phi _{22}(x)\end{array} \right)
\left( \begin{array}{cc}\phi _{11}(x^*)& \phi _{12}(x^*)\\
\phi _{21}(x^*) & \phi _{22}(x^*)\end{array} \right),$$
hence
$$\phi _{11}(xx^*)=\phi _{11}(x)\phi _{11}(x^*)+\phi _{12}(x)\phi _{12}(x^*).$$
$\square$.\\

Observe that
$$BDF(X)\rightarrow Hom(K^{-1}X,\mathbb{Z})$$
and
$$0\rightarrow K(H)\rightarrow A\rightarrow C(X)\rightarrow 0.$$
The corresponding long exact sequence contains  a map 

$$\begin{CD}
K_{1}(C(X)) @>>> K_0(K(H))\\
  @|        @|\\
K^{-1}X @>>> \mathbb{Z}
\end{CD}$$\\
\\

Then the basic result of Brown, Douglas, Fillmore is the following\\

{\bf Theorem 2.} 
$$BDF(X)=K_{1}X.$$

{\bf Proof:} We give a sketch of the proof:
$$\alpha :C(X)\rightarrow Q(H)$$
$$\phi :C(X)\rightarrow B(H\oplus H),$$
where
$$\phi =\left( \begin{array}{cc}\phi _{11}& \phi _{12}\\
\phi _{21} & \phi _{22}\end{array} \right).$$
Let
$$p:H\oplus H\rightarrow H$$
be the projection to the first factor. Then
$$[\phi ,p]=
\left( \begin{array}{cc}0 & -\phi _{12}\\
\phi _{21} & 0\end{array} \right)\in K(H)$$
is exactly a K-cycle.\\

In general, $X$ is some polyhedron in ${\mathbb{R}}^3$ and $BDF(X)\rightarrow Hom(K^{-1}X;{\mathbb{Z}})$.\\

Next we continue with some diagram chasing:\\

$$\begin{CD}
K_1A @>>> K_1X  @>>> K_{1}(X,A)  @>>> K_0A @>>>\\ 
@VVV      @VVV      @VVV  @VVV\\
Hom(K^{-1}A; {\mathbb{Z}}) @>>> Hom(K^{-1}X; {\mathbb{Z}}) @>>> ...
\end{CD}$$
 $\square$.\\

Finally we get:
$$0\rightarrow Ext(K^{0}X; {\mathbb{Z}})\rightarrow K_{1}X\rightarrow Hom(K^{-1}X; {\mathbb{Z}})\rightarrow 0$$
and
$$K^{-1}X\rightarrow{\mathbb{Z}} (=K_{0}(K(H)))\rightarrow K_{0}A\rightarrow K^{0}X\rightarrow 0 (=K_{1}(K(H))).$$

We have also used the relations
$$K_{1}(C(X))=K^{-1}X$$
and
$$K_{1}X=BDF(X)=Ext(C(X))=Hom(K^{-1}X; {\mathbb{Z}}).$$

\newpage

\section{Higher Algebraic K-Theory}

Algebraic K-Theory defines and applies a sequence of functors denoted $K_n$ from rings to abelian groups for all integers $n$. The rings can be either commutative or noncommutative, usually unital but this later assumption can be relaxed at the expense of increased complexity.\\

For historic reasons, the lower algebraic K-groups $K_0$ and $K_1$ are thought of in some different terms from the higher algebraic K-groups $K_n$ for $n\geq 2$. Indeed, the lower groups are more accessible, and have more applications (to this day), than the higher groups. The theory of the higher K-groups is noticeably deeper, and certainly much harder to compute (even when the ring is the ring of integers).\\

The group $K_0$ generalises the construction of the ideal class group of a ring, (originated by the failure of the unique factorization property of integers) using projective modules. Its development in the 1960's and 1970's was linked to attempts to solve a conjecture of Serre on projective modules that now is the Quillen-Suslin theorem; numerous other connections with classical algebraic problems were found in this era. Similarly, $K_1$ is a modification of the group of units in a ring, using elementary matrix theory. Intuitively one could say that $K_1$ encodes information about the homotopy classes of $\infty\times\infty$ matrices with entries from a ring.\\

 The first K-group $K_1(R)$ is important in topology, especially when $R$ is a group ring, because its quotient the Whitehead group contains the Whitehead torsion used to study problems in simple homotopy theory and surgery theory; the group $K_0(R)$ also contains other invariants such as the finiteness invariant. Since the 1980's, algebraic K-theory has increasingly had applications to algebraic geometry. For example, motivic cohomology is closely related to algebraic K-theory.\\ 

The 0th K-Group was defined by Alexander Grothendieck in about 1957 in his studies on Algebraic Geometry, the 1st K-Group was defined by the famous topologist  John Whitehead shortly afterwards whereas the 2nd K-Group was defined by Robert Steinberg. Progress stopped until 1970's when D.G. Quillen appeared; with his trully ingeneous \emph{plus construction} he gave the general definition of the nth K-Group  which includes all previous definitions of Grothendieck, Whitehead and Steinberg for $n=0,1,2$. Some years later Quillen discovered his famous $Q$-construction which is more general with improved functorial properties and can be applied to any exact category (instead of just rings). Quillen's results excited Grothendieck who posted a 600 pages handwritten letter concerning the future of Algebraic Geometry and K-Theory. During 1980's, a topological version of Quillen's Higher Algebraic K-Theory was developed by Waldhausen known as Waldhausen K-Theory. Unfortunately, as it is usually the case with modern mathematics, computations are hard and thus even today there are only a few known examples: complete computations exist only for finite fields (Quillen) and even for the ring of integers, the computations are only known "mudulo torsion" (Borel).\\

\newpage

\section{Lecture 1 (Algebraic Preliminaries and the the Grothendieck Group)}

The key notion in this chapter will be the notion of an $A$-\emph{module} where $A$ is a ring (either commutative or not, usually but not necessarily unital). In the case where $A$ is actually a field, then an $A$-module is just an $A$-vector space. Motivated by the usual notation in the study of vector spaces, we shall consider \textsl{left} $A$-modules.\\

Classical results in linear algebra (finite dimensional vector spaces) state that every vector space has a basis and the number of vectors in all bases of some fixed vector space is constant; this enables one to define the notion of the \emph{dimension} of a vector space. Similar (but not identical) things hold for modules: A \emph{basis} for an $A$-mod $P$ is a subset $\{e_i\}_{i\in I}$ (not necessarily of finite cardinality, if this happens then $P$ is called \textsl{finitely generated}, abreviated to f.g., see below), of $P$ such that every element of $P$ can be written uniquely as a \textsl{finite} sum $\sum _{i}a_{i}e_{i}$ where $a_i\in A$. If $P$ has a fixed basis we call it a \emph{based free module} and we define the \emph{dimension (or rank)} of $P$ to be the \textsl{cardinality} of its given basis. A module is called \emph{free} if there is a basis which makes it into a based free module. The typical example of a free module is $A^n$ which consists of $n$-tuples of elements of the ring $A$. There are rings however for which $A^n\simeq A^{n+m}$ where $m\neq 0$. To avoid this pathology, we shall assume that all our rings have the so-called \emph{invariant basis property} which means precisely that $A^n$ and $A^m$ are not isomorphic unless $n=m$. In this case the dimension of a free $A$-mod $P$ is invariant, independent of the choice of a basis in $P$. (An example of a module which does not have the invariant basis property is the infinite matrix ring $End_{F}(F^{\infty})$ of endomorphisms of an infinite dim vector space $F^{\infty}$ over some field $F$). After these brief algebraic preliminaries, we proceed to the basic material of this section.\\

{\bf Definition 1.} Let $A$ be a unital ring. An $A-mod$ (A-module) $P$ is called \emph{projctive} if there is a module $Q$ such that the direct summand
$$F=P\oplus Q$$
is a free module $F$.\\

It follows from the above definition that $P$ is the image of a projection $p$ on $F$; the module endomorphism in $F$ which is the identity on $P$ and $0$ on $Q$, namely 
$$Q=ker(p),$$
is idempotent and projects $F$ to $P$. In other words $P$ is projective if $P$ is a direct summand of a free module.\\

Using category language, one can define projective modules relying on the lifting property: An $A$-mod $P$ is projective if and only if for every surjective module homomorphism $f: N\twoheadrightarrow M$ and every module homomorphism $g:P\rightarrow M$, there exists a homomorphism $h:P\rightarrow N$ such that $fh=g$ (we do not require the lifting homomorphism $h$ to be unique, this is not a universal property). The following diagram describes the situation:\\

$$\begin{CD}
P @>h>> N\\
  @|    @VVfV\\
P @>g>> M
\end{CD}$$
\\

The notion of a finitely generated module generalises the notion of a finite dimensional vector space (in the case where the ring $A$ is actually a field):\\

{\bf Definition 2.} A (left) $A-mod$ $P$ is called \emph{finitely generated} if and only if there exists a  finite number of elements  $x_1, x_2,...,x_n\in P$ such that $\forall x\in P$ there exist elements $a_1, a_2,...,a_n\in A$ so that
$$x=a_1 x_1 + a_2 x_2 +...+a_n x_n.$$
The set $\{x_1,x_2,...,x_n\}$ is called a set of generators for $P$.\\

{\bf Theorem 1.} The following are equivalent:\\
{\bf 1.} $P$ is a finitely generated projective (f.g.p. for short) $A-mod$.\\
{\bf 2.} There exists $n\in\mathbb{N}^*$ along with another $A-mod$ $Q$ such that
$$P\oplus Q=A^{n}.$$
{\bf 3.} There exists some $n\in\mathbb{N}^*$ along with some element $e\in M_{n}A$ such that
$$e=e^2$$
and
$$P=A^{n}e$$
(where $M_{n}A$ denotes the set of $n\times n$ matrices with entries from $A$).\\

For the proof of this theorem we refer to the books on Algebra listed in the introduction (background references).\\

We denote by $P_A$ the category of f.g.p. $A-mod's$ and by $Iso(P_A)$ we denote the set of isomorphism classes of f.g.p. $A-mod's$. The second set is in fact an abelian monoid with addition given by
$$[P]+[Q]=[P+Q].$$

{\bf Definition 3.} We define the 0th K-group of the ring $A$, denoted $K_{0}A$, the abelian group completion (or the Grothendieck group $Gr$ we used in chapter 1) of $Iso(P_A)$ ), namely\\
 
$K_{0}A:=$ abelian group completion of $Iso (P_A)$.\\

The 0th K-group is often called the \emph{Grothendieck group} of the ring $A$ to honour the great French mathematician (of German origin)
Alexander Grothendieck who defined it in 1957.\\

In general, if $I$ is any abelian monoid, then there exists an Abelian group $I^*$ along with an homomorphism
$\phi :I\rightarrow I^*$, which is a universal homomorphism from $I$ to some abelian group. The pair $(I^* ,\phi )$ is only defined up to a canonical isomorphism. \\

There are three ways to construct $I^*$:\\

{\bf 1.} $I^*$ is the free abelian group with generators $[a], a\in I$ and relations $[a+b]=[a]+[b]$.\\

{\bf 2.} $I^* = I\times I/\sim$, where the equivalence relation $\sim$ is defined as follows: $(a,b)\sim (a',b')\Leftrightarrow
\exists\gamma :a+b'+\gamma =a'+b+\gamma$ (clearly this is indeed an equivalence relation). Next we define
$(a,b)+(a',b')=(a+a',b+b')$ whereas $0=(0,0)$ and $-(a,b)=(b,a)$. Thus $(a,b)=[a]-[b]$. (This is the construction we used in section 1.1 to define $KX$, the Grothendieck group of $Vect(X)$).\\

{\bf 3.} Assume there exists an $a_0\in I$ such that $\forall a'\in I, \exists n\in\mathbb{N}$ and $b\in I$ so that 
$a+b=nd_0$.  Then
$$I^*=I\times\mathbb{N}/\sim$$
where $(a,n)\sim (a',n')\Leftrightarrow\exists n\in\mathbb{N}$ such that $a+n'a_0=a'+na_0$.\\

All the above abelian groups have the desired universal property.\\

Hence, to summarise, the 0th K-group $K_{0}A$ of a unital ring $A$ is the free abelian group with generators $[p]\forall p\in P_{A}$ and relations $[p\oplus q]=[p]+[q]$, which is equal to the group of formal differences $[p]-[q]$ and which is also equal to the group of formal differences $[p]-[A^n]$.\\
 
{\bf Examples:}\\

{\bf 1.} Let $\mathbb{F}$ be a skew-field, namely a ring with $0\neq 1$ where division is possible, in other words each element has a multiplication inverse (or equivalently a field where multiplication is not necessarily commutative). Then the category $P_{\mathbb{F}}$ of f.g.p. $\mathbb{F}$-mod's consists of finite dimensional vector spaces over the field $\mathbb{F}$. Then the abelian monoid  
$Iso (P_{\mathbb{F}})$ is equal to $\mathbb{N}$, and hence 
$$K_{0}\mathbb{F}=\mathbb{Z}.$$\\

{\bf 2.} If $A=\mathbb{Z}$, then the category $P_{\mathbb{Z}}^n$ is the finitely generated free abelian group
${\mathbb{Z}}^n$ for all $n\in\mathbb{N}$ whereas 
$$Iso (P_{\mathbb{Z}})=\mathbb{N},$$ 
hence
$$K_0 \mathbb{Z}=\mathbb{Z}.$$ 
The same holds for principal ideal domains.\\

{\bf 3.} \textsl{Dedekind domains}, eg take $\mathbb{F}$ some number field $[F,Q]<\infty$ and  
$A$ the integral quotient of $\mathbb{Z}$ in $\mathbb{F}$. Let $Pic (A)$  denote the ideal class group of $A$-fractional ideals divided by the principal ones. Then one has the following:\\

\textsl{Theorem 2:} If $P\in P_{A}$, then $P\simeq a_1 \oplus a_2 \oplus ... \oplus a_n$, where the $a_i$'s are fractional. 
Hence $K_{0}A=\mathbb{Z}\oplus$ $Pic (A).$\\

The proof of this theorem is left as an excercise to the reader (see for example \cite{weibel}).\\

{\bf 4. Serre-Swan Theorem}. Let $X$ be a compact Hausdorff topological space and let $A=C(X)$ denote the set of continuous, complex valued functions on $X$ (which is a unital ring). Then the category  
$P_A$ is equivalent to the category of complex vector bundles over $X$.\\

[\emph{Aside Note:} This theorem was one of the motivations for the development of \textsl{Noncommutative Geometry} by the French mathematician Alain Connes which, in some sense, is an attempt to unify differential and algebraic geometry].\\

{\bf Proof of Serre-Swan Theorem:} This is an important theorem, thus we shall give a full proof of it. We have to prove the existence of a bijection between $Iso(P_A)$ (which is the set of isomorphism classes of f.g.p. $A$-mod's) and $Vect(X)$ (which is the set of isomorphism classes of complex vector bundles over $X$) where $A=C(X)$.\\

Recall the definitions of vector bundles and local sections of
vector bundles in section 1 of chapter 1. Given any vector bundle $E\rightarrow X$ we denote by $\Gamma (E,X)$ (or simply $\Gamma (E)$) the set of all local sections of the vector bundle $E\rightarrow X$. It is clear that $\Gamma (E)$ becomes a $C(X)$-module if we define:\\
$\bullet$ $(s_1+s_2)(x)=s_1(x)+s_2(x) \forall s_1,s_2\in\Gamma (E)$\\
$\bullet$ $(as)(x)=a(x)s(x) \forall a\in C(X)$ and $s\in\Gamma (E)$\\
It is clear that for the trivial vector bundle ${\mathbb{C}}^{n}_{X}$, the set $\Gamma ({\mathbb{C}}^{n}_{X})$ is a \textsl{free} $C(X)$-module on $n$ generators. 
Moreover $\Gamma (-)$ is an additive functor from the category of complex vector bundles over $X$ to the category $P_A$. Thus to prove the theorem, we have to prove three claims:\\

{\bf Claim 1.} The functor $\Gamma (-)$ "respects" isomorphisms, namely that if $E\simeq E'$ as vector bundles over $X$, then $\Gamma (E)\simeq\Gamma (E')$ as $C(X)$-modules.\\ 

{\bf Claim 2.} $\Gamma (E)$ is a f.g.p. $C(X)$-mod for any vector bundle $E$ over $X$.\\

{\bf Claim 3.} All f.g.p. $C(X)$-mod's can be "realised" as a $\Gamma (E)$ for some vector bundle $E$ over $X$.\\

The proof of Claim 1 requires a number of propositions (we prove it for the more general case where $X$ is a normal topological space and to be more precise we shall prove that if $X$ is normal then $\Gamma (-)$ gives an isomorphism $Hom (E,E')\simeq Hom_{C(X)}(\Gamma (E),\Gamma (E'))$):\\

{\bf Lemma 1.} If $X$ is normal, suppose $U$ is a neighborhood of  $x\in X$ and let $s$ be a section of a vector bundle $E$ over $U$. Then there is a section $s'$ of $E$ over $X$ so that $s$ and $s'$ agree in some neighborhood of $x$.\\

{\bf Proof:} We shall construct an $s'$: Let $V,W$ be neighborhoods of $x$ so that $\bar{V}\subset U$ and $\bar{W}\subset V$. Let $a$ be a real valued function on $X$ such that $a|_{\bar{W}}=1$ and $a|_{(X-V)}=0$. We let $s'(y)=a(y)s(y)$ if $y\in U$ and $s'(y)=0$ if $y\notin U$. $\square$\\

{\bf Corollary 1.} If $X$ is normal, then for any $x\in X$ there are elements $s_1,s_2,...,s_n\in\Gamma (E)$ which form a local base at $x$.\\

{\bf Proof:} (Obvious). $\square$\\

{\bf Corollary 2.} If $X$ is normal, $f,g:E\rightarrow E'$ are two vector bundle maps and $\Gamma (f)=\Gamma (g):\Gamma (E)\rightarrow\Gamma (E')$, then $f=g$.\\

{\bf Proof:} Given $e\in E$ with $\pi (e)=x$ (where $\pi :E\rightarrow X$ is the vector bundle projection), there is a section $s$ over a neighborhood $U$ of $x$ with $s(x)=e$. By Lemma 1 there is a section $s'\in\Gamma (E)$ with $s'(x)=e$. Now $f(e)=fs'(x)=(\Gamma (f)s')(x)=(\Gamma (g)s')(x)=g(e)$. $\square$\\

{\bf Lemma 2.} If $X$ is normal and if $s\in\Gamma (E)$ with $s(x)=0$, then there exist elements $s_1,...,s_k\in\Gamma (E)$ and $a_1,...,a_k\in C(X)$ such that $a_i(x)=0$ for $i=1,...,k$ and $s=\sum_{i}a_is_i$.\\

{\bf Proof:} From Corollary 1 suppose $s_1,...,s_n\in\Gamma (E)$ be a local base at $x$ and let $s(y)=\sum_{i}b_i(y)s_i(y)$ near $x$ where $b_i(y)\in{\mathbb{C}}$. Let $a_i\in C(X)$ be such that $a_i$ and $b_i$ agree in a neighborhood of $x$ (these exist by Lemma 1 applied to $X\times {\mathbb{C}}$). Then $s'=s-\sum_{i}a_is_i$ vanishes in an neighborhood $U$ of $x$. Let $V$ be a neighborhood of $x$ so that $\bar{V}\subset U$. Let $a\in C(X)$ be zero at $x$ and $1$ on $X-V$. Then $s=as'+\sum_{i}a_is_i$. But $a(x)=0$ and $a_i(x)=b_i(x)=0$. $\square$\\

{\bf Corollary 3.} Let $I_x$ be the two-sided ideal of $C(X)$ consisting of all $a\in C(X)$ with $a(x)=0$. Then $\Gamma (E)/[I_{x}\Gamma (E)]\simeq \pi ^{-1}(x)$ (where $\pi :E\rightarrow X$ is the vector bundle projection), the isomorphism being given by $s\mapsto s(x)$.\\

{\bf Proof:} This follows from Lemma 2 and (the proof of) Corollary 2. $\square$\\

{\bf Proposition 1.} If $X$ is normal, then given any $C(X)$-map $F:\Gamma (E)\rightarrow\Gamma (E')$, there exists a unique (complex) vector bundle map $f:E\rightarrow E'$ so that $F=\Gamma (f)$.\\

{\bf Proof:} Uniqueness follows from Corollary 2. Now $F$ induces a map $f_x:\Gamma (E)/[I_{x}\Gamma (E)]\rightarrow \Gamma (E')/[I_{x}\Gamma (E')]$. The totality of these yield a map $f:E\rightarrow E'$ which is linear on fibres. If $s\in\Gamma (E)$, then $(fs)(x)=f_{x}s(x)=(F(s))(x)$ by construction and thus $F=\Gamma (f)$. The final step is to check continuity: Let $s_1,...,s_m\in\Gamma (E)$ be a local base at $x$. If $e\in E$ and $\pi (e)$ is near $x$, one has $e=\sum_ia_i(e)s_i(\pi (e))$ where the $a_i$'s are continuous complex valued functions. Now $f(e)=\sum_ia_i(e)fs_i(\pi (e))$. Since $fs_i=F(s_i)$, this implies that $fs_i$ is a continuous section of $E'$; yet  all terms in the sum are continuous in $e$, hence $f$ is continuous. 
$\square$\\

{\bf Corollary 4.}   Let $X$ be normal and let $E$ and $E'$ be to complex vector bundles over $X$. Then $E\simeq E'\Leftrightarrow\Gamma (E)\simeq\Gamma (E')$ as $C(X)$-modules.\\

{\bf Proof:} Obvious from Proposition 1. $\square$\\

To prove Claim 2 we need two Lemmas:\\

{\bf Lemma 3.} If $X$ is compact and Hausdorff, let $E$ be any complex vector bundle over $X$. Then there is a trivial complex vector bundle
${\mathbb{C}}^{n}_{X}$ (for some $n\in {\mathbb{N}}$) along with a surjection $f:{\mathbb{C}}^{n}_{X}\rightarrow E$.\\
 
{\bf Proof:} Recall that Propositon 1.1.1 states that if the base space $X$ is compact and Hausdorff, then for every complex vector bundle $E$ over $X$ there exists another complex vector bundle $E'$ over $X$ so that their direct sum is the trivial bundle, namely 
$$E\oplus E'\simeq {\mathbb{C}}^{n}_{X}$$ 
for some $n\in {\mathbb{N}}$. (We proved the theorem for the real case but this can be carried over to the complex case in a straightforward way). Then Lemma 3 is a direct consequence of this proposition where the surjection $f$ is just the natural projection
$$f:{\mathbb{C}}^{n}_{X}\simeq E\oplus E'\rightarrow E.$$
$\square$\\

{\bf Lemma 4.} For any complex vector bundle $E\rightarrow X$ over $X$, where $X$ is compact and Hausdorff, $\Gamma (E)$ is a f.g.p. $C(X)$-module.\\

{\bf Proof:}  It is clear that for the trivial vector bundle ${\mathbb{C}}^{n}_{X}$, the set $\Gamma ({\mathbb{C}}^{n}_{X})$ is a \textsl{free} $C(X)$-module on $n$ generators. Since $E\oplus E'\simeq {\mathbb{C}}^{n}_{X}$, one then has that (from Corollary 4) $\Gamma (E)\oplus\Gamma (E')\simeq \Gamma ({\mathbb{C}}^{n}_{X})$ as $C(X)$-modules, and therefore $\Gamma (E)$ is a f.g.p. $C(X)$-module. $\square$\\

Finally, for the proof of Claim 3 we need the following Lemma and Proposition:\\

{\bf Lemma 5.} Let $t_1,...,t_k$ be sections of a vector bundle $E$ over a neighborhood $U$ of $x$ so that $t_1(x),...,t_k(x)$ are linearly independent. Then there exists a neighborhood $V$ of $x$ so that $t_1(y),...,t_k(y)$ are linearly independent $\forall y\in V$.\\

{\bf Proof:}  Let $s_1,...,s_n$ be a local base at $x$ and let $t_i(y)=\sum_j a_{ij}(y)s_j(y)$ where $y\mapsto a_{ij}(y)$ is a continuous map $U\rightarrow GL_n({\mathbb{C}})$. Then by hypothesis a $k\times k$ submatrix of $a_{ij}(y)$ must be nonsingular and this should hold for all $y$ sufficiently close to $x$. From this fact the result follows. $\square$\\

{\bf Proposition 2.} Let $f:E\rightarrow E'$ be a vector bundle map. Then T.F.A.E.:\\
{\bf (1)} $Imf$ is a subbundle of $E'$.\\
{\bf (2)} $Kerf$ is a subbundle of $E$.\\
{\bf (3)} The dimensions of the fibres of $Imf$ are locally constant.\\
{\bf (4)} The dimensions of the fibres of $Kerf$ are locally constant.\\

{\bf Proof:} Let us start by making a comment: We know from linear algebra that for any linear map between finite dimensional vector spaces $f:V\rightarrow V'$, $Kerf$ is a vector subspace of $V$ and $Imf$ is a vector subspace of $V'$. This cannot be carried over to vector bundles, namely if $f:E\rightarrow E'$ is a vector bundle map, then $Kerf$ and $Imf$ are not necessarily vector subbundles of $E$ and $E'$ respectively. For example, let $X=I=[0,1]$ the unit interval, $E=I\times {\mathbb{C}}$ and $\pi (x,y)=x$. Let $f:E\rightarrow E$ be given by $f(x,y)=(x,xy)$. Then the image of $f$ has a fibre of dim 1 everywhere except at $x=0$ where the fibre is zero, thus $Imf$ cannot be a vector bundle. Neither is $Kerf$. This, however is the only thing which can go wrong.\\
Back to our proof, it is clear that (3) and (4) are equivalent and they are implied by either (1) or (2). To see that (3) implies (1), let $x\in X$, choose a local base $s_1,...,s_m$ for $E$ at $x$ and a local base $t_1,...,t_n$ for $E'$ at $x$. Let $k$ be the dimension
 of the fibre of $Imf$ at $x$. After a possible renumbering, we can assume that $fs_1(x),...,fs_k(x)$ span the fibre of $Imf$ at $x$ and so thy are linearly independent. By another possible renumbering, we can assume that $fs_1(x),...,fs_k(x),t_{k+1}(x),...,t_n(x)$ are linearly independent and hence by local constancy of the dimension of the fibre of $E'$ and Lemma 5 one has that $fs_1,...,fs_k,t_{k+1},...,t_n$ form a local base for $E'$ at $x$. By the hypothesis of Lemma 5 $fs_1,...,fs_k$ form a local base for $Imf$ at $x$ which means that $Imf$ is a subbundle of $E'$.\\
To see that (3)$\Rightarrow$(2), let $s_1,...,s_m$ be as above. For all $y$ near $x$ we can write $fs_i(y)=\sum_{j=1}^{k}a_{ij}(y)fs_j(y)$ for $i>k$. Let $s'_{i}(y)=s_i(y)-\sum_{j=1}^{k}a_{ij}(y)s_j(y)$. Then $s'_{k+1},...,s'_{m}$ are local sections of $Kerf$ and they are linearly independent near $x$. Since there are exactly the correct number of them, they form a local base for $Kerf$ and hence $Kerf$ is a subbundle of $E$.\\ 
\emph{Remark:} Without any hypothesis, this proof shows that if $dimF_x (Imf)=n$, then $dimF_y (Imf)\geq n$ for all $y$ in some neighborhood of $x$. $\square$\\
   
{\bf Proof of Claim 3:} Suppose that $P$ is f.g.p. Then $P$ is a direct summand of a f.g. free $C(X)$-module $F$. Therefore there exists an idempotent endomorphism $g:F\rightarrow F$ with $P\simeq Img$. Now $F=\Gamma ({\mathbb{C}}_{X}^{n})$ for some $n$. By Proposition 1 above, $g=\Gamma (f)$ where $f:{\mathbb{C}}_{X}^{n}\rightarrow {\mathbb{C}}_{X}^{n}$. Since $g^2=g$, Proposition 1 implies that $f^2=f$ as well. If one knew that $E=Imf$ was a subbundle of ${\mathbb{C}}_{X}^{n})$ (for some $n$, namely the trivial bundle) one would have by Lemma 1.1.1 that ${\mathbb{C}}_{X}^{n})\simeq E\oplus E'$ where $E'\simeq Ker(f)$ and so $P\simeq Im\Gamma (f)=\Gamma (E)$ since $\Gamma (-)$ is an additive functor. By Proposition 2 above it saffices to show that $dimF_{x}(E)$ is locally constant (where recall that $F_{x}(E)$ denotes the fibre of $E$ over $x$). Since $f^2=f$, then $E=Kerf=Im(1-f)$ and $F_x({\mathbb{C}}_{X}^{n})=F_x(E)\oplus F_x(E')$. Suppose $dimF_x(E)=h$ and $dimF_x(E')=k$. Apply Remark 1 at the end of the proof of Proposition 2 to $f$ and $1-f$ respectively and get that   
$dimF_y(E)\geq h$ and $dimF_y(E')\geq k$ for all $y$ in some neighborhood of $x$. Yet 
$$dimF_y(E) + dimF_y(E')=h+k$$        
is a constant, thus the dimension of the trivial bundle is locally constant. $\square$\\

There is thus a complete analogy between vector bundles and f.g.p. modules; that enables one to carry over many of the notions we met in chapter 1 (topological K-Theory) to this algebraic setting (e.g. define stably isomorphic modules etc).\\

\newpage

\section{Lecture 2 (The Whitehead and the Steinberg Groups)}

We denote by $GL_{n}A$ the set (multiplicative group in fact) of $n\times n$ invertible matrices with entries from the ring $A$. Clearly 
$$GL_{n}A\subset GL_{n+1}A\subset ...\subset GLA=\cup _n GL_{n}A,$$
where the inclusion is given by 
$$a\hookrightarrow  
\left( \begin{array}{cc}a & 0\\
0 & 1\end{array} \right).$$
We denote by $e_{ij}$ the matrix with entry 1 in the position $(i,j)$ and zero everywhere else $(i\neq j)$. Clearly
$$e_{ij}e_{kl}=\delta _{jk}e_{il}.$$
Let $a\in A$. Then we set
$$e^{a}_{ij}=1+ae_{ij}$$
whereas
$$e^{a}_{ij}e^{b}_{ij}=e^{a+b}_{ij}.$$
Recall that in a multiplicative group the commutator of two elements is defined by
$$[x,y]=xyx^{-1}y^{-1}$$
and the inverse commutator is defined as 
$$[x,y]^{-1}=[y,x].$$
Using the above definitions, we compute the commutator
$$[e^{a}_{ij},e^{b}_{kl}]= \left\{ \begin{array}{rcl}
1, & \mbox{if} & j\neq k, i\neq l\\ e^{ab}_{il}, & \mbox{if} & j=k, i\neq l\\ e^{-ba}_{kj}, & \mbox{if} & j\neq k, i=l
 \end{array}\right.$$

Moreover we compute the following quantities:

$$e^{a}_{ij}e^{b}_{jk}=(1+ae_{ij})(1+be_{jk})=1+ae_{ij}+be_{jk}+abe_{ik},$$
$$e^{a}_{ij}e^{b}_{jk}e^{-a}_{ij}=(1+ae_{ij}+be_{jk}+abe_{ik})(1-ae_{ij})=1+ae_{ij}+be_{jk}+abe_{ik}-ae_{ij},$$
ïðüôå
$$[e^{a}_{ij},e^{b}_{jk}]=(1+ae_{ij}+be_{jk}+abe_{ik})(1-be_{jk})=1+be_{jk}+abe_{ik}-be_{jk}=e^{ab}_{ik}$$
êáé
$$[e^{a}_{ij},e^{b}_{ki}]^{-1}=[e^{b}_{ki},e^{a}_{ij}]=e^{ba}_{kj}\Rightarrow [e^{a}_{ij},e^{b}_{ki}]=e^{-ba}_{kj}.$$

{\bf Definition 1.} We denote by $E_{n}A$ the subgroup of $GL_{n}A$ which is generated by the elements $e^{a}_{ij}$ 
for $1\leq i,j\leq n$, $i\neq j$ and some $a\in A$. Moreover we denote  $EA$ the union of all $E_{n}A$:
$$EA=\cup _n E_{n}A.$$

{\bf Definition 2.} A group $G$ is called \emph{perfect} if it is equal to its commutator subgroup
$$G=[G,G],$$
where $[G,G]$ is the subgroup of $G$ generated by elements of the form $[g,g']\forall g,g' \in G$.\\

{\bf Definition 3.} We define the maximal abelian quotient group $G_{(ab)}$ of $G$ as follows:
$$G_{(ab)}=G/[G,G].$$  
($ab$ stands for \textsl{"abelian"} and they are not indices).\\

{\bf Proposition 1.} The group $E_{n}A$ is perfect for $n\geq 3$.\\

{\bf Proof:} Since 
$$e^{a}_{ik}=[e^{a}_{ij},e^{1}_{jk}]\in [E_{n}A,E_{n}A],$$
given $i,k$, we choose $j\neq i$ and $j\neq k$ and we see that all generators are commutators. $\square$.\\

{\bf Lemma 1.} (\emph{Whitehead Lemma}).
$$EA=[GLA, GLA]=[EA,EA].$$
{\bf Proof:} Let $a\in GL_{n}A$, thus
$$\left( \begin{array}{cc}1 & a\\
0 & 1\end{array} \right)=\Pi e^{a_{ij}}_{ij}\in GL_{2n}A$$
where $1\leq i\leq n$ and $n+1\leq j\leq 2n$.
Then
$$\left( \begin{array}{cc}1 & a\\
0 & 1\end{array} \right)
\left( \begin{array}{cc}1 & 0\\
-a^{-1} & 1\end{array} \right)  
\left( \begin{array}{cc}1 & a\\
0 & 1\end{array} \right)=
\left( \begin{array}{cc}0 & a\\
-a^{-1} & 0\end{array} \right)\in GL_{2n}A.$$
Moreover
$$\left( \begin{array}{cc}a & 0\\
0 & a^{-1}\end{array} \right)=
\left( \begin{array}{cc}0 & a\\
-a^{-1} & 0\end{array} \right)
\left( \begin{array}{cc}0 & -1\\
1 & 0\end{array} \right)\in E_{2n}A.$$
Next, let $a,b\in GL_{n}A$. Then, for the above diagonal (in block form) matrices we have:
$$diag([a,b],1,1)=diag(a,a^{-1},1)diag(b,1,b^{-1})[diag(a^{-1},a^{-1},1)]^{-1}[diag(b,1,b^{-1})]^{-1}$$
(they all belong to the group $E_{3n}A$).\\

Hence $[GL_{n}A,GL_{n}A]$ in $GL_{3n}A$ is contained in  $E_{3n}A$. Taking the union for all $n$ we find that 
$[GL_{n}A,GL_{n}A]\subset EA$. Therefore
$$EA=[EA,EA]\subset [GLA,GLA]\subset EA.$$
$\square$.\\

{\bf Definition 4.} The first K-Group $K_{1}A$ of the ring $A$ is defined as the maximal abelian quotient group of $GLA$, namely
$$K_{1}A=GLA_{(ab)}=GLA/EA.$$

The first K-Group is called the \textsl{Whitehead group} to honour the great british topologist John Whitehead who defined it.\\

{\bf Examples:}\\
{\bf 1.} Let $A=\mathbb{F}$ a field. Left multiplication by the matrices $e^{a}_{ij}$ adds $j$-row 
to the $i$-row. It is known that\\

$E_{n}(\mathbb{F})=$ ker $\{GL_{n}({\mathbb{F}})\xrightarrow{det} {\mathbb{F}}^{x}\}$\\

where ${\mathbb{F}}^{x}$ denotes the non-zero elements under the action of $x$ and the above map is given by the matrix \textsl{determinant}. Next we consider the quotient
$$GL_{n}({\mathbb{F}})/E_{n}({\mathbb{F}})={\mathbb{F}}^{x}=K_{1}(\mathbb{F}).$$
 
{\bf Definition 5.} The \textsl{Steinberg Group} $St_{n}A$ ($StA$ is the case where $n\rightarrow\infty$) is the group with generators
$x^{a}_{ij}$, $i\neq j$ and $a\in A$ under the relations
$$x^{a}_{ij}x^{b}_{ij}=x^{a+b}_{ij}.$$ 
The commutators are given by the following relations:

$$[x^{a}_{ij},x^{b}_{kl}]= \left\{ \begin{array}{rcl}
1, & \mbox{if} & j\neq k, i\neq l\\ x^{ab}_{ik}, & \mbox{if} & j=k, i\neq l\\ x^{-ba}_{kj}, & \mbox{if} & j\neq k, i=l
 \end{array}\right.$$

There is a canonical surjection
$$\phi :StA\rightarrow EA$$
such that
$$\phi (x^{a}_{ij})=e^{a}_{ij}.$$

{\bf Proposition 2.} The Kernel of $\phi$ above is the centre of the group $StA$.\\

{\bf Definition 6.} The second K-group $K_{2}A$ is defined as the kernel of the above canonical surjection $\phi$:
$$K_{2}A=Ker\{\phi :StA\rightarrow EA\}.$$

{\bf Proof of Proposition 2.} Let $C_n$ be the subgroup of $StA$ generated by $x^{a}_{in}$ for $i\neq n$ and $a\in A$.
Let $\phi (C_n)$ be the subgroup of $EA$ defiend by $e^{a}_{in}$ for $i\neq n$ and $a\in A$.\\
{\bf Claim:}  $\phi :C_{n}\rightarrow\phi (C_n)$ is an isomorphism, namely the restriction of $\phi$ to $C_n$ 
is injective.\\

{\emph{Proof of the Claim:}\\
 Since the $x^{a}_{in}$'s commute, in the map
$$\oplus _{i\neq n}A\rightarrow C_n$$
defined by
$$(a_i)_{i\neq n}\mapsto\Pi _{i\neq n}x^{a_{i}}_{in},$$
the product does not depend on the order and moreover it is a homomorphism because:
$$(a_i + a'_{i})\mapsto\Pi _{i\neq n}x^{a_i + a'_{i}}_{in}=
\Pi _{i\neq n}x^{a_i}_{in}x^{a'_{i}}_{in}=
\Pi _{i\neq n}x^{a_i}_{in}\Pi _{i\neq n}x^{a'_{i}}_{in}.$$

Next we define the group $R_n$ (which is a subgroup of the group $StA$) as the group generated by 
$x^{a}_{ni}$ for $i\neq n$ and $a\in A$. A similar argument proves that the restriction of $\phi$ in $R_n$ is injective.\\

To complete the proof of Proposition 2, that the kernel of $\phi$  is the centre of $StA$, we consider an arbitrary element $\alpha\in Ker\phi$ and we write it as a finite product of factors $x^{a}_{ij}$. Choose some $n$ different from $i,j$ occuring in the specific representation of $\alpha$. Then $\alpha$ 
normalises $C_n$, namely 
$$\alpha C_{n}\alpha ^{-1}\subset C_{n}$$
whereas
$$x^{a}_{ij}x^{b}_{kn}x^{-a}_{ij}= \left\{ \begin{array}{rcl}
x^{b}_{kn}, & \mbox{åÜí} & k\neq i,j, \\x^{ab}_{in}x^{-b}_{jn}, & \mbox{åÜí} & k=j
 \end{array}\right.$$
and
$$x^{a}_{ij}x^{b}_{jn}x^{-a}_{ij}x^{-b}_{jn}=x^{ab}_{in}.$$
Let $\gamma\in C_n$, then $\alpha\gamma\alpha ^{-1}\in C_n$ and moreover 
$$\phi (\alpha\gamma\alpha ^{-1})=\phi (\gamma )$$
since $\alpha\in Ker\phi$. But the restriction of $\phi$ in $C_n$ is injective and hence $\alpha\gamma\alpha ^{-1}=\gamma$, 
thus $\alpha$ centralises $C_n$. Similarly $\alpha$ centrilises $R_n$. But the union
$C_n\cup R_n$ generates $StA$ and the centre $ÅÁ$ is $1$. $\square$.\\

\newpage

\section{Lecture 3 (Central Extensions of Groups)}

We saw that
$$K_{1}A=GLA/EA.$$

For topological spaces $X$ one has  
$$K^{-1}X=[X,GL\mathbb{C}],$$
where the RHS (as we know from Homotopy theory) denotes the set of homotopy classes of continuous maps from $X$ to $GL\mathbb{C}$,
and
$$Hom _{(spaces)}(X,GL_{n}{\mathbb{C}})=GL_{n}(C(X)),$$
where $C(X)$ denotes the set (unital ring) of continuous maps from $X$ to $\mathbb{C}$.\\

Let 
$$\left( \begin{array}{cc}1 & ta\\
* & 1\end{array} \right)$$
be a homotopy from
$$\left( \begin{array}{cc}1 & 0\\
0 & 1\end{array} \right)$$
to
$$\left( \begin{array}{cc}1 & a\\
* & 1\end{array} \right).$$

Let $X$ be a topological space and $SX$ its suspension. We know that 
$$K^{-2}X=K^{-1}SX=[SX,GL{\mathbb{C}}]=[S^1, [X,GL\mathbb{C}]]$$
where
$$[X,GL{\mathbb{C}}]=Hom(X,GL{\mathbb{C}})=GL(C(X)).$$

The group $StA$ has generators $x^{a}_{ij}$ with $1\leq i,j<\infty$ and $a\in A$ along with relations
$$x^{a}_{ij}x^{b}_{ij}=x^{a+b}_{ij}$$
and commutators
$$[x^{a}_{ij},x^{b}_{kl}]= \left\{ \begin{array}{rcl}
1, & \mbox{if} & j\neq k, i\neq l\\ x^{ab}_{ik}, & \mbox{if} & j=k, i\neq l\\ x^{-ba}_{kj}, & \mbox{if} & j\neq k, i=l
 \end{array}\right.$$

One has the following central extension 
$$1\rightarrow K_{2}A\rightarrow StA\rightarrow EA\rightarrow 1,$$
namely $K_{2}A$ is contained in the centre of $StA$, with
$$x^{a}_{ij}\mapsto e^{a}_{ij}.$$

{\bf Theorem 1.} Given a central extension
$$1\rightarrow C\rightarrow Y\xrightarrow{\phi} StA\rightarrow 1,$$
it splits, namely there exists a homomorphism $s:StA\rightarrow Y$ such that
$\phi \circ s=identity$\\

{\bf Corollary 1.} If $Y=[Y,Y]$, then $Y$ is isomorphic to $StA$.\\

 {\bf Proof of Theorem 1.} The basic idea is the following: Suppose  $y_1,y_2 \in Y$ are such that $\phi (y_1)=\phi (y_2)$, 
namely $y_1 =cy_2$ with $c\in C$. Then
$$[y_1,y']=[cy_2,y']=cy_{2}y'(cy_2)'y'^{-1}=cy_2 y'y_{2}^{-1}c^{-1}y'^{-1}=[y_2,y'],$$
namely for $x\in StA$, then $[\phi ^{-1}(x),y']$ is a well-defined element of $Y$ and similarly
$[\phi ^{-1}(x),\phi ^{-1}(x')]$ is a well-defined element of $Y$. Then
$$x^{a}_{ij}=[x^{a}_{in},x^{1}_{nj}]$$
for $n\neq i,j$, so we shall try to define $s$ by
$$s(x^{a}_{ij})=[\phi ^{-1}x^{a}_{in},\phi ^{-1}x^{1}_{nj}].$$
The hard point is to prove independence of $n$ which we shall do using 2 Lemmas below:\\

{\bf Lemma 1.}   
$$[\phi ^{-1}x^{a}_{ij},\phi ^{-1}x^{b}_{kl}]=1$$
if $j\neq k$ and $i\neq l$.\\

{\bf Proof:} Choose $n$ different from $k,l,i,j$. Then
$$[\phi ^{-1}x^{b}_{kn},\phi ^{-1}x^{1}_{nl}]\subset \phi ^{-1}[x^{b}_{kn},x^{1}_{nl}]=\phi ^{-1}x^{b}_{kl}.$$
Pick  $v\in\phi ^{-1}x^{b}_{kn}$ and $w\in\phi ^{-1}x^{1}_{nl}$.
Next choose $u\in \phi ^{-1}x^{a}_{ij}$ and it saffices to prove that  $[u,[v,w]]=1$. One has:
$$u[v,w]u^{-1}=[uvu^{-1},uwu^{-1}]=[v,w]$$
if $uvu^{-1}$ and $v$ along with $uwu^{-1}$ and $w$ are congruent. But $\phi (u), \phi (v)$ commute and also $\phi (u), \phi (w)$ commute.
$\square$.\\

Recall the following identities which apply to any group:
$$[x,[y,z]]=[xy,z][z,x][z,y]$$
and
$$[xy,z]=[x,[y,z]][y,z][x,z].$$
Using them we come to the second Lemma:\\

{\bf Lemma 2.}   Let $h,i,j,k$ be distinct and $a,b,c\in A$. Then
$$[\phi ^{-1}x^{a}_{hi},[\phi ^{-1}x^{b}_{ij},\phi ^{-1}x^{c}_{jk}]]=
[[\phi ^{-1}x^{a}_{hi},\phi ^{-1}x^{b}_{ij}],\phi ^{-1}x^{c}_{jk}].$$

{\bf Proof:} Pick some $u\in\phi ^{-1}x^{a}_{hi}$, $v\in\phi ^{-1}x^{b}_{ij}$ and $w\in\phi ^{-1}x^{c}_{jk}$. From Lemma 1 we have
$$[u,w]=1$$
whereas
$$[u,v]\subset\phi ^{-1}[x^{a}_{hi},x^{b}_{ij}]=\phi ^{-1}x^{ab}_{hj}$$
which commutes with $u,v$ whereas $[v,w]$ commutes with $v,w$. Moreover 
$$[u,[v,w]]\subset [\phi ^{-1}x^{a}_{hi},\phi ^{-1}x^{bc}_{ik}]\subset\phi ^{-1}x^{abc}_{ik}$$
which commutes with $u,v,w$. Similarly $[[u,v],w]$ commutes with  $u,v,w$.\\

One also has:
$$[x,[y,z]]=[xy,z][z,x][z,y]$$
and
$$[u,[v,w]]=[uv,w][w,u][w,v]=[[u,v]vu,w][w,v]=[vu[u,v],w][w,v]$$
whereas
$$[xy,z]=[x,[y,z]][y,z][x,z].$$
Then 
$$[u,[v,w]]=[uv,w][w,u][w,v]=[[u,v]vu,w][w,v]=[vu[u,v],w][w,v]$$
is equal to
$$[u,[v,w]]=[uv,w][w,u][w,v]=[[u,v]vu,w][w,v]=[vu[u,v],w][w,v]=$$
$$=[vu,[[u,v],w]][[u,v],w][vu,w][w,v]=[[u,v],w]$$
while
$$s^{a}_{ij}:=s(x^{a}_{ij})=[\phi ^{-1}x^{a}_{in},\phi ^{-1}x^{1}_{nj}].$$
One can rewrite Lemma 2 as:
$$[\phi ^{-1}x^{a}_{hi},\phi ^{-1}x^{bc}_{ik}]=[\phi ^{-1}x^{ab}_{hj},\phi ^{-1}x^{c}_{jk}]$$
and
$$[\phi ^{-1}x^{a}_{hi},\phi ^{-1}x^{b}_{ik}]=[\phi ^{-1}x^{a}_{hj},\phi ^{-1}x^{1}_{jk}],$$
namely it is independent of $n$ $\square$.\\

{\bf Theorem 2.} Every perfect group has a "universal" central  extension:\\

$$\begin{CD}
1 @>>> H_{2}(G) @>>> \tilde{G} @>>> G @>>> 1\\
  @.   @VVV     @VVV           @|     @.\\
1 @>>> C @>>> Y @>>>         G @>>>   1
\end{CD}$$\\
\\
where "universal" means that there exists a unique map $H_{2}(G)\rightarrow C$.\\

The preceeding result shows that the group $StA$ is the universal central extension of the group $EA$. In particular
$$H_{2}(EA)=Ker(StA\rightarrow EA):=K_{2}A$$
whereas
$$K_{1}A=GLA/[GLA,GLA]=H_{1}(GLA,\mathbb{Z}),$$
$$K_{2}A=H_{2}(EA,\mathbb{Z})$$
and
$$K_{3}A=H_{3}(StA, \mathbb{Z}).$$

\newpage

\section{Lecture 4 (Classifying Spaces and Group (Co)Homology)}

Let $G$ be a group and $M$ a $G$-module. We know that for cohomology one has
$$H^{i}(G,M)=H^{i}(C^{.}(G,M))$$
and similar things hold for homology. If $C$ is an abelian group with trivial $G$-action we have
$$H^{0}(G,C)=C,$$
$$H^{1}(G,C)=Hom(G,C)=Hom(G_{(ab)},C)$$
(where $G_{(ab)}=G/[G,G]$) and\\
$H^{2}(G,C)= \{$set of isomorphism classes of central extensions of $G$ by $C\}$, namely

$$\begin{CD}
1 @>>> C @>>> E @>>> G @>>> 1\\
  @.   @|     @VVV           @|     @.\\
1 @>>> C @>>> E' @>>>         G @>>>   1
\end{CD}$$\\
\\
where $C$ is a subset of the centre of $E$.\\

The topological interpretation of group (co)homology is the following: $G$ has a classifying space $BG$ which is a pointed nice space such that:\\
{\bf 1.} $\pi _1 (BG)=G$.\\
{\bf 2.} The universal covering space of $BG$ is contractible.\\

For example if $G=\mathbb{Z}$, then $BG=S^{1}$.\\

$BG$ is unique up to homotopy equivalence.\\

{\bf Fact:} $H^{i}(G,C)=H^{i}(BG,C)$. This is also true for homology.\\

From the \textsl{universal coefficients' theorem} (see appendix), it follows that
$$0\rightarrow Ext^{1}_{\mathbb{Z}}(H_{i-1}(BG,{\mathbb{Z}}),C)\rightarrow H^{i}(BG,C)\rightarrow Hom(H_{i}(BG,{\mathbb{Z}}),C)
\rightarrow 0,$$
where $H_{i-1}(BG,{\mathbb{Z}})=H_{i-1}(G,{\mathbb{Z}})=H_{i-1}G$.\\

Moreover
$$H^{1}(G,C)\simeq Hom(H_{1}G,C)$$
from which we deduce that
$$H_{1}G=G_{(ab)}.$$
Moreover
$$0\rightarrow Ext^{1}_{\mathbb{Z}}(H_{1}G,C)\rightarrow H^{2}(G,C)\rightarrow Hom(H_{2}G,C)\rightarrow 0.$$
In particular, if $G$ is a perfect group, namely $G_{(ab)}=0$, then
$$H^{2}(G,C)=Hom(H_{2}G,C).$$

The central extension which corresponds to the identity homomorphism $id\in Hom(H_{2}G,H_{1}G)$ is the following

$$\begin{CD}
1 @>>> H_2 G @>i>> \tilde{G} @>>> G @>>> 1\\
  @.   @VuVV     @VVV           @|     @.\\
1 @>>> C @>>> u_{*}\tilde{G} @>>>         G @>>>   1
\end{CD}$$\\
\\

where $u:H_{2}G\rightarrow C$.\\

Given any homomorphism $u$, one has a \textsl{push-out} which gives a central extension of $G$ by $C$, where
$$u_{*}\tilde{G}=\frac{C\times\tilde{G}}{\{(-u(x),i(x))|x\in H_{2}G\}}$$
and $i:H_{2}G\hookrightarrow\tilde{G}$.\\

Thus any central extension of $G$ by $C$ is induced by a unique homomorphism $u:H_{2}G\rightarrow C$.\\

{\bf Proposition 1.} The group $\tilde{G}$ is perfect.\\

{\bf Proof:} One has the following diagram

$$\begin{CD}
1 @>>> H_2 G @>>> \tilde{G} @>>> G @>>> 1\\
  @.   @VuVV     @VVV           @|     @.\\
1 @>>> B @>i>> [\tilde{G},\tilde{G}] @>>>         G @>>>   1\\
  @.   @ViVV     @ViVV           @|     @.\\
1 @>>> H_2 G @>i>> \tilde{G} @>>> G @>>> 1\\
\end{CD}$$\\
\\

where $B=H_{2}G\cap [\tilde{G},\tilde{G}]$ and $i$ denotes the inclusion. Then since $u:H_{2}G\rightarrow B$, we shall obtain $i\circ u=1_{H_{2}G}$ since 
$i_{*}\circ u_{*}(G)=i_{*}[\tilde{G},\tilde{G}]=\tilde{G}$. Hence 
$B=H_{2}G\Rightarrow [\tilde{G},\tilde{G}]=\tilde{G}$. $\square$\\

{\bf Proposition 2.} One has that 
$$H_{2}\tilde{G}=0,$$ 
namely $\tilde{G}$ has no non-trivial central extension, in other words every central extension splits.\\

{\bf Proof:}  Given the sequence
$$E\xrightarrow{q} \tilde{G}\xrightarrow{p} G$$
where $E$ is a perfect central extension of $\tilde{G}$ and both $p$ and $q$ are surjective, we claim that $E$ is a central extension $G$.\\

Then
$$1\rightarrow Kerq\rightarrow Ker(pq)\rightarrow Kerp\rightarrow 1.$$
$E$ acts on the above exact sequence of abelian groups, it has trivial action on $Kerq$ and $Kerp$, hence we get a homomorphism $E\rightarrow Hom(Kerp, Kerq)$. Thus the action of $E$ on $Ker(pq)$ is also trivial.\\

So if $E=\tilde{\tilde{G}}$, then it is a perfect central extension of $\tilde{G}$. One has

$$\begin{CD}
1 @>>> H_2 \tilde{G} @>>> \tilde{\tilde{G}} @>>> \tilde{G} @>>> 1\\
  @.   @.     @AAA           @VVV     @.\\
@. @. \tilde{G} @. G @.          @.   @.
\end{CD}$$\\
\\

and the universal property of $\tilde{G}$ implies  that $\tilde{G}$ lifts into $\tilde{\tilde{G}}$. Hence 
$\tilde{\tilde{G}}=\tilde{G}\times H_{2}\tilde{G}\Rightarrow H_{2}\tilde{G}=0$. $\square$\\

{\bf Examples:}\\
{\bf 1.} Let $G=A_5$, the simple non-Abelian group of order 60 (which gives the rotational symmetries of the 20-hedron). What is $\tilde{G}$ and $H_{2}A_{5}$?\\

We start by considering the exact sequence
$$1\rightarrow \{\pm 1\}\rightarrow SU_{2}\rightarrow SO_{3}\rightarrow 1.$$
Recall that topologically $SU_{2}=S^3$ and this is a subgroup of the group of automorphisms of ${\mathbb{C}}^2$ whereas the group $SO_{3}$ is a subgroup of the group of automorphisms of ${\mathbb{C}}P^{1}$ and topologically ${\mathbb{C}}P^{1}=S^2$. Moreover $A_5$ is a subgroup of $SO_3$. We have then the following diagrams:

$$\begin{CD}
1 @>>> {\mathbb{Z}}_{2} @>>> S^3 @>>> SO_{3} @>>> 1\\
  @.   @|     @.           @.     @|\\
1 @>>> {\mathbb{Z}}_{2} @>>> E @>>>         G @>>>   1
\end{CD}$$\\
\\

We know that $E\subset S^3$ and $G\subset SO_3$. We define
$$M:=S^3/E=SO_3/G$$
and $M$ is an orientable Poincare Homology 3-sphere.\\

Because
$$\pi _{i}S^3=0$$
for $i=0,1,2$, then $M$ is "close to" being $BE$. For $i=1,2$, we have that
$$\pi _{i}(M)=\pi _{i}S^3=0$$
from which we deduce that
$$H_{i}M=H_{i}E.$$ 
$M$ is an orientable 3-manifold, hence by Poincre duality we have that $H_2M=0$ since $H_1M=H_1E=0$. Thus 
$H_2E=0$. Finally
$$H_1E=H_2E=0\Rightarrow E=A_5$$ 
and
$$H_1E=H_2E=0\Rightarrow H_{2}A_{5}={\mathbb{Z}}_2.$$ 
$\square$\\

{\bf 2.} We shall mention some analogies between K-Theory and homotopy groups of Lie groups. We know that
$$S_3\twoheadrightarrow SO_3\subset O_3$$

where $\pi _{0}(O_3)=SO_3$ the connected component of the identity of  $O_3$ and $S^3$ is the universal covering of $SO_3$.\\

If we replace the word "connected" by "perfect" we have the following analogy:
$$StA\rightarrow EA\subset GLA$$
where $EA$ is the largest perfect subgroup of $GLA$ whereas $StA=\tilde{EA}$ the universal covering.\\

The infinite orthogonal group $O_\infty =O$ is even closer to K-Theory:
$$Spin\twoheadrightarrow SO\subset O$$
and $\pi _0 (O)=K^{-1}O(pt)$ is the kernel of the inclusion $SO\hookrightarrow O$ whereas $\pi _{1}(O)=K^{-2}O(pt)$ 
is the kernel of the surjection  $Spin\twoheadrightarrow SO$.\\

\newpage

\section{Lecture 5 (The Plus Construction and the general K-group definition)}

By the word \emph{"space"} we mean a $CW$-complex with a base point and we denote by $[X,Y]$ the set of homotopy 
classes of base point preserving maps $X\rightarrow Y$.\\

Let $X$ be a space such that its fundamental group $\pi _1 (X)$ is perfect, namely 
$\pi _{1}(X)_{ab}=0=H_1(X;\mathbb{Z})$.\\

{\bf Problem:} We would like to construct a space $X^{+}$ such that $\pi _{1}(X^{+})=0$ along with a map 
$i:X\rightarrow X^{+}$ so that the induced map $i_{*}:H_{*}(X)\rightarrow H_{*}(X^{+})$ is an \emph{isomorphism} 
in homology.\\

{\bf Step 1.} We choose elements $\gamma _i\in\pi _{1}(X), i\in I$, such that 
the normal subgroup they generate is the whole group $\pi _{1}(X)$.\\

Let's assume the following special case: one $\gamma$ is realised by a loop $S^1\xrightarrow{u} X$ and let
$$Y=X\cup D^2$$
where $D^2$ is a 2-disc (a $2-cell$) whose boundary is the loop $\gamma$.\\

From \textsl{Van Kampen's theorem} (see Appendix) one has that
$$\pi _{1}(Y)=\pi _{1}(X)*\pi _{1}(D^2)=\pi _{1}(X)/\Gamma =0$$
since $\pi _{1}(D^2)=0$ and where $\Gamma$ denotes the normal subgroup generated by $\gamma$. Then:
$$H_2 (X)\rightarrow H_2 (Y)\rightarrow H_2 (X,Y)\rightarrow H_1 (X)\rightarrow H_1 (Y)\rightarrow 0.$$  
But $H_2 (X,Y)=\mathbb{Z}$ whereas $H_{1}(X)=0$. Thus
$$H_n (X)\simeq H_n (Y), n\geq 3.$$
One has the sequence
$$0\rightarrow H_2 (X)\rightarrow H_2 (Y)\xrightarrow{p}{\mathbb{Z}}\rightarrow 0.$$
We apply \textsl{Hurewicz theorem} (see Appendix) and get: $\pi _1 (Y)=0\Rightarrow \pi _2 (Y)=H_2 (Y)$. There is a map 
$$v:S^2\rightarrow Y$$
such that
$$H_2 (S^2)\rightarrow H_2 (Y)\xrightarrow{p}\mathbb{Z}.$$
Put
$$X^{+}:=Y\cup _{v}C^3$$
where $C^3$ is a $3-cell$. Hence we get the following sequence:
$$H_3 (X^{+})\rightarrow H_3 (X^{+},Y)\xrightarrow{\delta} H_2 (Y)\rightarrow H_2 (X^{+})\rightarrow H_2 (X^{+},Y)$$
where
$$\delta :H_3 (X^{+},Y)\rightarrow H_2 (Y)$$
and the image of $\delta$ is the class of $v$. Hence
$$H_2 (X)\rightarrow H_2 (X^{+})$$
and
$$H_n (X)=H_n (Y)=H_n (X^{+}), n\geq 3.$$

{\bf Proposition 1.} Assuming that $\pi _1 (X)$ is perfect, let $i:X\rightarrow X^{+}$ be such that 
$\pi _1 (X^{+})=0$ and the induced map $i_{*}:H_{*}(X)\rightarrow H_{*}(X^{+})$ is an isomorphism. Then for any $Y$ with  
$\pi _1 (Y)=0$ we have that
$$i^{*}:[X^{+},Y]\rightarrow [X,Y]$$
is also an isomorphism.\\

{\bf Corollary 1.} $(X^{+},i)$ are defined only up to homotopy equivalence.\\

{\bf Proof of Proposition 1:} To prove surjectivity, consider the following diagram:

$$\begin{CD}
X @>i>> X^{+} \\
  @VuVV   @VVu'V\\
Y @>>i'> Y\cup _{X}X^{+}=Z\\
\end{CD}$$\\
\\

If $i$ is a homology isomorphism, then so is $i'$.\\

From\textsl{Van Kampen's Theorem} we have that $\pi _{1}(Z)=0$.\\
By the \textsl{Whitehead Theorem} (the 1-connected version, see Appendix) we get that 
$i':Y\rightarrow Z$ is a homotopy equivalence. Now let $\tau :Z\rightarrow Y$ be a homotopy inverse for $i'$.
Then
$$(\tau u')i=\tau (i'u)=(\tau i')u\sim (id)u=u.$$
That proves that $i^{*}:[X^{+},Y]\rightarrow [X,Y]$ is onto.\\

To prove injectivity, consider the following diagram:

$$\begin{CD}
X @>i>> X^{+}\\
  @VuVV   @VVu'V\\
Y @>>i'> Z
\end{CD}$$\\
\\
 
Next consider two maps  
$$g_{0},g_{1}:X^{+}\rightarrow Y$$
such that $g_0 i$ and $g_1 i$ are both homotopic to $u$. Then we apply the homotopy extension theorem and we deform 
$g_0$ and $g_1$ so that $g_0 i=g_1 i=u$. Then we define
$$\tau _0:Z\xrightarrow{(id_Y,g_0)}Y$$
and
$$\tau _1:Z\xrightarrow{(id_Y,g_1)}Y$$
so that
$$\tau _0 i'=\tau _1 i'=id_Y.$$
But $i'$ is a homotopy equivalence, hence $\tau _0 \sim\tau_1$. Since $\tau _0 u'=g_0$ and $\tau _1 u'=g_1$, we deduce that
$g_0\sim g_1$. 
 
$\square$.\\

Suppose that $G$ is perfect and let
$$1\rightarrow H_{2}G\rightarrow\tilde{G}\rightarrow G\rightarrow 1$$
be its universal central extension. Recall that $H_{1}G=H_{2}\tilde{G}=0$. Then we have the following:\\

{\bf Lemma 1.} (Put $C=H_{2}G$). One has a map of fibrations:   

$$\begin{CD}
BC @>>> B\tilde{G} @>>> BG\\
  @|   @VVV     @VVV\\
BC @>>> B\tilde{G}^{+} @>>> BG^{+}
\end{CD}$$\\
\\

{\bf Proof:} (Reacall that $\pi _{1}(BG)=G$). The proof is by diagram chasing:
$$H^{2}(G,C)=H^{2}(BG,C)=[BG,K^{EM}(C,2)]$$
(where $K^{EM}$ is the Eilenberg-McLane space). Moreover

$$\begin{CD}
BC @= BC @= K^{EM}(C,1)\\
  @VVV   @VVV     @VVV\\
B\tilde{G} @>\beta >> P @>>> E(C,2)\\
  @VVV   @VVV     @VVV\\
BG @>>> BG^{+} @>>\alpha > K^{EM}(C,2)\\
\end{CD}$$\\
\\

where $\alpha :BG^{+}\rightarrow K^{EM}(C,2)$ exists by the previous proposition 1 and $\beta :B\tilde{G}\rightarrow P$.
We also have:
$$\pi _{2}(BG^{+})=H_{2}(BG^{+})=H_{2}(BG)=C.$$
But $\alpha$ induces an isomorphism on the $\pi _2$'s:
$$\pi _{2}(BG^{+})\rightarrow\pi _{1}(BC)\rightarrow\pi _{1}(P)\rightarrow\pi _{1}(BG^{+})$$
where $\pi _{1}(BC)=C$ and $\pi _{1}(BG^{+})=0$. Hence  $\pi _{1}(P)=0$.\\

Now $\beta$ is a homology isomorphism, $\pi _{1}(P)=0\Rightarrow B\tilde{G}^{+}\sim P$ are homotopy equivalent.\\

From the maps
$$BC\rightarrow B\tilde{G}^{+}\rightarrow BG^{+},$$
using Hurewicz theorem, we find that
$$\pi _{3}(BG^{+})=\pi _{3}(B\tilde{G}^{+})=H_{3}(B\tilde{G}^{+})$$
since $H_{1}\tilde{G}=H_{2}\tilde{G}=0$.\\

Next assume that $G=EA$ is perfect. We have that $\tilde{G}=StA$ and $H_{2}G=K_{2}A$.\\

{\bf Claim:}
$$\pi _{2}(BEA^{+})=K_{2}A=H_{2}(EA)$$
$$\pi _{3}(BEA^{+})=H_{3}(StA).$$

The proof of the claim follows from what we mentioned above:
$$\pi _{2}(BEA^{+})=H_{2}(BEA^{+})=H_{2}(BEA)=H_{2}(EA)=K_{2}A.$$
$$\pi _{3}(BEA^{+})=\pi _{3}(BStA^{+})=H_{3}(BStA)=H_{3}(StA).$$
$\square$.\\

{\bf Theorem 1.} Let $N\subset\pi _{1}(X)$ be a perfect normal subgroup. Then there exists a space $X^{+}$ (depending on 
$N$) along with a map $i:X\rightarrow X^{+}$ such that:\\

$\bullet$ The groups $\pi _{1}(X)/N\simeq \pi _{1}(X^{+})$ are isomorphic.\\

$\bullet$ For any $\pi _{1}(X^{+})-mod$ $L$, one has that the induced map in homology
$$i_{*}:H_{*}(X,i^{*}L)\rightarrow H_{*}(X,L)$$
is an isomorphism.\\

The pair $(X^{+},i)$ is only defined up to homotopy equivalence.\\

{\bf Proof:} Essentially we shall construct the space $X^{+}$ with the desired properties, this is \emph{Quillen's 
famous plus construction}.\\

Let $\tilde{X}$ be the covering space corresponding to the subgroup $N\subset\pi _{1}(X)$, namely 
$\pi _{1}(\tilde{X})=N$. Then we have: 

$$\begin{CD}
\tilde{X} @>>> \tilde{X}^{+}\\
 @VVV     @VVV\\
X @>>> X^{+}
\end{CD}$$\\
\\

Then $\pi _{1}(X^{+})=\pi _{1}(X)/N$.\\

For K-Theory, we make the following choices: $X=BGLA$, $\pi _{1}(X)=GLA$, $EA=N$ and obviously 
 $EA\subset GLA$. We thus get 
 $BGLA^{+}$ with $\pi _{1}(BGLA^{+})=GLA/EA=K_{1}A$.\\

We also have 
$$H_{*}(BGLA,L)\simeq H_{*}(BGLA^{+},L)$$
for all modules $L$ over $K_{1}A$ along with the following diagram 

$$\begin{CD}
BEA @>>> \tilde{BGLA}^{+}\\
 @VVV     @VVV\\
BGLA @>>> BGLA^{+}
\end{CD}$$\\
\\

where $\tilde{BGLA^{+}}$ is the universal covering and the vertical arrows $BEA\rightarrow BGLA$ and 
$\tilde{BGLA^{+}}\rightarrow BGLA^{+}$ are pull backs. Hence $BEA^{+}$ is the universal covering of 
$BGLA^{+}$, thus 
$$\pi _{n}(BGLA^{+})=\pi _{n}(BEA^{+}), n\geq 2.$$
For $n=2$ one has $\pi _{2}(BGLA^{+})=K_{2}A$.\\

For $n=3$ one has $\pi _{3}(BGLA^{+})=H_{3}(StA)$.\\

\emph{Hence we have the general K-Group definition given by Quillen:}\\

{\bf Definition 1:}
$$K_{n}A:=\pi _{n}(BGLA^{+}), n\geq 1.$$

\newpage

\section{Lecture 6 (Some Theorems and Examples)}

We give the following proposition without proof which we shall use below to get a basic definition:\\

{\bf Proposition 1.} For a map  $f:X\rightarrow Y$ (where $X,Y$ are $CW$-complexes with base points), 
the following are equivalent:\\

{\bf (a.)} The homotopy fibre $F$ of $f$ is acyclic, namely $\tilde{H}_{*}F=0$, where homotopy fibre means replace 
 $f$ by a fibration and take the actual fibre.\\

{\bf (b.)} $\pi _{1}(f):\pi _{1}(X)\rightarrow\pi _{1}(Y)$ is surjective and for any $\pi _{1}(Y)-mod$ $L$ we have 
that the induced map in homology
$$f_{*}:H_{*}(X,f^{*}L)\rightarrow H_{*}(Y,L)$$
is an isomorphism.\\

{\bf (c.)} The lift $f':X\times _{Y}\tilde{Y}\rightarrow\tilde{Y}$ of $f:X\rightarrow Y$ is a homology isomorphism
(where $\tilde{Y}$ is the universal cover of $Y$):\\ 

$$\begin{CD}
X\times _{Y}\tilde{Y} @>f'>> \tilde{Y}\\
 @VVV     @VVV\\
X @>>f> Y
\end{CD}$$\\
\\

{\bf Definition 1.} We call $f$ \emph{acyclic} when the above conditions hold.\\

{\bf Corollary 1.} Acyclic maps are closed under composition, homotopy pull-backs and homotopy push-outs.\\

{\bf Theorem 1.} Given a perfect normal subgroup $N\subset\pi _{1}(X)$, there is a unique (up to homotopy equivalence) 
acyclic map $f:X\rightarrow Y$ such that $N=ker[\pi _{1}(f)]$ (where $\pi _{1}(Y)=\pi _{1}(X)/N$). 

Moreover, for any space $T$, the induced map
$$f^{*}:[Y,T]\rightarrow \{a\in [X,T] | \pi _{1}(a):\pi _{1}(X)\rightarrow\pi _{1}(T) kills N\}$$
is an isomorphism where

$$\begin{CD}
X @>f>> Y\\
 @VVV     @|\\
T @<<< Y
\end{CD}$$\\
\\

{\bf Proof:} The key point is the construction of $Y$. Let $\tilde{X}$ be a covering space of $X$ with 
$\pi _{1}(\tilde{X})=N$. Since $f:X\rightarrow Y$ is acyclic, we use the acyclic push out
$\tilde{X}\rightarrow\tilde{X}^{+}$, namely

$$\begin{CD}
\tilde{X} @>>> \tilde{X}^{+}\\
 @VVV     @VVV\\
X @>>> Y
\end{CD}$$\\
\\

$\square$.\\

If we apply the above to $EA$ which is a perfect normal subgroup of $GLA=\pi _{1}(BGLA)$, we get:

$$\begin{CD}
BEA @>\tilde{f}>> \tilde{BGLA}^{+}\\
 @VVV     @VVV\\
BGLA @>>f> BGLA^{+} 
\end{CD}$$\\
\\

where $f:BGLA\rightarrow BGLA^{+}$ is the unique acyclic map with $ker[\pi_{1}(f)]=EA$, $\tilde{BGLA^{+}}$ is the 
universal cover and
 $\pi _{1}(BGLA^{+})=GLA/EA=K_{1}A$. If $f$ is acyclic, then $\tilde{f}: BEA\rightarrow\tilde{BGLA^{+}}$ is also acyclic.\\

From all the above we conclude that:\\

$\bullet$ $BGLA^{+}$ has universal cover $BEA^{+}$.\\

$\bullet$ $\pi _{1}(BGLA^{+})=K_{1}A$.\\

$\bullet$ $\pi _{2}(BGLA^{+})=\pi _{2}(BEA^{+})=H_{2}(BEA^{+})=H_{2}(EA)=K_{2}A$.\\

$\bullet$ $\pi _{3}(BGLA^{+})=H_{3}(StA)$.\\

$\bullet$ In general, $\pi _{n}(BGLA^{+})=K_{n}A, n\geq 1$.\\

Next we shalll see the relation with topological K-Theory (see Chapter 3).\\

Let
$$G_{\infty}=lim _{k\rightarrow\infty}Gr_{k}({\mathbb{C}}^{\infty})=BGL{\mathbb{C}}$$
where $Gr$ denotes the $Grassmannians$ and $Gr_{k}({\mathbb{C}}^{\infty})=BGL_{k}{\mathbb{C}}$. We assume that the group 
$GL\mathbb{C}$ has its natural topology.\\

To the direct sum of vector bundles $\xi\oplus\eta$ corresponds an $h-space$ structure on $G_{\infty}$. Note that 
 although $BGLA$ is not an $h-space$, $BGLA^{+}$ is one. There exists a continuous map 
$\mu :G_{\infty}\times G_{\infty}\rightarrow G_{\infty}$ (with $x_0$ the base point) such that $\mu (x_0,x)=x$ and 
$\mu (x,x_0)=x_0$.\\

Recall the following fact: Let $a\in GL_{n}\mathbb{C}$. The maps
$$a\mapsto   
\left( \begin{array}{cc}a & 0\\
0 & 1\end{array} \right)$$
and
 $$a\mapsto   
\left( \begin{array}{cc}1 & 0\\
0 & a\end{array} \right)$$
map the group $GL_{n}\mathbb{C}$ to $GL_{2n}\mathbb{C}$ and
they are homotopic in the natural topology but the maps
$$a\mapsto   
\left( \begin{array}{ccc}a & 0 & 0\\
0 & 1 & 0\\
0 & 0 & 1\end{array} \right)$$
and
 $$a\mapsto   
\left( \begin{array}{ccc}1 & 0 & 0\\
0 & a & 0\\
0 & 0 & 1\end{array} \right)$$
are $conjugate$ via an element of $EA$. Conjugation by elements of $EA$ on $BGLA$ is not trivial (provided base point 
preserving maps are considered) but since 
$$\pi _{1}(BGLA)=GLA,$$
conjugation by elements of $EA$ on $BGLA^{+}$ is trivial up to homotopy.\\

We know that homotopic spaces have isomorphic vector bundles in Topological K-Theory but it $fails$ in Algebraic K-Theory.\\

In general, Higher Algebraic K-Theory computations are hard. It's been over 30 years since these groups were defined by 
Quillen and only few examples are known. We shall present  some: Let $F[x_1,x_2,...,x_n]\sim F$ 
be the ring of polynomials with coefficients from some field $F$. 
Grothendieck proved that 
$$K_{0}(F[x_1,x_2,...,x_n])=K_{0}F,$$
whereas Quillen proved that every vector bundle over a polynomial ring is trivial.\\

Let ${\mathbb{F}}_q$ be a finite field of $q$ elements. The following groups have been computed for this case (by Quillen):\\
$$K_{n}({\mathbb{F}}_q)={\mathbb{Z}}, n=0.$$
$$K_{n}({\mathbb{F}}_q)={\mathbb{F}}_{q}^{+}\simeq{\mathbb{Z}}_{q-1}, n=1.$$
$$K_{n}({\mathbb{F}}_q)=0, n=2.$$
$$K_{n}({\mathbb{F}}_q)={\mathbb{Z}}_{q^2 -1}, n=3.$$
$$K_{n}({\mathbb{F}}_q)=0, n=4.$$
$$K_{n}({\mathbb{F}}_q)={\mathbb{Z}}_{q^3 -1}, n=5.$$
$$etc$$

In a more concise form one has:\\

$$K_{0}({\mathbb{F}}_q)={\mathbb{Z}},$$
$$K_{2n}({\mathbb{F}}_q)=0, n\in{\mathbb{N}^{*}}$$
and
$$K_{2n-1}({\mathbb{F}}_q)={\mathbb{Z}}_{q^n -1},   n\in{\mathbb{N}^{*}},$$
where ${\mathbb{Z}}_{r}$ denotes the cyclic group with $r$ elements.\\

We have the following maps:
$$BGL{\mathbb{F}}_{q}\rightarrow BU\rightarrow BU,$$
where the first map is the Brauer lift, $BU=G_{\infty}$ in topology and the second map is 
$\Psi ^{q}(L)=L^{q}$. There is a map from $BGL{\mathbb{F}}_{q}^{+}$ to the $h-fibre$ of 
$BU\xrightarrow{\Psi ^{q-1}}BU$.\\

We have used the following theorem:\\

{\bf Theorem 2.} The map $BGL{\mathbb{F}}_{q}^{+}\rightarrow h-fibre (BU\xrightarrow{\Psi ^{q-1}} BU)$ is a homotopy 
equivalence.\\

Because $BGLA^{+}$ is an $h$-space, one knows that $H_{*}(BGLA^{+};\mathbb{Q})$ is a Hopf algebra. 
Moreover one also knows that
$$\pi _{*}(BGLA^{+})\otimes{\mathbb{Q}}=Prim[H_{*}(GLA;{\mathbb{Q}})]$$
which is the \textsl{Milnor-Moore Theorem}.\\

For the case
$A={\mathbb{Z}}$, Borel proved that:
$$\dim _{\mathbb{Q}}[K_{n}{\mathbb{Z}}\otimes{\mathbb{Q}}]={\mathbb{Z}},n=0,5,9.$$
$$\dim _{\mathbb{Q}}[K_{n}{\mathbb{Z}}\otimes{\mathbb{Q}}]=0,n=1,2,3,4,6,7,8.$$

More concretely:\\

{\bf Theorem 3} (D.G. Quillen) If $A$ is the ring of algebraic integers in an algebraic number field ${\mathbb{F}}$ (a finite extension of the rationals), then  the groups $K_n(A)$ are finitely generated.\\

Borel used this to compute both $K_n(A)$ and $K_n({\mathbb{F}})$ modulo torsion. For example, for the integers Borel proved that (modulo torsion)
$$K_n({\mathbb{Z}})=0$$
for positive $n$ unless $n=4k+1$ with $k$ again positive and (modulo torsion)
$$K_{4k+1}({\mathbb{Z}})={\mathbb{Z}}$$
for positive $k$. The torsion subgroups of $K_{2n+1}({\mathbb{Z}})$ and the orders of the finite groups $K_{4k+2}({\mathbb{Z}})$ have recently been determined but whether the later groups are cyclic and whether the groups $K_{4k}({\mathbb{Z}})$ vanish depend on the so-called Vandiver's conjecture about the class groups of cyclotomic integers.\\

\newpage

\section{Appendix}

In this Appendix we briefly recall some definitions and results from other branches of mathematics to  help the reader in his study (for more details one can see the list of books with background meterial in the introduction).\\

We start with the \emph{Seifert-van Kampen theorem}, which expresses the fundamental group of some topological space say $X$ in terms of  the fundamental groups of two open and path connected topological subspaces  $X_1, X_2\subset X$ which cover $X$ (hence it can be used to compute the fundamental groups of spaces which are constructed out of simpler ones):\\

{\bf Van Kampen's Theorem}. Let $X$ be a topological space which is the union of the interiors of two path connected subspaces $X_1,X_2\subset X$. Suppose that $X_0=X_1\cap X_2$ is non empty and path connected. 
Let also $*$ be a base point of $X_0$ and suppose that
$i_k :\pi _1 (X_0 ,*)\rightarrow\pi _1 (X_k, *)$ and $j_k :\pi _1 (X_k ,*)\rightarrow\pi _1 (X, *)$ are the induced maps from the corresponding inclusions for $k=1,2$. Then $X$ is path connected and the natural map  
$$\pi _1 (X_1,*)\star _{\pi _1 (X_0, *)}\pi _1 (X_2,*)\rightarrow\pi _1 (X,*)$$
is an isomorphism, namely the fundamental group of $X$ is the free product of the fundamental groups of the subspaces $X_1$ and $X_2$ with the amalgamation of the fundamental group of $X_0$.\\

The key idea of the proof is that paths in $X$ can be analysed in parts of paths inside the intersection
$X_0$, inside $X_1$ but outside $X_2$ and inside $X_2$ but outside of  $X_1$. Usually the induced maps of inclusions are not injective and hence more precise versions of the theorem use the push-outs of groups. 
There is a generalisation for non-connected spaces in the category of groupoids (as is the case in noncommutative geometry) along with generalisations in "higher dimensions" (for example using 2-groups etc, for more details one can see J.P. May's notes on algebraic topology).\\

Next we recall the so called \emph{5-Lemma} in homological algebra. In fact this is true in any abelian category (like the category of abelian groups or the category of vector spaces over a given field):\\

{\bf 5-Lemma}. Suppose that in an abelian category the following commutative diagram (lader) is given:

$$\begin{CD}
A @>f>> B @>g>> C @>h>> D @>j>> E\\
  @VlVV   @VmVV   @VnVV  @VpVV  @VVqV\\
A' @>>r> B' @>>s> C' @>>t> D' @>>u> E'
\end{CD}.$$\\
\\

Then if the rows are exact sequences, $m$ and $p$ are isomorphisms, $l$ is an epimorphism and $q$ is a monomorphism, then $n$ is also an isomorphism.\\

The proof is by "diagram chasing".\\

The 5-Lemma is the combination of the two 4-Lemmas where one is the dual of the other. An interesting special case is the so-called \emph{short 5-Lemma} where the rows are short exact sequences; in this case one has that $A=A'=E=E'=0$, hence one has only three vertical maps $m,n,p$  and the short 5-Lemma states that if $m,p$ are isomorphisms, then so is $n$.\\

A closely related statement in homological algebra is the \emph{serpent (snake) lemma}:\\

{\bf Serpent Lemma}. In an abelian category (e.g. the category of vector spaces over some given field), consider the following commutative diagram:

$$\begin{CD}
0 @>>> A @>f>> B @>g>> C @>>> 0\\
  @.   @VaVV   @VbVV  @VVcV  @.\\
0 @>>> A' @>>f'> B' @>>g'> C' @>>> 0
\end{CD}$$\\
\\
where the rows are exact and $0$ is the zero object. Then there is an exact sequence relating the kernels and cokernels of the maps $a,b$ and $c$:
$$0\rightarrow Ker a\rightarrow Ker b\rightarrow Ker c\xrightarrow{d} coker a\rightarrow coker b\rightarrow coker c\rightarrow 0.$$

Furthermore if $f$ is a monomorphism, then so is the map $Ker a\rightarrow Ker b$ and if $g'$ is an epimorphism, then so is the map $coker b\rightarrow coker c$.\\

Next we recall the \emph{universal coefficient theorem} which gives a relation between the integral homology $H_i(X; {\mathbb{Z}})$ of some topological space $X$ and the corresponding homology
$H_i(X; A)$ with coefficients from an arbitrary abelian group $A$:\\

{\bf Universal Coefficient Theorem}.  With the above notation, consider the tensor product $H_i(X; {\mathbb{Z}})\otimes A$. 
Then there is a group homomorphism $i:H_i(X; {\mathbb{Z}})\otimes A\rightarrow H_i (X; A)$ which is injective and whose cokernel is $Tor(H_{i-1}(X;A),A)$.\\

In other words, there is a natural short exact sequence 
$$0\rightarrow H_i(X; {\mathbb{Z}})\otimes A\rightarrow H_i(X; A)\rightarrow Tor(H_{i-1}(X;A),A)\rightarrow 0.$$    
This sequence splits (but not in a natural way) and the torsion group can be considered as the obstruction to $i$ being an isomorphism. There is a corresponding (dual) theorem for cohomology.\\

The \emph{Hurewicz Theorem} is a central result in algebraic topology which generalises the theorem by Poincare which relates homotopy and homology. One version is this:\\

{\bf Hurewicz Theorem}. For any topological space $X$ and natural number $k$ there exists a group homomorphism 
$$h_* :\pi _k (X)\rightarrow H_k(X; {\mathbb{Z}})$$
which is called the \textsl{Hurewicz homomorphism}. There is a relative version along with a triadic version.\\

For $k=1$ the Hurewicz homomorphism is the canonical abelianisation map 
$$h_* :\pi _1 (X)\rightarrow \pi _1 (X)/[\pi _1 (X), \pi _1 (X)].$$
In particular the first (integral) homology group is isomorphic to the abelianisation of the fundamental group. Moreover if $X$ is $(n-1)$-connected, then the Hurewicz homomorphism is an isomorphism $\forall k\leq n$ (which means that the first homology group vanishes if the space is path-connected and the fundamental group is perfect).\\

The \emph{Whitehead theorem} justifies the use of $CW-complexes$ in algebraic topology since it states that if a continuous map between two topological spaces induces isomorphisms between all the corresponding homotopy groups (in other words it is a \textsl{weak homotopy equivalence}), then the map is in fact a homotopy equivalence if and only if the spaces are connected $CW-complexes$:\\

{\bf Whitehead Theorem}. Let $X,Y$ be two $CW-complexes$ with base points $x,y$ respectively and suppose $f:X\rightarrow Y$ is a continuous map such that $f(x)=y$. For any positive integer $n\neq 0$ we consider the induced homomorphisms
$$f_*:\pi _n (X,x)\rightarrow\pi _n (Y,y).$$ 
We say that $f$ is a weak homotopy equivalence if $f_*$ is an isomorphism for all $n$. Then if $f$ is a weak homotopy equivalence, then it is also a homotopy equivalence.\\

Generalisations of Whitehead's theorem for spaces which are not $CW-complexes$ are studied by another branch of topology called \textsl{shape theory}. Moreover $Quillen$ proved that in any \textsl{model category}, a weak homotopy equivalence between \textsl{fibrant} and \textsl{cofibrant} objects is in fact a homotopy equivalence. \\

\emph{Eilenberg-MacLane} spaces are the building blocks of homotopy theory:\\

{\bf Definition of Eilenberg-MacLane spaces}. Let $G$ be a group and $n$ a non-zero natural number. A connected topological space $X$ is called an \emph{Eilenberg-MacLane space of type} $K(G,n)$ if the $n-th$ homotopy group $\pi _n (X)$ is isomorphic to $G$ and all other homotopy groups of $X$ vanish.\\

For $n>1$, $G$ must be abelian and hence as a $CW-complex$ the Eilenberg-MacLane space exists and it is unique up to weak homotopy.\\

Examples:\\

$\bullet$ $K({\mathbb{Z}} ,1)=S^1$.\\

$\bullet$ $K({\mathbb{Z}} ,2)={\mathbb{C}}P^{\infty}$.\\

$\bullet$ $K({\mathbb{Z}}_{2},1)={\mathbb{R}}P^{\infty}$.\\

$\bullet$ $K(G,1)=\wedge _{i=1}^{k}S^1$, where $G$ is a free group with $k$ generators.\\

Every space $K(G,n)$ can be constructed as a $CW-complex$ starting from the smash product of factors $S^n$, one factor for every generator of $G$ and next by adding cells (possibly infinite in number) of higher dimension to kill of the remaining homotopy.\\

The basic property of the Eilenberg-MacLane spaces is that they give representations of homology with $G$ coefficients: For any abelian group $G$ and $CW-complex$ $X$, the set $[X,K(G,n)]$ of homotopy classes of maps $f:X\rightarrow K(G,n)$ is in a natural 1-1 correspondence with the set $H^n(X;G)$ 
(the n-th singular cohomology group of $X$ with coefficients from $G$).\\

{\bf Milnor-Moore Theorem}. Let $X$ be a 1-connected space with loop space $\Omega X$. Then the Hurewicz homomorphism induces a Hopf algebra isomorphism
$$U(\pi _{*}(\Omega X)\otimes{\mathbb{Q}})\rightarrow H_{*}(\Omega X;{\mathbb{Q}})$$
where $U$ denotes the universal enveloping algebra and the Lie commutator of 
$\pi _{*}(\Omega X)\otimes{\mathbb{Q}}$  
is given by the Samelson product.\\

{\bf Group Completion Theorem} (D. McDuff, G. Segal). Let $M$ be a topological monoid and $BM$ its classifying space. Let $M\rightarrow\Omega BM$ be the canonical map. Then the map
$$H_{*}(M)\rightarrow H_{*}(\Omega BM)$$
induces an isomorphism 
$$H_{*}(M)[\pi _{0}(M)^{-1}]\rightarrow H_{*}(\Omega BM).$$

The \emph{closed graph theorem} is an important result in functional analysis:\\

{\bf Closed Graph Theorem}. Let $T:X\rightarrow Y$ be a linear operator between two Banach spaces which is defined in the whole of $X$. The graph of $T$ is defined as
$$\{(x,y)\in X\times Y:Tx=y\}.$$
Then $T$ is continuous if and only if it is closed, namely if and only if its graph is a closed subset of the Cartesian product $X\times Y$ 
(equipped with the product topology).\\

Recall that the  \textsl{Weierstrass theorem} (1885)  in analysis (real and complex) states that every continuous function which is defined in a closed interval can be uniformly approximated by a polynomial function.\\

In 1937 $M.H. Stone$ generalised the above result in two directions: Instead of a closed interval he assumed an arbitrary compact Hausdorff space $X$ (this can be further generalised for non-compact $Tychonoff$ spaces) and instead of polynomial functions he considered general subalgebras of $C(X)$ 
(the algebra of continuous complex functions on $X$). More specifically we have:\\

{\bf Stone-Weierstrass Theorem}. Let $X$ be a compact Hausdorff space and let $S\subset C(X)$ be a subset which \emph{seperates points}    (namely for any pair of points $x,y\in X$ there exists a function $f\in S$ such that $f(x)\neq f(y)$). Then the complex unital $*$-algebra defined by $S$ is dense in $C(X)$.\\

Note that there is also a real version of this theorem along with a number of other generalisations (involving locally compact Hausdorff spaces etc).\\

We close this section with \textsl{Gelfand's theorem} which is an important motivation for the development of \emph{noncommutative geometry}:\\

{\bf Gelfand Theorem}. \emph{The following two categories are equivalent}:\\
{\bf 1.} The category with objects the unital commutative $C^*$-algebras and arrows the $*$-preserving homomorphisms.\\
{\bf 2.} The category with objects the compact Hausdorff spaces and arrows the homeomorphisms among them.\\ 
The first functor is denoted $C(-)$ and to each compact Hausdorff space $X$ it assigns the unital commutative $C^*$-algebra $C(X)$ of continuous complex functions on $X$; the second functor is denoted  $Spec(-)$ 
and to each commutative unital $C^*$-algebra $A$ it assigns the compact Hausdorff space $Spec(A)$ which is the spectrum of $A$.\\

The above equivalence can be extended between \textsl{locally compact} Hausdorff spaces and \textsl{non-unital} commutative $C^*$-algebras.\\

One of the starting points in the development of noncommutative geometry by the great French mathematician Alain Connes is precisely the \emph{"topological extension"} of this equivalence when one considers \emph{noncommutative} algebras.\\

\end{document}